\documentclass[11pt]{article}

\usepackage[margin=2.5cm]{geometry}
\usepackage[utf8]{inputenc}
\usepackage[T1]{fontenc}
\usepackage[utf8]{inputenc}
\usepackage{lmodern}
\usepackage{textcomp}
\usepackage{microtype}
\usepackage{mathtools}
\usepackage{verbatim}
\usepackage{titlesec}
\usepackage{url}
\usepackage{graphicx,wrapfig,lipsum}
\usepackage{enumitem} \setlist{noitemsep}
\usepackage{leftidx} 
\usepackage{float}
\usepackage{amsmath}
\usepackage{amsfonts}
\usepackage{amssymb}
\usepackage{amsthm}
\usepackage{makecell}
\usepackage{titlesec}
\usepackage{color}
\usepackage{xfrac}
\usepackage{faktor} 
\usepackage[utf8]{inputenc}
\usepackage{xcolor}
\usepackage{stmaryrd}
\usepackage{import}
\usepackage{graphicx}
\usepackage{siunitx}
\usepackage{import}
\usepackage{placeins}
\usepackage{graphicx}
\usepackage[skip=5pt]{caption}
\usepackage{subcaption}
\usepackage{hhline}
\usepackage[cmtip,all]{xy}
\usepackage{makecell}
\usepackage{array,amsmath,booktabs}
\usepackage{calc}
\usepackage{mathtools}
\usepackage{enumitem}

\newcommand{\longsquiggly}{\xymatrix{{}\ar@{~>}[r]&{}}}

\makeatletter
\def\bign#1{\mathclose{\hbox{$\left#1\vbox to8.5\p@{}\right.\n@space$}}\mathopen{}}
\def\Bign#1{\mathclose{\hbox{$\left#1\vbox to11.5\p@{}\right.\n@space$}}\mathopen{}}
\def\biggn#1{\mathclose{\hbox{$\left#1\vbox to14.5\p@{}\right.\n@space$}}\mathopen{}}
\def\Biggn#1{\mathclose{\hbox{$\left#1\vbox to17.5\p@{}\right.\n@space$}}\mathopen{}}
\makeatother

\usepackage[english]{babel}
\usepackage{csquotes}
\definecolor{darkred}{rgb}{0.65, 0.0, 0.0}
\usepackage[colorlinks=true,linkcolor=darkred,citecolor=blue]{hyperref}

\usepackage[numbers]{natbib}

\usepackage{thmtools}

\declaretheorem[
	name=Theorem,
	numberwithin=section
	]{thm}
\declaretheorem[
	name=Lemma,
	sibling=thm,
	]{lem}
\declaretheorem[
	name=Proposition,
	sibling=thm,
	]{prop}
\declaretheorem[
	name=Corollary,
	sibling=thm,
	]{cor}

\declaretheorem[
	name=Definition,
	style=definition,
	sibling=thm
	]{defin}

\declaretheorem[
	name=Remark,
	style=remark,
	sibling=thm
	]{rem}
\declaretheorem[
	name=Claim,
	numberwithin=thm
	]{claim}
\declaretheorem[
	name=Algorithm,
	style=remark,
	sibling=claim
	]{alg}	
\declaretheorem[
	name=Example,
	style=remark,
	sibling=thm
	]{exam}
\declaretheorem[
	name=Construction,
	style=remark,
	sibling=thm
	]{constr}
\declaretheorem[
	name=Acknowledgements,
	style=remark,
	numbered=no
	]{ackn}
\declaretheorem[
	name=Outline,
	style=remark,
	sibling=thm
	]{out}

\declaretheorem[
	name=Standing Assumption,
	style=definition,
	]{Assump}
\declaretheorem[
	name=Convention,
	style=remark,
	numbered=no
	]{nota}
	
\newenvironment{sproof}{%
  \proof}{\endproof}

\newcommand\restr[2]{{
  \left.\kern-\nulldelimiterspace 
  #1 
  \vphantom{\big|} 
  \right|_{#2} 
  }}

\usepackage{chngcntr}
\counterwithin{figure}{section}

\AtBeginDocument{%
   \def\MR#1{}
}

\begin{document}

\title{Quasi-Isometries for certain Right-Angled Coxeter Groups}
\author{Alexandra Edletzberger}
\date{}
\maketitle
\vspace{-0.7cm}
{\centering \small Universit\"at Wien, Fakult\"at f\"ur Mathematik\\
Oskar-Morgenstern-Platz 1, 1090 Wien, Austria \\ \vspace{0.2cm}}
{\centering \small \textit{Email address:}  \texttt{alexandra.edletzberger@univie.ac.at}\\ \vspace{0.2cm}}

\bigskip
\begin{abstract}
We construct the JSJ tree of cylinders $T_c$ for finitely presented, one-ended, two-dimensional right-angled Coxeter groups (RACGs) splitting over two-ended subgroups in terms of the defining graph of the group, generalizing the visual construction by Dani and Thomas \cite{DaniThomasBowditch} given for certain hyperbolic RACGs. Additionally, we prove that $T_c$ has two-ended edge stabilizers if and only if the defining graph does not contain a certain subdivided $K_4$. By use of the structure invariant of $T_c$ introduced by Cashen and Martin \cite{CashenMartinStructureInvar}, we obtain a quasi-isometry-invariant of these RACGs, essentially determined by the defining graph. Furthermore, we refine the structure invariant to make it a complete quasi-isometry-invariant in case the JSJ decomposition of the RACG does not have any rigid vertices.
\end{abstract}

\noindent%
{\it Keywords:} Coxeter groups, Visual decomposition, JSJ splitting, Tree of cylinders, Structure invariant

\section{Introduction}
In this paper, we give a construction of the JSJ tree of cylinders of a wide family of right-angled Coxeter groups (RACGs). It is visual, that is, it is determined in terms of the defining graph:

\begin{thm} \textup{[cf.\ Theorem \ref{SummaryMainThm}]} \label{MainToC}
For a one-ended, two-dimensional RACG W splitting over two-ended subgroups, the defining graph visually determines the JSJ tree of cylinders $T_c$: Subsets of vertices of the defining graph satisfying certain graph theoretic conditions are in bijection with W-orbits of vertices of $T_c$ and they generate the representatives of the conjugacy classes of the vertex stabilizers.
\end{thm}

With this construction, generalizing the one by Dani and Thomas \cite{DaniThomasBowditch} for certain such RACGs which are in addition hyperbolic, the JSJ tree of cylinders can be easily "read off" the defining graph. Throughout this article, we will illustrate the convenience of this method with a range of examples. In particular, the cylinder vertices are produced by a simple process, see Section \ref{SubsectionCylVertices}: Each comes from an uncrossed cut collection, that is a cut pair or a cut triple, of the defining graph and its common adjacent vertices. This implies that cylinder vertices occur only in three types: Two-ended, virtually $\mathbb{Z}^2$ or the direct product of a virtually non-abelian free group and an infinite dihedral group.\\

Additionally, we characterize the edge stabilizers of the JSJ tree of cylinders visually:
\begin{thm} \label{ThmEdgeStabTwoEnded}
All the edge stabilizers of the JSJ tree of cylinders of a one-ended, two-dimensional RACG W splitting over two-ended subgroups are two-ended if and only if in the defining graph no uncrossed cut collection contains opposite corners of a square, whose other two corners are connected by a subdivided diagonal.
\end{thm}

We are interested in the JSJ tree of cylinders mainly because it has two key features:
\begin{enumerate}
    \item It is a canonical representative of the space of all JSJ decompositions of a group.
    \item Quasi-isometric groups have isomorphic JSJ trees of cylinders.
\end{enumerate}

While Feature 1 can be of interest on its own, see \cite[Part IV]{JSJDecompGps} and Section \ref{SubsectionToC}, we aim to use it rather as a stepping stone towards the application of Feature 2 as this provides a new tool to classify a large class of RACGs up to quasi-isometry. This class we consider includes the hyperbolic RACGs classified in \cite{DaniThomasBowditch} and the RACGs on generalized theta graphs classified in \cite{hruska2020surface}.

With Feature 2, the structure invariant introduced by Cashen and Martin in \cite{CashenMartinStructureInvar}, see Section \ref{SectionStrucInvar}, comes into play: It can often detect that two JSJ trees of cylinders are non-isomorphic by using the types of the vertex stabilizers, implying that the RACGs are not quasi-isometric. With this technique occasionally one glance suffices to conclude that RACGs with rather basic defining graphs such as the following from Figure \ref{fig:EasyEx} are not quasi-isometric:

\begin{figure}[ht]
    \centering
    \scriptsize
    \def\svgwidth{230pt}
    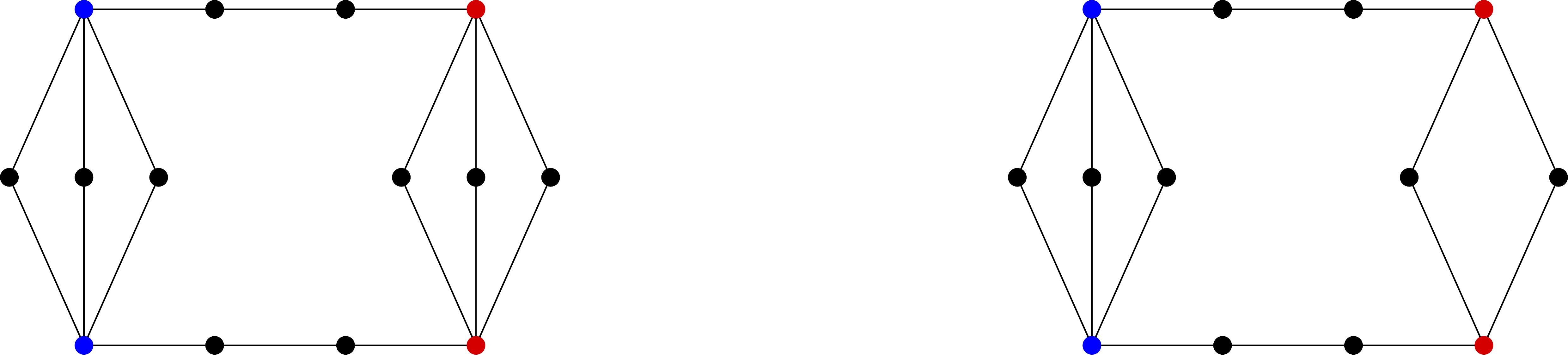
\end{figure}

The graph on the left has two uncrossed cut pairs, coloured in blue and red, which both have three common adjacent vertices. This implies that the corresponding cylinder vertices both have vertex groups that are the direct product of a virtually non-abelian free group and an infinite dihedral group. The red cut pair of the right graph, however, has only two common adjacent vertices. Thus the corresponding cylinder vertex group is virtually $\mathbb{Z}^2$. This is an obstruction for the existence of a quasi-isometry between the corresponding RACGs. \smallskip

In a second step, we adjust the structure invariant to our setting of RACGs by refining it in a way that also the converse of Feature 2 is true in certain cases, turning our (modified) structure invariant into a complete quasi-isometry-invariant:

\begin{thm} \textup{[cf.\ Theorem \ref{finalresult}]} \label{completeInvar}
Let W and $W^\prime$ be two finitely presented, one-ended RACGs with non-trivial JSJ decompositions over two-ended subgroups, both without rigid vertices. Define $T$ and $T^\prime$ to be the JSJ trees of cylinders of W and $W^\prime$ respectively. Then W and $W^\prime$ are quasi-isometric if and only if $T$ and $T^\prime$ have the same structure invariant up to reordering and quasi-isometry-equivalence of vertex groups.
\end{thm}

With this Theorem \ref{completeInvar} at hand, we can now immediately see that RACGs corresponding to defining graphs such as the following from Figure \ref{fig:EasyQIEx} are indeed quasi-isometric: 

\begin{figure}[ht]
    \centering
    \scriptsize
    \def\svgwidth{230pt}
    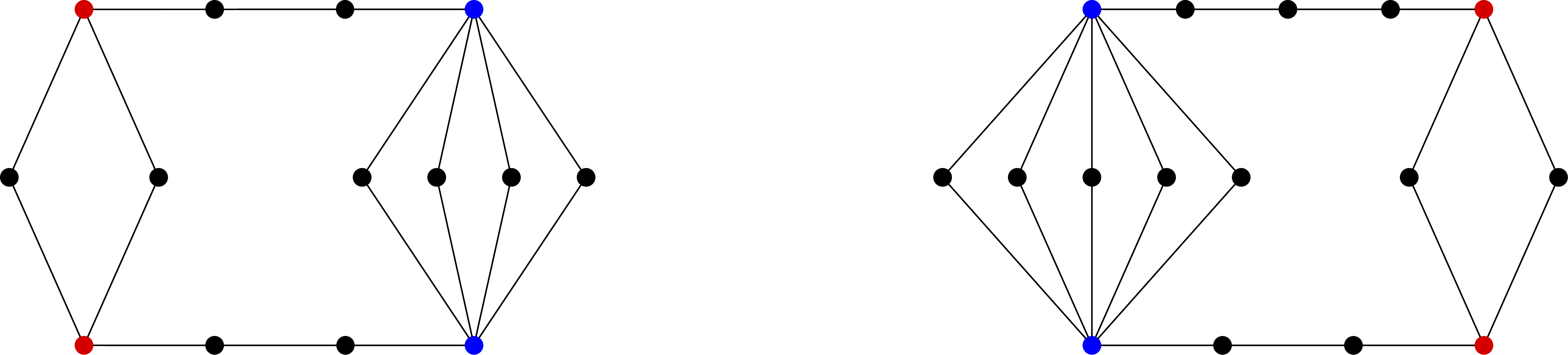
\end{figure}

Both graphs have one red uncrossed cut pair with two common adjacent vertices producing a virtually $\mathbb{Z}^2$ cylinder vertex group and a blue uncrossed cut pair with more than two common adjacent vertices producing a cylinder vertex group that is the direct product of a virtually non-abelian free group and an infinite dihedral group. So, the two defining graphs produce the same (modified) structure invariant.

Additionally, Theorem \ref{completeInvar} and its proof in Section \ref{SectionQIinvariant} can be exploited to obtain various examples of RACGs that are quasi-isometric, see Examples \ref{ExQI1} and \ref{ExQI2}: Starting from a defining graph, we perform reflections and duplications of subgraphs to produce new graphs whose corresponding RACGs are quasi-isometric to the original one. This method is even applicable to groups with rigid vertices, as long as these remain unaltered or have additional properties (see Remarks \ref{technicalities} and \ref{RemWithRig}).

\subsection{Background}

By introducing a geometric viewpoint in \cite{gromov1984infinite}, Gromov started the classification of groups in terms of their geometric equivalence. Groups are considered to be indistinguishable from a large-scale geometric perspective if there is a \textit{quasi-isometry (QI)} between them. That is a map which has two properties: It distorts distance at most by a scaling factor $C$ and an additive shift $D$ and it is almost surjective in the sense that in a uniform neighborhood of every point in target space we can find an image point of the map. If there is such a quasi-isometry between two geodesic metric spaces, we call the map a $(C,D)$\textit{-QI} and refer to the spaces as \textit{quasi-isometric}, short \textit{QI}, to each other.

A large class of interesting groups are the \textit{Coxeter groups}, introduced by Coxeter in \cite{coxeter1934discrete} as abstract reflection groups of geometric objects, see \cite{DavisCoxeter} for a recent survey. Their simplest examples are the \textit{right-angled Coxeter groups (RACGs)}, which are defined by a finite, simplicial, labelled graph, whose vertex labels are self-inverse generators of the group and whose edges determine commutation relations. They are called right-angled, because they act geometrically on a CAT(0) cube complex. This paper gives a QI-classification of a wide class of RACGs.

A strategy to produce QIs is to decompose groups into smaller pieces, whose QI-classification is understood. The interplay of the single pieces is captured by the \textit{graph of groups}, a graph equipped with vertex and edge groups. For an edge $e$, the corresponding edge group is contained in the vertex group of its initial vertex $o(e)$ and it embeds into the vertex group of its terminal vertex $t(e)$ via an attaching map. Stallings' theorem \cite{Sta71} states that a finitely generated group admits the simplest possible graph of groups decomposition as an \textit{HNN extension} or an \textit{amalgamated product} over a finite edge group if and only if it has more than one end. Therefore, since the number of ends of a group is a \textit{QI-invarant}, that is, preserved under a QI, so is the existence of such a decomposition.

Any finitely presented group admits a maximal decomposition over finite subgroups by Dunwoody's accessibility \cite{dunwoody1985accessibility}. Among the finitely presented groups with infinitely many ends, the collection of occurring QI-types of one-ended vertex groups in a maximal splitting is a QI-invariant by a result of Papasoglu and Whyte \cite[Theorem 0.4]{papasogluwhyte}. Thus we can focus on one-ended RACGs. In particular, the first obvious step is to consider one-ended groups that split over two-ended subgroups, a property which is a QI-invariant by \cite{papasoglu2005quasi} as long as the group is not commensurable to a surface group. In the case of RACGs, we restrict our attention to the finitely presented one-ended groups splitting over two-ended subgroups whose corresponding cube complex is additionally two-dimensional. Then, these properties can be easily ensured by restrictions on the defining graph, summarized in the \hyperref[StandingAssumption1]{Standing Assumption 1} in Section \ref{SectionPreliminaries}.

We consider splittings which are non-trivial and maximal in a certain sense. They are called \textit{JSJ decompositions}, produced from \textit{JSJ trees}. This terminology is borrowed from the decomposition of 3-manifolds. Its genesis is traced in \cite{JSJDecompGps}. JSJ decompositions are very robust under QIs: The QI-type of non-elementary vertex groups together with the pattern coming from the incident edge groups are preserved \cite[Section 2.3.2]{CashenMartinStructureInvar}. Thus, they can be used to distinguish groups up to QI. Unfortunately, one group might have plenty of JSJ decompositions. However, there is a canonical object, the \textit{JSJ tree of cylinders}, which can be built from any JSJ decomposition and thus captures the structure of the group. It has three different types of vertices: \textit{cylinder}, \textit{hanging} and \textit{rigid}. A QI between two groups induces an isomorphism between their JSJ trees of cylinders \cite[cf.][]{GuirardelLevittToC, CashenMartinStructureInvar}. In fact, the tree isomorphism even preserves additional information about the vertex groups such as the vertex type and about the structure provided by the adjacent edge groups. This was exploited in \cite{CashenMartinStructureInvar} by the introduction of the \textit{structure invariant}.\medskip

For \textit{hyperbolic}, one-ended, two-dimensional RACGs splitting over $D_\infty$-subgroups, Dani and Thomas give a QI-classification in \cite{DaniThomasBowditch}. While they claim to consider such RACGs which split over any two-ended subgroup, they implicitly use the assumption that the group does not split over $D_\infty \times \mathbb{Z}_2$. This error also occurred in an earlier version of this paper. Now this issue is addressed in Section \ref{SubsectionHypToC} and Theorem \ref{ToChypRACG} is a corrected version of the main Theorem 3.37 of \cite{DaniThomasBowditch}. In particular, Theorem \ref{ToChypRACG} gives an explicit \textit{visual} construction of the JSJ tree of cylinders of certain hyperbolic RACGs, that is, a construction which can be expressed only in terms of subgraphs of the defining graph. This implies that the defining graph not only determines the group presentation, but fully encodes its whole structure. Thus, essentially, certain hyperbolic RACGs can be distinguished up to QI just by looking at their defining graphs. More precisely, the distinction can be seen in the JSJ trees of cylinders or its quotient by the group action and follows from the finite valencies at the cylinder vertices of the tree.

Extending this QI-classification of certain hyperbolic RACGs from \cite{DaniThomasBowditch}, Hruska, Stark and Tran give a QI-classification for (not necessarily hyperbolic) RACGs whose defining graphs are \textit{generalized theta graphs} in \cite[Theorem 1.6]{hruska2020surface}. These results are combined in \cite[Theorem 5.20]{dani2018large} to a QI-classification of RACGs whose defining graphs are included in the much larger class of graphs dealt with in this paper.   

The construction followed in \cite{DaniThomasBowditch} is the one for \textit{Bowditch's JSJ tree}, introduced in \cite{BowditchJSJtree}, a special case of the JSJ tree. Its construction preceded the more general one defined in \cite{JSJDecompGps}. It works only for hyperbolic groups and by \cite[Theorem 9.18]{JSJDecompGps} it coincides with the JSJ tree of cylinders of the group. In this paper, the general construction and its broader set of tools are used.

\subsection{Outline}
After a short preface on right-angled Coxeter groups in Section \ref{SubsectionRACGs}, the general construction of JSJ trees of cylinders is introduced in Section \ref{SubsectionToC}. Then, Section \ref{SubsectionHypToC} analyzes the specific construction for hyperbolic RACGs from \cite{DaniThomasBowditch}. We conclude Section \ref{SectionPreliminaries} with a careful comparison of the two constructions. The proof of Theorem \ref{MainToC} on how to visually obtain the JSJ tree of cylinders for any one-ended, two-dimensional RACG splitting over two-ended subgroups stretches across all of Section \ref{SectionTreeOfCylRACGs} and is summarized in all detail as Theorem \ref{SummaryMainThm}. The proof has three main ingredients:
\begin{itemize}
    \item Section \ref{SubsectionCylVertices}: For the construction of cylinder vertices we use Proposition \ref{ThmUncrossedCutPairs}, essentially stating that they all come from the \textit{uncrossed cut collections} of the defining graph. In Lemma \ref{LemTypesofCylinderGroups} we show that the cylinder vertex groups are either virtually cyclic, virtually $\mathbb{Z}^2$ or the direct product of a virtually non-abelian free group and an infinite dihedral group.
    \item Section \ref{SubsectionNonCylVert}: The hanging and rigid vertices are produced by the analogy between the two constructions introduced in Section \ref{SectionPreliminaries}.
    \item Section \ref{SubsectionEdges}: The characterization of two-ended edge groups in terms of the defining graph from Theorem \ref{ThmEdgeStabTwoEnded} is a combination of Lemma \ref{HangCylEdges2Ended} and Theorem \ref{ThmK4Assump}.
\end{itemize}

Section \ref{SectionQIinvariant} is dedicated to the QI-classification of the RACGs. We can distinguish some of them up to QI by use of the structure invariant for JSJ trees of cylinders, whose construction and key features are illustrated in Section \ref{SectionStrucInvar}. Then, we are guided by the natural question: If two groups have equivalent structure invariants, when does this imply that they are indeed QI to each other?

The blueprint for producing such a QI is set up in
\cite{CashenMartinStructureInvar}: We need to understand the local QIs between vertex groups which are matched up by the structure invariant. These local QIs must also respect the structure coming from the incident edges. Then we try to patch those together inductively to a global QI.

The local QIs are produced in Section \ref{SubsectionLocalQIs}: For two-ended vertex groups, they are already dealt with in \cite{CashenMartinStructureInvar}; essentially they are determined by the finite valence of the corresponding vertex in the JSJ tree of cylinders. In Proposition \ref{VFD}, however, we see that the QIs between vertex groups that are the direct product of a virtually non-abelian free group and an infinite dihedral group can be chosen very flexibly. For the virtually abelian vertex groups, the QIs are of intermediate versatility, as proven in Proposition \ref{1VACyl}. The constraints they cause are implemented to the structure invariant by the \textit{density refinement} in Section \ref{SubsectionDensityRef}. That leads to a complete QI-invariant for certain RACGs, as outlined in Theorem \ref{completeInvar} and stated with all precision as Theorem \ref{finalresult} in Section \ref{SectionCompleteQI}. The proof works along the lines of the proofs of Sections 5 and 7 of \cite{CashenMartinStructureInvar}. As described in Outline \ref{methodQIEx}, Theorem \ref{finalresult} can be used as a tool to produce new examples of one-ended RACGs which are QI to each other. This is illustrated in Examples \ref{ExQI1} and \ref{ExQI2}. In fact, by Lemma \ref{LemCommens} the groups in these examples are even not commensurable. Thus, as far as the author is aware, they provide the first examples of one-ended, non-hyperbolic, non-commensurable, quasi-isometric RACGs.

\begin{ackn}  
First and foremost, I would like to thank my PhD supervisor Christopher Cashen for his expertise, his support and his guidance on and towards equally exciting as challenging mathematical questions. I am grateful to my co-advisor Goulnara Arzhantseva for providing an encouraging and enriching working environment in her research group as well as for our research seminar, which expands my mathematical horizon every week. Also, I would like to thank Emily Stark for offering comments on a preliminary version of this paper and the anonymous referee for
their many useful suggestions. Moreover, I greatly appreciate the opportunity to be supported by the Austrian Science Fund (FWF):  P34214-N.
\end{ackn}

\section{Preliminaries} \label{SectionPreliminaries}
We assume familiarity with standard group theoretic concepts such as Cayley graphs, quasi-isometries, ends and hyperbolicity of groups as wells as graphs of groups and Bass-Serre theory. For background consult for instance \cite{BogopolskiOleg2008Itgt, loeh2017geometric, SerreJeanPierre1980T}.

\subsection{RACGs} \label{SubsectionRACGs}
In the following section, we introduce the key properties of Right-Angled Coxeter groups.

\begin{defin}
For a finite, simplicial graph $\Gamma$, the \textit{defining graph}, on a vertex set $S$, we define the \textit{Right-Angled Coxeter Group (RACG)} $W_{\Gamma}$ as the group given by the following presentation
\begin{align*}
    W_{\Gamma} = \langle s \in S \mid s^2 = 1 \text{ for all } s \in S \, , (st)^2 = 1 \text{ if } (s,t) \in E(\Gamma)\rangle \, .
\end{align*}
\end{defin}

\begin{rem}
General Coxeter groups are often defined on the \textit{Coxeter graph} instead, which in the case of RACGs corresponds to the complement graph of the defining graph.
\end{rem}

\begin{nota}
Throughout the article, $W$ denotes a RACG and $\Gamma$ is a simplicial graph with vertex set $S = V(\Gamma)$. In order to emphasize the generating set $S$, we often denote the corresponding group as $W_S$ instead of $W_{\Gamma}$. Both notations are used without any further comment, depending on whichever is more suitable for the context.
\end{nota}

\begin{exam}
We obtain the following ‘extrema' as standard examples:
\begin{itemize}[noitemsep,topsep=2pt]
    \item If $\Gamma$ is a complete graph, then $W_{\Gamma} = \mathbb{Z}_2^{|S|}$. $W_\Gamma$ is finite if and only if $\Gamma$ is complete.
    \item If $\Gamma$ does not have any edges, then $W_{\Gamma} = *_{|S|} \, \mathbb{Z}_2$, so in particular, the infinite dihedral group $D_{\infty} = \mathbb{Z}_2 * \mathbb{Z}_2$ is a RACG.
\end{itemize}
\end{exam}

\begin{exam}
The class of RACGs is closed under taking direct products by taking the join of defining graphs and under taking free products by taking the disjoint union of defining graphs.
\end{exam}\vspace{6pt}

Certain subgroups can be "read off" the defining graph:

\begin{defin}
Given a RACG $W_S$ on $S = V(\Gamma)$, the subgroup $W_T$ generated by $T \subseteq S$ is called a \textit{special subgroup} of $W_S$.
\end{defin}

In fact, by Theorem 4.1.6 of \cite{DavisCoxeter}, $W_T$ is itself a (right-angled) Coxeter group on the defining graph $\Gamma_T$, which is the induced subgraph of $\Gamma$ on the vertices labelled by $T$. Moreover, the intersection of two special subgroups $W_{T} \cap W_{T^\prime}$ is the special subgroup generated by the intersection $T \cap T^\prime$.

The geometry of a Coxeter group $W_S$ is encoded in a complex, the so-called \textit{Davis complex}. Its construction and properties can be found in \cite{DavisCoxeter} and \cite[Section 2.1]{DaniThomasBowditch}. We outline the following facts relevant for this paper: The Davis complex of a special subgroup $W_T \subseteq W_S$ embeds isometrically into the Davis complex of $W_S$. For RACGs, the Davis complex is a CAT(0) cube complex. Its 1-skeleton is precisely the Cayley graph $Cay(W_S,S)$ of $W_S$ with respect to the generating set $S = V(\Gamma)$. Note that in case $W_S$ is infinite, it contains $D_\infty = W_{\{a,b\}}$ as a subgroup, where $a$ and $b$ are non-adjacent vertices in $S$. Then we find a geodesic of arbitrary length in the Cayley graph of $W_S$ that is labelled alternately by $a$ and $b$. We call such a geodesic \textit{bi-labelled}. \vspace{6pt} 

We also need some graph theoretic terminology:

\begin{defin}
A vertex $v$ of $\Gamma$ is \textit{essential} if it has at least valence 3. We denote the set of all essential vertices in $\Gamma$ by $EV(\Gamma)$. An embedded path between essential vertices, which does not contain any essential vertices in its interior, is a \textit{branch.}

A pair $\{a,b\}$ of vertices of $\Gamma$ is a \textit{cut pair} if it separates $\Gamma$, that is $\Gamma \setminus \{a,b\}$ has at least two connected components. If both vertices are essential, we call it an \textit{essential cut pair}.

A set $\{a,b,c\}$ of vertices of $\Gamma$ is called a \textit{cut triple} if $a$ and $b$ are not a cut pair, $c$ is a common adjacent vertex of $a$ and $b$ and the subgraph induced by $\{a,b,c\}$ separates $\Gamma$. 
\end{defin}

\begin{nota}
For economy of notation we will use the term \textit{cut collection} when referring to both cut pair and cut triple at once and use the notation $\{a-b\}$. The $-$ represents the possibly existing common adjacent vertex $c$ of $a$ and $b$ contributing to the triple.
\end{nota}

\begin{exam} \label{cutcollection}
In the left graph $\Gamma_1$ of Figure \ref{fig:ExampleCutPair}, the set $T_1 = \{a,b\}$ is a cut pair. Since $a$ and $b$ are not connected by an edge in $\Gamma_1$, the $T_1$-induced subgraph contains only two disconnected vertices, and thus the special subgroup generated by $T_1$ is $W_{\{a,b\}} = D_\infty$. The graph $\Gamma_2$ on the right has a cut triple $T_2 = \{a,b,c\}$. The special subgroup on the $T_2$-induced subgraph is $W_{\{a,b,c\}} = D_\infty \times \mathbb{Z}_2$.

\begin{figure}[ht]
    \centering
    \scriptsize
    \def\svgwidth{230pt}
    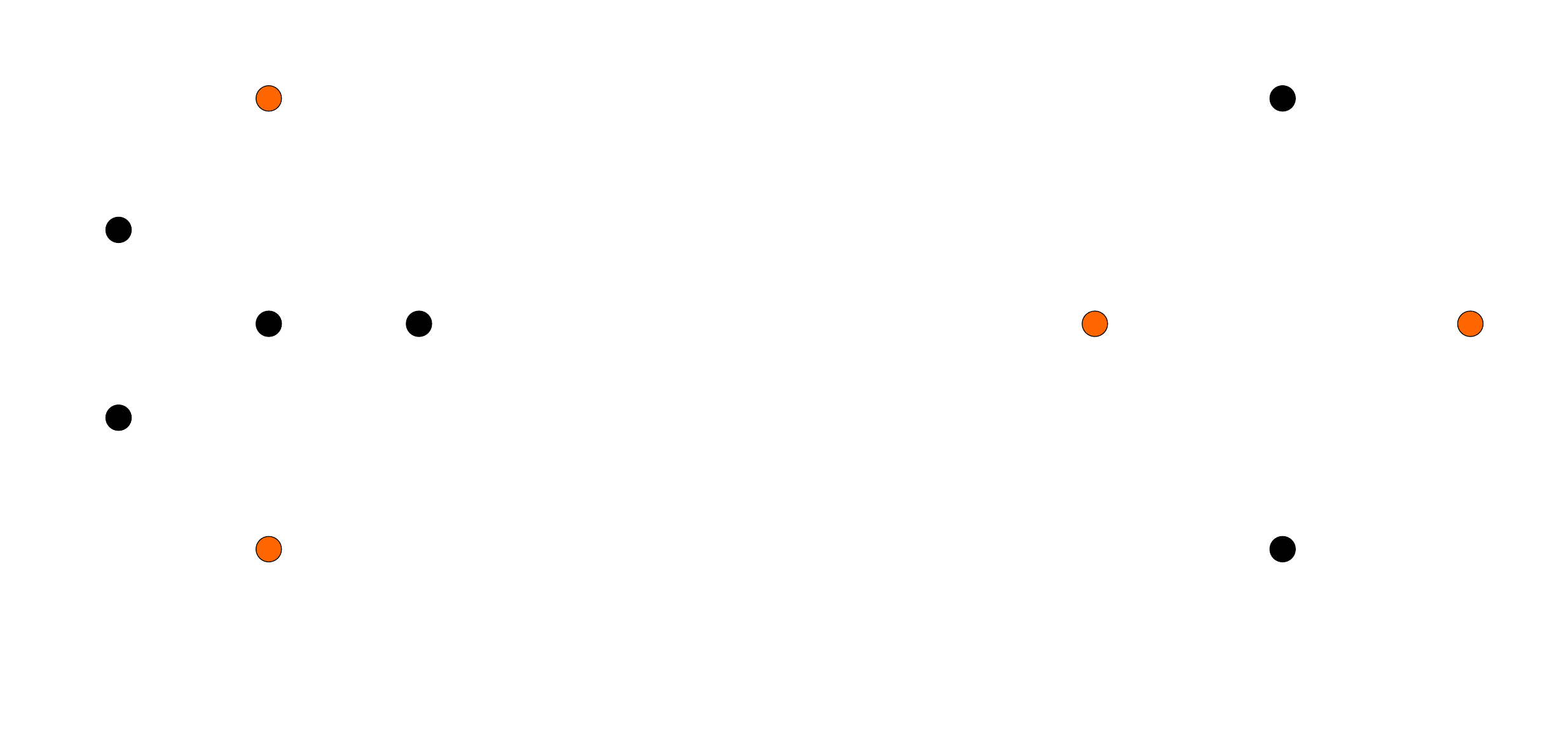
    \caption{The orange vertices form a cut pair and a cut triple respectively.}
   \label{fig:ExampleCutPair}
\end{figure}
\end{exam}

Recall that in the search for QIs, we want to limit ourselves to one-ended RACGs that split over two-ended subgroups. In order to translate these properties into accessible conditions on the defining graph $\Gamma$, we make use of the auxiliary assumption that the RACGs are two-dimensional. While we expect that this restriction can be dropped, the generalization is not immediate. This issue is also addressed in \cite[Section 1]{DaniThomasBowditch} and \cite[Question 5.17]{dani2018large}. Therefore, we fix the following properties of the defining graph $\Gamma$:

\begin{Assump} \label{StandingAssumption1}
The defining graph $\Gamma$
\begin{enumerate}[label=(\arabic*)]
    \item has no triangles: This corresponds to the fact that the Davis complex is two-dimensional, which simplifies the geometry encoded by the group.
    \item is connected and has neither a separating vertex nor a separating edge: $W_{\Gamma}$ is one-ended by \textup{\cite[Theorem 8.7.2]{DavisCoxeter}} under assumption (1) that $\Gamma$ has no triangles.
    \item has a cut collection $\{a-b\}$: By Theorem \ref{thmvisualsplitting}, recalling \textup{\cite[Theorem 1]{MihalikTschantzVisual}} in our setting under assumptions (1) and (2), this ensures the existence of a splitting over a two-ended subgroup. Indeed, if there is a cut collection $\{a-b\}$ all $k$ parts of $\Gamma \setminus \{a-b\}$ attach along the two-ended special subgroup $W_{\{a,b\}}=D_{\infty}$ or $W_{\{a,b,c\}}=D_{\infty}\times \mathbb{Z}_2$ as a $k$-fold amalgamated product.
    \item \label{StandingAssumption1Part4} is not a cycle of length $\geq 5$: By \cite[Theorem 4.2]{DaniThomasBowditch}, $\Gamma$ is a cycle of length $\geq 5$ if and only if $W_\Gamma$ is \textit{cocompact Fuchsian}. That means that it acts \textit{geometrically} on the hyperbolic plane. However, then the \textit{Švarc-Milnor-Lemma} implies that all groups with such a defining graph are QI to each other, thus their QI-Problem is understood.
\end{enumerate}
\end{Assump}

\begin{rem} Observe the following:
\begin{itemize}[noitemsep,topsep=2pt]
    \item Under \hyperref[StandingAssumption1]{Standing Assumtion 1}, a cut pair $\{a,b\}$ always consists of non-adjacent vertices and a cut triple $\{a,b,c\}$ forms a segment where $a$ and $b$ are both adjacent to $c$ and not adjacent to each other. Thus the special subgroup generated by both a cut pair and a cut triple is two-ended and the elements $a$ and $b$ generate a copy of $D_\infty$.
    \item For a cut triple $\{a-b\}$ the common adjacent vertex of $a$ and $b$ might not be unique: See for instance Figure \ref{fig:Overlap}, where $\{x,y,b\}$, $\{x,y,c\}$ and $\{x,y,d\}$ are cut triples. We say that the cut triples \textit{overlap}.
    \vspace{-10pt}
\begin{figure}[ht]
    \centering
    \scriptsize
    \def\svgwidth{365pt}
    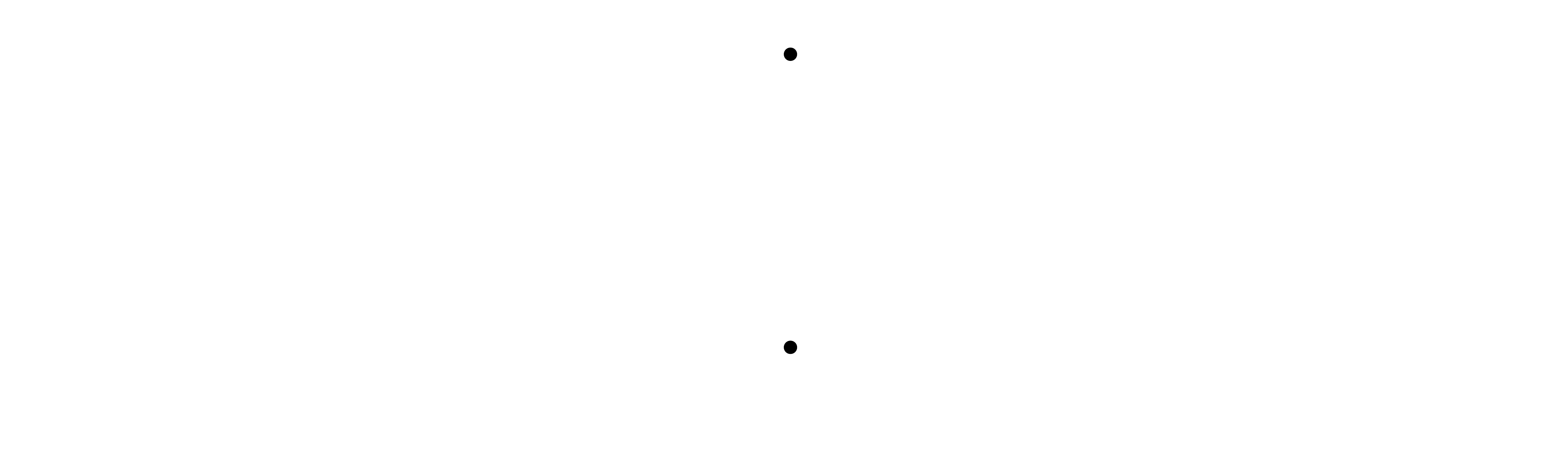
    \captionsetup{belowskip=-0.8cm}
    \caption{The vertex sets $\{x,y,b\}$, $\{x,y,c\}$ and $\{x,y,d\}$ form overlapping cut triples.}
    \label{fig:Overlap}
\end{figure}

    However, when there are overlapping cut triples the graph necessarily has an induced square, so this complication does not arise in the hyperbolic case (by \cite[Corollary 12.6.3]{DavisCoxeter}), but we have to deal with it in our more general setting. In Section \ref{TwoEndedEdgeGroups} we will make additional assumptions (to guarantee that the graph of cylinders has two-ended edge groups, see Remark \ref{RemCollapsedToC}) which exclude overlapping cut triples, see Remark \ref{NoOverlappingCutTriples}.
\end{itemize}
\end{rem}

\subsection{Trees of cylinders}\label{SubsectionToC}

Throughout this section, let $T$ be a simplicial tree and $G$ a finitely generated group acting on $T$ by isometries and without edge inversions. The stabilizer of any element $t$ in $T$ is denoted as $G_t$, geodesic paths in $T$ starting at vertex $a$ and ending at vertex $b$ are denoted as $[a,b]$. Let $\mathcal{A}$ be a class of infinite subgroups of $G$ that is stable under conjugation. $T$ is an $\mathcal{A}$\textit{-tree} if all the edge stabilizers $G_e$ of $T$ are contained in $\mathcal{A}$.

\begin{exam}
Since we split RACGs over two-ended subgroups, the class of subgroups we have in mind as $\mathcal{A}$ is the class $\mathcal{VC}$ of virtually infinite cyclic (or equivalently two-ended) subgroups. Note that $\mathcal{VC}$ is invariant under conjugation, but not under taking subgroups. 
\end{exam}\vspace{6pt}

Our main tool is a universal tree on which $G$ acts with vertex stabilizers as small as possible:

\begin{defin} \leavevmode
    \begin{enumerate}[nosep]
        \item A subgroup $H$ of $G$ is \textit{elliptic} in $T$ if it fixes a point in $T$. It is a \textit{universally elliptic} subgroup if it fixes a point in any $\mathcal{A}$-tree. An $\mathcal{A}$-tree is \textit{universally elliptic} if all its edge stabilizers are universally elliptic subgroups of $G$. 
        \item An $\mathcal{A}$-tree $T$ \textit{dominates} another $\mathcal{A}$-tree $T^\prime$ if every vertex stabilizer of $T$ is elliptic in $T^\prime$.
        \item A \textit{JSJ tree} of $G$ is an $\mathcal{A}$-tree $T$ that is universally elliptic and that dominates any other universally elliptic $\mathcal{A}$-tree $T^\prime$. The quotient graph $\Lambda = T/G$ is called a \textit{JSJ decomposition} or \textit{JSJ splitting} of $G$. 
    \end{enumerate}
\end{defin}

JSJ trees are extensively surveyed in \cite{JSJDecompGps}. Unfortunately, the JSJ tree is not as universal as we would like it to be. It does not even always exist, nor is it unique if it does. It rather happens that we find a collection of universally elliptic trees, which are pairwise dominating each other. This collection then is called the \textit{JSJ deformation space} \cite[Section 2.3]{JSJDecompGps}.

\begin{rem} \label{allunivellipiticsubgps}
If in a graph of groups $\Lambda = T/G$ of $G$ whose edge groups are all universally elliptic, also up to conjugation all universally elliptic subgroups of $G$ occur as edge groups, $\Lambda$ is a JSJ decomposition of $G$. Indeed, if all edge groups in $\Lambda$ are universally elliptic, so are the stabilizers of $T$, thus $T$ is universally elliptic. Furthermore, given any other universally elliptic tree $T^\prime$, we can refine it to $T$, and $T^\prime$ is therefore dominated by $T$ \cite[Lemma 2.15]{JSJDecompGps}. Thus $T$ is a JSJ tree.
\end{rem}

We aim to obtain a more accessible equivalent definition, when restricting to one-ended groups splitting over two-ended subgroups. For that we need to introduce the following terminology:

\begin{defin} \citep[cf.][Definition 5.13]{JSJDecompGps} \label{defHanging}
A vertex $v$ of a graph of groups $\Lambda$ over two-ended edge groups and its vertex group $G_v$ are called \textit{hanging} if $G_v$ maps onto the fundamental group $\pi_1(\Sigma_v)$ of a hyperbolic, compact, two-dimensional orbifold $\Sigma_v$ and the image of every edge group incident to $G_v$ in $\pi_1(\Sigma_v)$ is either finite or contained in a boundary subgroup of $\pi_1(\Sigma_v)$. We call $v$ and $G_v$ \textit{maximal hanging} if there is no other hanging vertex group $G_w$ such that the corresponding orbifold $\Sigma_w$ can be glued to $\Sigma_v$ along identical boundary components to obtain a new splitting of the group.
\end{defin}

\begin{rem} 
While the interpretation of a \textit{hanging} subgroup is not universal, in the setting of RACGs all versions are equivalent: For instance, suppose that, following \cite{BowditchJSJtree}, we see a vertex group $G_v$, which is non-elementary, finitely generated and which acts properly discontinuously on the hyperbolic plane $\mathbb{H}^2$. This is equivalent to saying that $G_v$ surjects with finite kernel onto the fundamental group of a hyperbolic, compact, two-dimensional orbifold $\Sigma_v$ \cite[cf.][Defnition 3.2.]{barrett2018computing}. If additionally all the maximal two-ended subgroups of $G_v$ map onto the fundamental groups of the boundary components of $\Sigma_v$, Bowditch calls $G_v$ \textit{hanging Fuchsian}. However, then $G_v$ meets the Definition \ref{defHanging} of a hanging vertex group as well.

Also, it is worth noting that in their Definition 5.13 in \cite{JSJDecompGps}, Guirardel and Levitt define the vertex and vertex group we call hanging as \textit{quadratically hanging (QH)}, to extend the definition of quadratically hanging subgroups given by Rips and Sela in \cite{RipsE1997CSoF}. Moreover, various authors call vertex groups meeting the properties of Definition \ref{defHanging} along similar lines as the hanging Fuchsian groups, \textit{hanging surface groups} for instance.
\end{rem}

\begin{defin}
A vertex $v$ of a graph of groups $\Lambda$ over two-ended edge groups and its vertex group $G_v$ are called \textit{rigid} if $G_v$ is not two-ended, not hanging and does not split over a two-ended subgroup relative to its incident edge groups.
\end{defin}

By piecing together Theorem 6.5, Corollary 6.3, Section 2.6 and Proposition 5 of \cite{JSJDecompGps}, which rely on work of Fujiwara and Papasoglu \cite{fujiwara2006jsj}, we can describe certain JSJ decompositions neatly in terms of graphs of groups:

\begin{lem} \label{EqDefJSJdecomp}
If $G$ is a finitely presented, one-ended group, then a graph of groups decomposition with two-ended edge groups is a JSJ decomposition if and only if the following conditions hold:
\begin{itemize}[nosep]
    \item Each vertex group is either two-ended, hanging or rigid.
    \item Any valence one vertex $v$ with two-ended vertex group does not have an incident edge group surjecting onto $G_v$.
    \item All hanging vertex groups are maximal.
\end{itemize}
\end{lem}

Even though JSJ decompositions are not unique, under certain conditions, we can produce a canonical representative of the JSJ deformation space, the so-called \textit{tree of cylinders} $T_c$. The rest of this subsection gives a short overview of its construction. For all details consult \cite{GuirardelLevittToC}.

\begin{defin}
An equivalence relation $\sim$ on $\mathcal{A}$ is called \textit{admissible} if for all $A, B \in \mathcal{A}$ the following axioms hold:
\begin{enumerate}[noitemsep,topsep=2pt]
    \item If $A \sim B$ and $g \in G$ then $gAg^{-1} \sim gBg^{-1}$.
    \item If $A \subseteq B$, then $A \sim B$.
    \item Given an $\mathcal{A}$-tree $T$ and $a,b \in V(T)$ that are fixed by $A,B \in \mathcal{A}$ respectively, then for every edge $e \subseteq [a,b]$ we have $A \sim G_e \sim B$. \label{AdmisRelAxiom3}
\end{enumerate}
\end{defin}

\begin{defin}
Two subgroups $H$ and $K$ of a group $G$ are called \textit{commensurable} if their intersection $H \cap K$ has finite index in both $H$ and $K$. The \textit{commensurator} of a subgroup $H$ in $G$ is the set
\begin{gather*}
    Comm_{\,G}(H) = \{ g \in G \mid  gHg^{-1} \text{ and } H \text{ are commensurable}\} \, .
\end{gather*}
Commensurability is an equivalence relation on subgroups. We denote the equivalence class of $A \in \mathcal{A}$ by $[A]$. The stabilizer of $[A]$ under the action of $G$ on $\mathcal{A}\bign/\sim$ by conjugation is denoted as $G_{[A]}$.
\end{defin}

\begin{exam} \label{ExCommAdmRel}
On the class $\mathcal{VC}$ of two-ended subgroups of $G$, commensurability is an admissible equivalence relation. For $A \in \mathcal{VC}$ we obtain $G_{[A]} = Comm_{\,G}(A)$. 
\end{exam}

\begin{constr} \label{constructionToc} 
Given an $\mathcal{A}$-tree $T$, we construct the object of interest, the \textit{cylinder}, in the following way:
\begin{itemize}[noitemsep,topsep=2pt]
    \item Start with an admissible equivalence relation $\sim$ on $\mathcal{A}$.
    \item Define two edges $e, f \in E(T)$ to be equivalent if their edge stabilizers $G_e$ and $G_f$ are equivalent, i.e. $e \sim f$ if $G_e \sim G_f$.
    \item If $G_e$ fixes the edge $e \in E(T)$, in particular it fixes its endpoints $o(e), t(e) \in V(T)$. Thus by axiom \hyperref[AdmisRelAxiom3]{(3)} for an admissible relation, all the edges on a path between two equivalent edges will be in the same equivalence class as well. Thus this equivalence class forms a subtree $Y$ of $T$, called a \textit{cylinder} of $T$.
    \item By construction, two distinct cylinders can intersect at most in one common vertex.
    \item We refer to the equivalence class in $\mathcal{A}\bign/\sim$ containing all edge stabilizers of edges in $Y$ as $[Y]$.
\end{itemize}
\end{constr}

\begin{defin}
Given an admissible equivalence relation on $\mathcal{A}$ and an $\mathcal{A}$-tree $T$, its \textit{tree of cylinders} $T_c$ is the following bipartite tree with vertex set $V_1 \sqcup V_2$: The vertex set $V_1$ contains one vertex $v_Y$ per cylinder $Y$, the \textit{cylinder vertices}. The vertex set $V_2$ contains all the vertices of $T$ that belong to at least two cylinders. There is an edge $(v_Y,v) \in E(T_c)$ between $v_Y$ and every vertex $v$ contained in $Y$. The graph of groups decomposition of $G$ coming from the quotient of the action of $G$ on $T_c$ is the \textit{graph of cylinders} $\Lambda_c$.
\end{defin}

The stabilizer $G_Y$ of a cylinder vertex $v_Y$ in $V_1$ is $G_{[Y]}$. The stabilizer $G_v$ of a vertex $v$ in $V_2$ is the stabilizer $G_v$ of $v$ viewed as a vertex of $T$. An edge $(v_Y,v)$ in $E(T_c)$ is stabilized by the intersection of $G_{[Y]}$ and $G_v$.

\begin{exam}
Consider the Baumslag-Solitar group $BS(m,n) = \langle a, b \mid b^{-1}a^mb = a^n \rangle$ defined for the integers $m,n \in \mathbb{Z}\setminus\{0\}$. We view it as an HNN-extension with stable letter $b$ and consider its action on the corresponding Bass-Serre tree $T$. All the edge stabilizers are of the form $g\langle a^m \rangle g^{-1}$ for $g \in BS(m,n)$, thus they are contained in $\mathcal{VC}$. By use of the inductive consequence
$$b^{-k}a^{m^ky}b^k = a^{n^ky} \quad \text{ for any } k, y \in \mathbb{N} $$  of the relation, one shows that $\langle a^m \rangle$ is commensurable to $g\langle a^m \rangle g^{-1}$ for any $g \in BS(m,n)$. Hence, all edges are part of the same commensurability-cylinder and $T_c$ consists of only one vertex.
\end{exam}

\begin{rem} \label{RemCollapsedToC}
Note that $T_c$ is not necessarily an $\mathcal{A}$-tree. This problem is resolved  by collapsing all edges that have edge stabilizers not in $\mathcal{A}$ to obtain the \textit{collapsed tree of cylinders} $T_c^*$. However, in our application of the construction, we aim to bypass this complication.
\end{rem}

\begin{nota}
Henceforth, when the set $\mathcal{A}$ and the admissible equivalence relation on it are not specified, it is fixed to be $\mathcal{VC}$ with the relation to be commensurability, as in Example \ref{ExCommAdmRel}.
\end{nota}

The question left to answer is how the construction of the tree of cylinders gives a canonical object encoding the structure of the group. Starting from a finitely presented, one-ended group $G$, we pick some JSJ tree $T$ of the JSJ deformation space, which exists by \cite[Theorem 1]{JSJDecompGps}. For $T$ we construct the tree of cylinders $T_c$, which by \cite[Theorem 1]{GuirardelLevittToC} does not depend on the choice of $T$ but only on the deformation space itself. Thus it makes sense to call it \textit{the JSJ tree of cylinders} and the corresponding graph of cylinders $\Lambda_c$ \textit{the JSJ graph of cylinders}. While for instance for hyperbolic groups, $\Lambda_c$ is itself a JSJ decomposition \cite[Theorem 9.18]{JSJDecompGps}, this is not true in general. However, by construction its Bass-Serre tree is $G$-equivariantly isomorphic to the tree of cylinders of any JSJ tree. Hence, from the JSJ graph of cylinders $\Lambda_c$, we can essentially recover the deformation space of JSJ splittings. \medskip

Moreover the JSJ tree of cylinders produces a QI-invariant for groups, by a result of Cashen and Martin based on work of Papasoglu \cite[Theorem 7.1]{papasoglu2005quasi}  with a correction made by Shepherd and Woodhouse in \cite{shepherd2020quasi}:
\begin{thm} \textup{\citetext{\citealp[Theorem 2.9]{CashenMartinStructureInvar}; \citealp[Theorem 2.8]{shepherd2020quasi}}} \label{QIinvar} Given two finitely presented, one-ended groups $G$ and $G^\prime$ splitting over two-ended subgroups which are quasi-isometric via $\phi: G \rightarrow G^\prime$, then $\phi$ induces a tree isomorphism $\phi_*: T_c \rightarrow T_c^\prime$. Moreover, $\phi_*$ is vertex-type preserving and for every vertex $v \in V(T)$ with vertex group $G_v$ there is a real constant $C_v \geq 0$ such that $\phi$ maps $G_v$ within distance $C_v$ of $G_{\phi_*(v)}^\prime$.
\end{thm} 

Thus, ideally, we construct the JSJ graphs of cylinders directly from the groups, in our case from the defining graphs of the RACGs. 
Deducing from them that the corresponding JSJ trees of cylinders are not isomorphic then implies that the groups we started with are not QI. On the other hand, if there is an isomorphism between the JSJ trees of cylinders, we try to promote it to a QI of the groups. 

\begin{out} \label{ExFramework}
To summarize, the framework we focus on is the following: The group $G$ is finitely presented, one-ended and splits over the set of two-ended subgroups $\mathcal{VC}$. We obtain a JSJ splitting $\Lambda$, in which all vertex groups are either two-ended, hanging or rigid by Lemma \ref{EqDefJSJdecomp}. By considering the commensurability relation on the corresponding JSJ tree, we produce the JSJ graph of cylinders $\Lambda_c$, whose cylinder vertex groups are the commensurators of the two-ended groups of $\Lambda$ and whose non-cylinder vertex groups are precisely the hanging and rigid vertices of $\Lambda$.
\end{out}

\subsection{JSJ trees of cylinders of hyperbolic RACGs} \label{SubsectionHypToC}
For one-ended, two-dimensional, hyperbolic RACGs whose defining graphs do not have any cut triples, a way to construct the JSJ graph of cylinders directly from the defining graph $\Gamma$ is given in \cite{DaniThomasBowditch}. By \cite[Corollary 12.6.3]{DavisCoxeter} a RACG $W_\Gamma$ is hyperbolic if and only if $\Gamma$ has no squares. Although Dani and Thomas's construction follows the one for \textit{Bowditch's JSJ tree} described in \cite{BowditchJSJtree}, it turns out that the tree they produce in their (main) Theorem 3.37 corresponds to the JSJ tree of cylinders of $W_\Gamma$. More precisely, since $W_\Gamma$ is hyperbolic, it follows from \cite[Theorem 9.18]{JSJDecompGps} that both trees and thus their corresponding decompositions are $W_\Gamma$-equivariantly isomorphic.

Dani and Thomas claim in \cite{DaniThomasBowditch} that they give a construction of Bowditch's JSJ tree for \textit{all} one-ended, two-dimensional, hyperbolic RACGs splitting over two-ended subgroups. However, they miss the fact that a RACG can not only split over a two-ended $D_\infty$-subgroup coming from a cut pair but also over a two-ended $D_\infty \times \mathbb{Z}_2$-subgroup coming from a cut triple. The origin of this problem is a miscitation of Theorem 1 of \cite{MihalikTschantzVisual} as Theorem 2.1 in \cite{DaniThomasBowditch} claiming that every splitting over a two-ended subgroup corresponds to a cut pair. Example \ref{cutcollection} gives a counterexample to this claim. 

However, under the mild additional assumption that the defining graph $\Gamma$ does not have any cut triples, all the results in \cite{DaniThomasBowditch} remain valid. We will add this assumption whenever referring to results in \cite{DaniThomasBowditch}. This additional assumption was also implicitly used in an earlier version of this paper, however the error has been removed as Theorem \ref{SummaryMainThm} now also includes the construction of the JSJ tree of cylinders of RACGs splitting over two-ended subgroups coming from cut triples in both the hyperbolic and non-hyperbolic case. In particular removing this assumption does not affect the strategy and large-scale geometry results developed in \cite{DaniThomasBowditch} and in this paper, but only certain descriptions of the subgraphs of $\Gamma$ corresponding to the large-scale structures of interest.

\medskip

Most of the proofs in \cite{DaniThomasBowditch} depend on the hyperbolicity, in particular the existence of the Gromov boundary, of the group $W_\Gamma$. Before we can produce the broader result, we want to understand the correspondence between the two constructions. This subsection is dedicated to this task.

In our terminology, the JSJ tree of cylinders of one-ended, two-dimensional, hyperbolic RACGs splitting over $D_\infty$-subgroups is produced by the following theorem:

\begin{thm}\textup{\citep[cf.][Theorem 3.37]{DaniThomasBowditch}} \label{ToChypRACG}
Let $W_\Gamma$ be a hyperbolic RACG with $\Gamma$ satisfying the \hyperref[StandingAssumption1]{Standing Assumption 1} and in addition let $\Gamma$ have no cut triples. Then its JSJ tree of cylinders $T_c$ has vertices and associated vertex groups in the JSJ graph of cylinders $\Lambda_c$ as follows:
\begin{enumerate}
    \item Type 1 vertex:
    \begin{enumerate}[label=(\alph*)]
        \item For any cut pair $\{a,b\}$ such that $\Gamma \setminus \{a,b\}$ has $k \geq 3$ connected components, none of which consists of only one single vertex, there is a valence $k$ vertex in $T_c$. The associated vertex group in $\Lambda_c$ is the subgroup of $W_\Gamma$ generated by $\{a, b\}$, unless $a$ and $b$ have a common adjacent vertex $c$, then it is generated by $\{ a, b, c \}$. \label{1a}
        \item For any cut pair  $\{a,b\}$ such that $\Gamma \setminus \{a,b\}$ has $k \geq 3$ connected components, one of which consists of only one vertex $c$, there is a valence $2 \cdot (k-1)$ vertex in $T_c$. The associated vertex group is the subgroup of $W_\Gamma$ generated by $\{ a, b, c \}$. \label{1b}
        \item For any set $A \subseteq V(\Gamma)$ satisfying the properties \hyperref[A1]{(A1)}, \hyperref[A2]{(A2)} and \hyperref[A3]{(A3)} and which generates a two-ended subgroup not occurring in \hyperref[1a]{1.(a)} nor in \hyperref[1b]{1.(b)}, there is a valence $2$ vertex in $T_c$, where the properties \hyperref[A1]{(A1)}, \hyperref[A2]{(A2)} and \hyperref[A3]{(A3)} are the following:
        \begin{itemize}[leftmargin=1.5cm]
            \item[(A1)] \label{A1} Elements of $A$ pairwise separate the geometric realization $|\Gamma|$. 
            \item[(A2)] If any subgraph $\Gamma^\prime$ of $\Gamma$ that is a subdivided $K_4$ contains more than 2 vertices of $A$, all vertices of $A$ lie on the same branch of $\Gamma^\prime$. \label{A2}
            \item[(A3)] The set $A$ is maximal among all sets satisfying \hyperref[A1]{(A1)} and \hyperref[A2]{(A2)}. \label{A3}
        \end{itemize}
        If either $|A| = 2$ and there is no third vertex $c$ adjacent to both elements in $A$ or $|A| = 3$, the associated vertex group in $\Lambda_c$ is the subgroup of $W_\Gamma$ generated by $A$. If $|A| = 2$ and the two elements in $A$ have a common adjacent vertex $c$, then the associated vertex group in $\Lambda_c$ is the subgroup of $W_\Gamma$ generated by $A \cup \{c\}$. \label{1c}
        \item On any edge between a type 2 and a type 3 vertex there is a valence $2$ vertex added in $T_c$. The associated vertex group in $\Lambda_c$ is the intersection of the vertex groups of its neighbors. \label{1d}
    \end{enumerate}
    \item Type 2 vertex: For any set $A \subseteq V(\Gamma)$ satisfying properties \hyperref[A1]{(A1)}, \hyperref[A2]{(A2)} and \hyperref[A3]{(A3)} such that the subgroup generated by $A$ is infinite but not two-ended, there is a vertex in $T_c$ with associated vertex group $W_A$ in $\Lambda_c$. \label{2}
    \item Type 3 vertex: For any set $B \subseteq EV(\Gamma)$ of essential vertices in $\Gamma$ satisfying the properties \hyperref[B1]{(B1)}, \hyperref[B2]{(B2)} and \hyperref[B3]{(B3)}, there is a vertex in $T_c$ whose associated vertex group in $\Lambda_c$ is the subgroup $W_B$ generated by $B$, where the properties \hyperref[B1]{(B1)}, \hyperref[B2]{(B2)} and \hyperref[B3]{(B3)} are the following:
        \begin{itemize}[leftmargin=1.5cm]
            \item[(B1)] For any pair $C = \{c,d\} \subseteq EV(\Gamma)$ of essential vertices, $B \setminus C$ is contained in one single connected component of $\Gamma \setminus C$. \label{B1}
            \item[(B2)] The set $B$ is maximal among all sets satisfying \hyperref[B1]{(B1)}. \label{B2}
            \item[(B3)] $|B| \geq 4$. \label{B3}
        \end{itemize} \label{3}
\end{enumerate}
Between a vertex $v$ of type 1 and a vertex $v^\prime$ of type 2 or 3 in $V(T_c)$, there is an edge connecting them if and only if their corresponding vertex groups intersect in an infinite subgroup.
\end{thm}

\begin{nota}
Whenever we illustrate a JSJ graph of cylinders $\Lambda_c$ of a RACG, see Figure \ref{fig:JSJHypCase} for instance, for economy of notation we omit the brackets of the vertex and edge groups and just write down the collection of generating vertices. For convenience, we mark cylinder vertices in green, hanging vertices in red and rigid vertices in blue.
\end{nota}

\begin{rem}
Not only type \hyperref[1a]{1.(a)} or type  \hyperref[1b]{1.(b)} vertices correspond to essential cut pairs, but all type 1 vertices in Theorem \ref{ToChypRACG} do:

Any set $A$ satisfying properties \hyperref[A1]{(A1)}, \hyperref[A2]{(A2)} and \hyperref[A3]{(A3)} must contain an essential cut pair as shown in Lemma \ref{LemVerticesInA}. 
But for a vertex of type \hyperref[1c]{1.(c)}, we need that $W_A$ is two-ended. By Theorem \ref{VCylic} we see that the only two options for a special subgroup of $\Gamma$ satisfying \hyperref[StandingAssumption1]{Standing Assumption 1} to be two-ended is that it consists either of two single non-adjacent vertices or of two vertices connected via one common adjacent vertex. So either $|A|=2$, then it is precisely an essential cut pair. Or $|A|=3$, thus it contains an essential cut pair and one common adjacent vertex in between. 

By \cite[Lemma 3.30]{DaniThomasBowditch}, the intersection of a set $A$ satisfying properties \hyperref[A1]{(A1)}, \hyperref[A2]{(A2)} and \hyperref[A3]{(A3)} and a set $B$ satisfying properties \hyperref[B1]{(B1)}, \hyperref[B2]{(B2)} and \hyperref[B3]{(B3)} contains at most two vertices. Thus $A$ and $B$ can intersect at most in an essential cut pair. But in case their associated vertex groups intersect non-trivially, this intersection cannot be finite, implying that it must contain precisely the essential cut pair. The vertex of type \hyperref[1d]{1.(d)} can therefore be detected from an essential cut pair as well.

However, it is important to note that not all essential cut pairs contribute to a type 1 vertex, as illustrated in Example \ref{ExamplesK4EssentialCutPair}. The question on how to distinguish the ones contributing from the ones that do not is addressed in Section \ref{SectionTreeOfCylRACGs} in Proposition \ref{ThmUncrossedCutPairs}. 
\end{rem}

\begin{exam} \label{ExamplesK4EssentialCutPair}
In Figure \ref{fig:JSJHypCase} we see on the left side a square-free graph $\Gamma$ satisfying the \hyperref[StandingAssumption1]{Standing Assumption 1}. On the right side, the JSJ graph of cylinders $\Lambda_c$ of $W_\Gamma$ is illustrated. It is obtained by Theorem \ref{ToChypRACG} with the following considerations: There is no cut pair of type \hyperref[1a]{1.(a)} and the cut pairs $\{u,y\}$ and $\{v,y\}$ give vertices of type \hyperref[1b]{1.(b)}. From the cut pairs $\{v,w\}$ and $\{w,x\}$ we obtain a vertex of type \hyperref[1c]{1.(c)} and $\{v,x\}$, $\{w,z\}$ and $\{y,z\}$ add vertices of type \hyperref[1d]{1.(d)}. Of type \hyperref[2]{2}, there are the five vertex sets $\{v, n_1, n_2, x\}$, $\{w, r_1, r_2, z\}$, $\{y, s_1, s_2, z\}$, $\{u,p_1,p_2,y\}$ and $\{v, l_1, l_2, u, y\}$. The only vertex of type \hyperref[3]{3} is given by $\{v,w,x,y,z\}$. Note that for instance the set $\{v,x,w\}$ does not give a vertex of type \hyperref[2]{2} as property \hyperref[A2]{(A2)} is violated by the subdivided $K_4$ with corners $w, v, x$ and $z$. Furthermore, while the set $\{v, l_1, l_2, u\}$ is a pairwise separating branch, it does not satisfy \hyperref[A3]{(A3)}. Thus, even though $\{u,v\}$ is an essential cut pair and thus gives a two-ended subgroup over which $W_\Gamma$ splits, it does not give a type 1 vertex. As proved in Proposition \ref{ThmUncrossedCutPairs} this is due to the fact that there are other cut pairs, for instance $\{y,l_1\}$,  separating $u$ and $v$ (see also Example \ref{ExCrossedCutPair}). This implies that $W_{\{u,v\}}$ is not universally elliptic and therefore contained within a hanging subgroup.  

\begin{figure}[ht]
    \centering
    \scriptsize
    \def\svgwidth{350pt}
    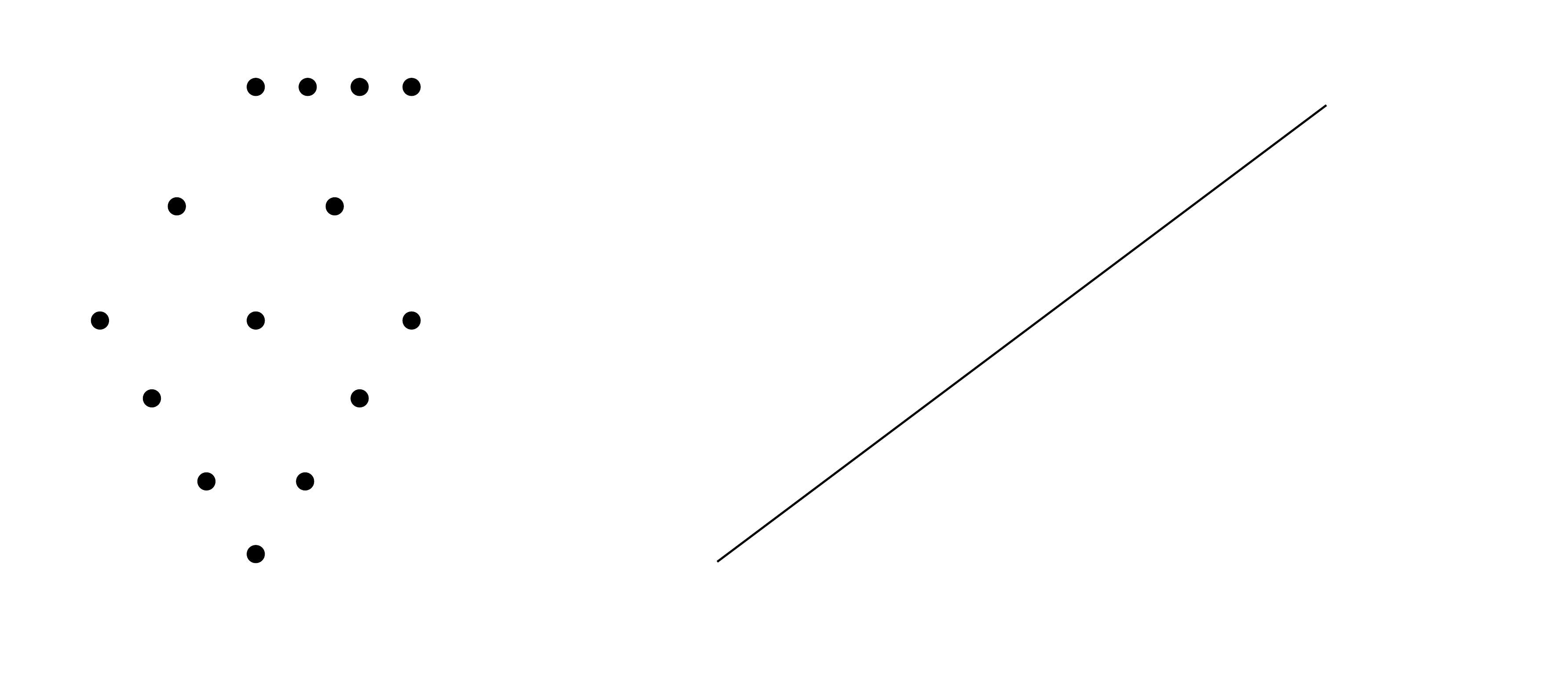
    \caption{$\Lambda_c$ is the JSJ graph of cylinders of the RACG $W_\Gamma$ obtained by Theorem \ref{ToChypRACG}.}
   \label{fig:JSJHypCase}
\end{figure}
\end{exam}

Also, the type of a vertex determines a key property:
\begin{thm} \textup{\citep[cf.][Theorem 5.28]{BowditchJSJtree}}
\label{BowditchJSJProperties}
Let $W_\Gamma$ be a hyperbolic RACG with $\Gamma$ satisfying the \hyperref[StandingAssumption1]{Standing Assumption 1}. Let  $\Lambda_c$ be the JSJ graph of cylinders given by Theorem \ref{ToChypRACG}. Then:
\begin{itemize}[nosep]
    \item The vertex group associated to a type 1 vertex is two-ended.
    \item The vertex group associated to a type 2 vertex is hanging. \label{Type2}
    \item The vertex group associated to a type 3 vertex is rigid. \label{Type3}
\end{itemize}
\end{thm}

Comparing this Theorem \ref{BowditchJSJProperties} and Outline \ref{ExFramework}, we can now establish the following correspondence between the JSJ tree of cylinders given by Construction \ref{constructionToc} and the JSJ tree of cylinders constructed in Theorem \ref{ToChypRACG} for hyperbolic RACGs:
\begin{itemize}
    \item \textit{Type 1 vertices correspond to cylinder vertices:} Type 1 vertices in $\Lambda_c$ lift to vertices of finite valence in $T_c$ with two-ended vertex stabilizer. These properties can only hold for cylinder vertices. Furthermore, by existence of vertices of type 1.(d), the JSJ tree of cylinders constructed in Theorem \ref{ToChypRACG} is bipartite with $V = V(1) \sqcup V(2,3)$, where $V(1)$ are all the vertices of type 1 and $V(2,3)$ contains vertices of type 2 and 3. Thus, as the JSJ tree of cylinders is also bipartite, no other than the type 1 vertices can correspond to cylinder vertices. 
    \item \textit{Type \hyperref[2]{2} vertices correspond to hanging vertices:} By Theorem \ref{BowditchJSJProperties} type 2 vertices are hanging, thus they are the hanging non-cylinder vertices. 
    \item \textit{Type \hyperref[3]{3} vertices correspond to rigid vertices:} Again,  by Theorem \ref{BowditchJSJProperties} type 3 vertices are rigid, thus they are the rigid non-cylinder vertices. 
\end{itemize}

\section{JSJ trees of cylinders of RACGs} \label{SectionTreeOfCylRACGs}
Since in the non-hyperbolic case, there is no universal construction of a JSJ tree and its tree of cylinders like the one given by Bowditch, for arbitrary RACGs we need to start from scratch: We first determine how to find a JSJ decomposition in terms of the defining graph $\Gamma$ and then produce a construction of the JSJ graph of cylinders from there.

In fact, any decomposition of a (right-angled) Coxeter group over two-ended subgroups is visible in the defining graph $\Gamma$:

\begin{thm} \textup{\citetext{\citealp[Theorem 1]{MihalikTschantzVisual}}} \label{thmvisualsplitting}
For a simplicial graph $\Gamma$ which is triangle-free and which has no separating vertices or edges (i.e. satisfies \hyperref[StandingAssumption1]{Standing  Assumptions 1.1 and 1.2}), $W_\Gamma$ splits over a two-ended subgroup $H$ if and only if $\Gamma$ has a cut collection $\{a-b\}$.

Moreover, given some decomposition $\Lambda$ of $W_\Gamma$ with two-ended edge groups there is a visual decomposition $\Psi$ of $W_\Gamma$ such that:
\begin{itemize}[nosep]
    \item all occurring subgroups in $\Psi$ are special;
    \item each vertex group of $\Psi$ is a subgroup of a conjugate of a vertex group of $\Lambda$;
    \item each edge group of $\Psi$ is a subgroup of a conjugate of an edge group of $\Lambda$;
    \item in particular, for each two-ended edge group $H$ of $\Lambda$ there is a unique cut collection $\{a-b\}$ such that some conjugate of $H$ contains $W_{\{a-b\}}$.
\end{itemize}
\end{thm}

Thus, in order to produce a splitting over two-ended special subgroups, by Theorem \ref{thmvisualsplitting} we need to collect all cut collections of $\Gamma$. Then by Remark \ref{allunivellipiticsubgps} we are left to eliminate the cut collections that produce a subgroup that is not universally elliptic and thus belong inside a hanging subgroup. To be able to do that, we need the following terminology:

\begin{defin}
A cut pair $\{a, b\} \in V(\Gamma)$ is said to be \textit{crossed} by another, disjoint cut pair $\{c,d\}\in V(\Gamma)\setminus\{a,b\}$ if $a$ and $b$ lie in different connected components of $\Gamma \setminus \{c,d\}$. We say $\{c,d\}$ is \textit{crossing} $\{a,b\}$. If there is no cut pair crossing $\{a,b\}$, then $\{a,b\}$ is \textit{uncrossed}.

A cut triple $\{a,b,c\}$, where $c$ is the common adjacent vertex of the non-adjacent vertices $a$ and $b$ is said to be \textit{crossed} by another cut triple $\{d,e,f\}$, where $f$ is the common adjacent vertex of the non-adjacent vertices $d$ and $e$, if $c$ is equal to $f$ and $a$ and $b$ lie in different connected components of $\Gamma \setminus \{d,e,c\}$. We say $\{d,e,f\}$ is \textit{crossing} $\{a,b,c\}$. If there is no cut triple crossing $\{a,b,c\}$, then $\{a,b,c\}$ is \textit{uncrossed}.
\end{defin}

\begin{exam} \label{ExCrossedCutPair}
In Figure \ref{fig:JSJHypCase} of Example \ref{ExamplesK4EssentialCutPair}, while for instance $\{w,z\}$ is an uncrossed cut pair, $\{u, v\}$ is not as it is crossed by $\{l_1, y\}$ for example. In the right graph of Figure \ref{fig:ExampleCutPair} considered in Example \ref{cutcollection} the cut triple $\{a,b,c\}$ is crossed by the cut triple $\{d,e,c\}$.
\end{exam}

\begin{rem}
Any uncrossed cut pair is essential, but not every essential cut pair is uncrossed. Moreover, it is not necessary to define a notion of a crossing between a cut pair and a cut triple because it is obvious that this situation cannot happen.
\end{rem}

It turns out that all the two-ended edge groups of a JSJ splitting are detected by the uncrossed cut collections of $\Gamma$:

\begin{prop}\label{ThmUncrossedCutPairs}
Let $\Gamma$ be a graph which satisfies \hyperref[StandingAssumption1]{Standing Assumption 1} and which has at least one uncrossed cut collection. Then the following hold:
\begin{enumerate}[label=(\alph*)]
    \item For every special subgroup $W_{\{a-b\}}$ generated by an uncrossed cut collection $\{a-b\}$ of $\Gamma$ there is a JSJ splitting $\Lambda$ such that $W_{\{a-b\}}$ is contained in a special, two-ended edge group of $\Lambda$.
    \item Given a two-ended edge group of a JSJ splitting $\Lambda$ of $W_\Gamma$ that is special and contains $W_{\{a-b\}}$ where $\{a-b\}$ is a cut collection, then $\{a-b\}$ is uncrossed.
\end{enumerate}
\end{prop}

\begin{proof}
For $(a)$ let $\{a-b\}$ be a cut collection of $\Gamma$. Let $\Lambda_1$ be a splitting of $W_\Gamma$ over a two-ended subgroup containing $W_{\{a-b\}}$, which exists by Theorem \ref{thmvisualsplitting}. Suppose that $\Lambda_1$ is not a JSJ splitting. If the Bass-Serre tree $T_1$ of $\Lambda_1$ is universally elliptic, but not dominating every other universally elliptic tree, then by \cite[Lemma 2.15]{GuirardelLevittToC}, $T_1$ can be refined to a JSJ tree $T_1^\prime$ with a two-ended edge stabilizer containing $W_{\{a-b\}}$ and the claim follows. If on the other hand the subgroup containing $W_{\{a-b\}}$ is not universally elliptic, there must be another splitting $\Lambda_2$ of $W_\Gamma$ in whose Bass-Serre tree the group $W_{\{a-b\}}$ is not elliptic. Hence, the infinite-order element $ab \in W_{\{a-b\}}$ cannot fix a point in it. Now we can refine $\Lambda_2$ by Theorem \ref{thmvisualsplitting} to a visual splitting $\Psi$ of $W_\Gamma$. Because $\Psi$ is visual, we know that we can find a unique cut collection $\{c-d\}$ in $\Gamma$ such that $W_{\{c-d\}}$ is contained in a two-ended edge group of $\Psi$. Also, every edge group of $\Psi$ is a subgroup of a conjugate of an edge group of $\Lambda_2$ and every vertex group of $\Psi$ is a subgroup of a conjugate of a vertex group of $\Lambda_2$. Thus, the element $ab$ does not fix a point in the Bass-Serre tree of $\Psi$ either, implying that the elements $a$ and $b$ must be in different vertex groups of $\Psi$. This implies that the vertices $a$ and $b$ must be separated in $\Gamma$ by the cut collection $\{c-d\}$, thus the cut collection $\{a-b\}$ is not uncrossed.

Assume conversely for $(b)$ that we have a cut collection $\{a-b\}$ crossed by another cut collection $\{c-d\}$ respectively with splittings $\Lambda_1$ and $\Lambda_2$ over two-ended subgroups containing $W_{\{a-b\}}$ and $W_{\{c-d\}}$ respectively. Since the cut collection $\{a-b\}$ is crossed by $\{c-d\}$, the elements $a$ and $b$ are in different vertex groups of the splitting $\Lambda_2$. Thus the infinite-order element $ab \in W_{\{a-b\}}$ cannot fix a point in the Bass-Serre tree of $\Lambda_2$, implying that the edge group of $\Lambda_1$ containing $W_{\{a-b\}}$ is not universally elliptic and therefore $\Lambda_1$ is not a JSJ decomposition.
\end{proof}

\begin{rem}
In Proposition \ref{ThmUncrossedCutPairs} the assumption that $\Gamma$ must contain an uncrossed cut collection excludes the case where $\Gamma$ is a square. This is important, because for $\Gamma$ a square, the corresponding RACG $W_\Gamma = D_\infty \times D_\infty$ is virtually $\mathbb{Z}^2$. Such $W_\Gamma$ is commensurable to the fundamental group of a surface, in this case a torus, which have to be treated separately \cite[cf.][]{papasoglu2005quasi}. However, this is the only case we need to rule out additionally since by the \hyperref[StandingAssumption1]{Standing Assumption 1}, $W_\Gamma$ is not cocompact Fuchsian and thus never commensurable to a surface group of higher genus.

Also, excluding the case that $\Gamma$ is not a square is no obstacle for the QI-classification, since the property of being virtually $\mathbb{Z}^2$ determines the QI-type of the group. Thus we refine the standing assumption by modifying \hyperref[StandingAssumption1Part4]{(4)} of \hyperref[StandingAssumption1]{Standing Assumption 1} to exclude squares:
\end{rem}

\begin{Assump} \label{StandingAssumption2}
The defining graph $\Gamma$
\begin{enumerate}[label=(\arabic*)]
    \item has no triangles.
    \item is connected and has neither a separating vertex nor a separating edge.
    \item has a cut collection. 
    \item is not a cycle.
\end{enumerate}
\end{Assump}

Now, starting from a visual JSJ decomposition over all uncrossed cut collections, we can determine how to produce the different vertices and the edges of the JSJ graph of cylinders.

\subsection{Cylinder vertices} \label{SubsectionCylVertices}
Given a JSJ decomposition $\Lambda$, by Outline \ref{ExFramework} we know that we can find all cylinder vertices of the JSJ graph of groups $\Lambda_c$ and their vertex groups $G_Y$ by running through all edge groups visible in $\Lambda$.

Thus in light of Proposition \ref{ThmUncrossedCutPairs}, we pick up all uncrossed cut collections in $\Gamma$ and compute the commensurators of the special subgroups they generate. It turns out that the commensurator of a special subgroup is also visible from the defining graph $\Gamma$:

\begin{thm} \textup{\cite[Theorem 2.1]{Paris1997commensurators}} \label{commensurator} Let W be a RACG on the defining graph $\Gamma$ with finite generating set $S$ and let $T \subseteq S$  be a subset of $S$. Consider the maximal decomposition $W_T = W_{T_1} \times \, \cdots \, \times W_{T_n}$ of $W_T$ as a direct product of subgroups, where $W_{T_1}, \,  \dots \,, W_{T_r}$ are finite for some $r \in \{1, \, \dots \, , n\}$ and all the other subgroups are infinite. Then the commensurator of $W_T$ in W is given by
\begin{gather*}
    Comm_{\,W}(W_T) = W_{T^\infty} \times W_{Y^\infty} \\[0.3cm]
    \text{with} \quad T^\infty = T_{r+1} \cup \dots \cup T_n \quad
    \text{and} \quad Y^\infty = \{s \in S \mid e = (t,s) \in E(\Gamma) \textit{ for all } t \in T^\infty\} \, .
\end{gather*}
\end{thm}

\begin{nota}
To simplify terminology we refer to the vertices of the defining graph of the commensurator of the special subgroup given by a cut collection as \textit{commensurator of the cut collection}.
\end{nota}

\begin{rem} \label{IllustrationCommensurator}
We encounter the following situations:
\begin{itemize}
    \item For a cut pair $\{a,b\}$ in $\Gamma$ this means that the commensurator is generated by $\{a,b\} \cup \mathcal{C}$, where $\mathcal{C}$ contains all the common adjacent vertices of $a$ and $b$. That is $$Comm_{\,W}(W_{\{a,b\}}) = W_{\{a,b\}} \times W_{\mathcal{C}} \, .$$
    \item For a cut triple $\{a,b,c\}$, where $a$ and $b$ are non-adjacent and $c$ is a common adjacent vertex, the special subgroup $W_{\{a,b,c\}}$ decomposes as $W_{\{a,b\}} \times W_{\{c\}}$. Thus also in this case, the commensurator is generated by $\{a,b\} \cup \mathcal{C}$, where $\mathcal{C}$ contains all the common adjacent vertices of $a$ and $b$, in particular $\mathcal{C}$ contains $c$.
    \item In case there are two overlapping cut triples $\{a,b,c\}$ and $\{a,b,c^\prime\}$ sharing the same non-adjacent vertex pair $\{a,b\}$, for both cut triples we obtain the same commensurator. Hence their corresponding edges in a JSJ decomposition are equivalent under commensurability, thus they lie in the same cylinder. Therefore such two cut triples only give one cylinder.
    \item It is immediate from Theorem \ref{commensurator} that a cut collection $\{a-b\}$ in a hyperbolic RACG always has a two-ended commensurator. This is because $a$ and $b$ can have at most one common adjacent vertex $c$, as otherwise two common adjacent vertices and $a$ and $b$ give a square in contradiction to hyperbolicity. But both $W_{\{a,b\}} \simeq D_\infty$ and $W_{\{a,b,c\}} \simeq D_\infty \times \mathbb{Z}_2$ are two-ended (cf.\ Theorem \ref{VCylic}), thus hyperbolic RACGs have two-ended cylinder vertices.
\end{itemize}
\end{rem}

\begin{exam}
In the non-hyperbolic defining graph $\Gamma$ in Figure \ref{fig:ExampleNonHypGraph} the orange cut pair $\{v,x\}$ has three purple common adjacent vertices $\mathcal{C} = \{w,m,y\}$, thus $$Comm_{\, W_\Gamma}(W_{\{v,x\}}) = W_{\{v,x\}} \times W_{\{w,m,y\}} = W_{\{v,w,m,y,x\}} \, .$$
The other two cut pairs $\{w,z\}$ and $\{y,z\}$ correspond to special subgroups with commensurators $W_{\{w,z,n,x\}}$ and $W_{\{y,z,o,x\}}$ respectively.

The commensurator of the special subgroup corresponding to the cut triple $\{w,x,y\}$ on the other hand is $W_{\{v,w,x,y\}}$ since $v$ and $x$ are the common adjacent vertices of $w$ and $y$. This is the only cut triple in $\Gamma$: Recall that for instance the vertices $\{w,m,y\}$ are not a cut triple, despite separating $v$ from the the rest of $\Gamma$, because $W_{\{w,m,y\}}$ is not two-ended.

\begin{figure}[ht]
    \centering
    \scriptsize
    \def\svgwidth{270pt}
    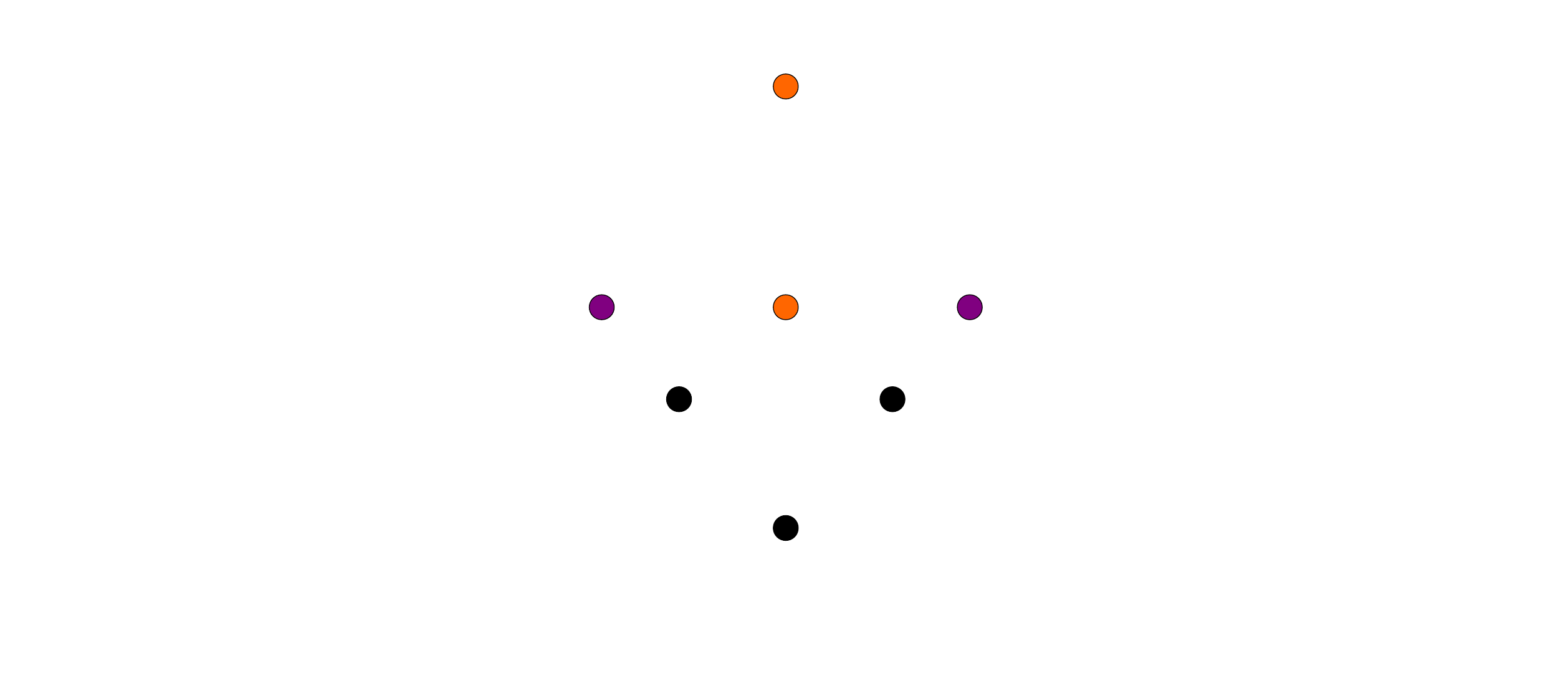
    \caption{The special subgroup generated by the orange cut pair $\{v,x\}$ has three purple common adjacent vertices $\{w,m,y\}$. Thus its commensurator is generated by $\{v,w,m,y,x\}$.}
   \label{fig:ExampleNonHypGraph}
\end{figure}
\end{exam}

For the sake of completeness we summarize this insight as a Proposition:

\begin{prop} \label{PropUncrossedCutPairsCylVert}
Let $\mathcal{S}$ be the following set: For every uncrossed cut collection $\{a-b\}$ of $\Gamma$, the set $\{a,b\} \cup \mathcal{C}$, where $\mathcal{C}$ is the set of common adjacent vertices of $a$ and $b$, is an element in $\mathcal{S}.$ Then every set $S$ in $\mathcal{S}$ corresponds to a cylinder vertex in the JSJ graph of cylinders of $W_\Gamma$ with vertex group the special subgroup generated by $S$.
\end{prop}

\begin{lem} \label{LemTypesofCylinderGroups}
Every cylinder vertex group of the JSJ graph of cylinders of a RACG $W_\Gamma$, where $\Gamma$ satisfies the \hyperref[StandingAssumption2]{Standing Assumption 2}, in particular is triangle-free, is either
\begin{itemize}[noitemsep,topsep=2pt]
    \item virtually cyclic;
    \item virtually $\mathbb{Z}^2$; or
    \item the direct product of a virtually non-abelian free group and the infinite dihedral group $D_\infty$.
\end{itemize}
\end{lem}

For the proof we need the following characterizations:

\begin{thm} \textup{\cite[Theorem 17.2.1]{DavisCoxeter}} \label{VAbelian}
A RACG $W_\Gamma$ is virtually abelian if and only if it decomposes as the direct product of finitely many infinite dihedral groups $D_\infty$ and a finite RACG.
\end{thm}

\begin{thm} \textup{\cite[Theorem 8.34]{MihalikTschantzVisual}} \label{VFree} A RACG $W_\Gamma$ is virtually free if and only if no induced subgraph is a circuit of more than three vertices.
\end{thm}

We can detect the intersection of the above two classes of groups by the following theorem:

\begin{thm} \cite[Theorem 8.7.3]{DavisCoxeter} \label{VCylic} A RACG $W_\Gamma$ is two-ended if and only if it is the direct product of one infinite dihedral group $D_\infty$ and a finite RACG. In terms of the defining graph this means that $\Gamma$ is a two-point suspension of a complete graph.
\end{thm}

\begin{proof} \textit{(of Lemma \ref{LemTypesofCylinderGroups})}
Consider a cut collection $\{a-b\}$ of $\Gamma$ with $a$ and $b$ non-adjacent, then $\Gamma \setminus \{a-b\}$ has $i \geq 2$ connected components $\Gamma_1, \, \dots \, , \Gamma_i$. We distinguish the contribution of some component $\Gamma_j$ for $j \in \{1,\dots,i\}$ to the commensurator $G_Y$ of $W_{\{a-b\}}$:
\begin{itemize}[noitemsep,topsep=2pt]
    \item If $\Gamma_j$ does not contain any common adjacent vertex of $\{a,b\}$, no vertex contributes to $G_Y$.
    \item If $\Gamma_j$ contains one common adjacent vertex $c$ of $\{a,b\}$, the contribution to $G_Y$ is a direct product with $\mathbb{Z}_2$.
    \item If $\Gamma_j$ contains at least two common adjacent vertices $c_1$ and $ c_2$ of $\{a,b\}$, then they must be connected by a path not passing through $a$ or $b$. Otherwise they would not lie in the same connected component of $\Gamma \setminus \{a-b\}$. However, there cannot be an edge between $c_1$ and $c_2$, as otherwise $\{a,c_1,c_2\}$ would form a triangle. Thus there is no relation between $c_1$ and $c_2$ in $G_Y$.
\end{itemize}
Moreover, if $\{a-b\}$ is a cut triple, the third vertex $c$ of the triple contributes a direct product with $\mathbb{Z}_2$ to $G_Y$. In conclusion the commensurator $G_Y$ of the cut collection $\{a-b\}$ is the RACG given by a defining graph $\Gamma_Y$ consisting of $a$ and $b$ with $k$ common adjacent vertices $\{c_1, \, \dots \, , c_k\} =: \mathcal{C}$, which are all only connected to $a$ and $b$ in $G_Y$, see Figure \ref{fig:PictureCommensurator}.
Thus we can consider the following cases:
\begin{itemize}[noitemsep,topsep=2pt]
    \item $\mathcal{C} = \{\}$: $G_Y = W_{\{a,b\}} \simeq D_\infty$, thus virtually cyclic.
    \item $k = 1$: $G_Y = W_{\{a,b, c_1\}} \simeq D_\infty \times \mathbb{Z}_2$, thus virtually cyclic.
    \item $k = 2$: $G_Y = W_{\{a,b, c_1, c_2\}} \simeq D_\infty \times D_\infty$, thus virtually abelian, in particular virtually $\mathbb{Z}^2$.
    \item $k > 2$: $G_Y = W_{\{a,b, c_1, \, \dots \, , c_k \}} \simeq D_\infty \times F$, where $F$ is virtually a non-abelian free group.
\end{itemize}

\begin{figure}[ht]
    \centering
    \scriptsize
    \def\svgwidth{140pt}
    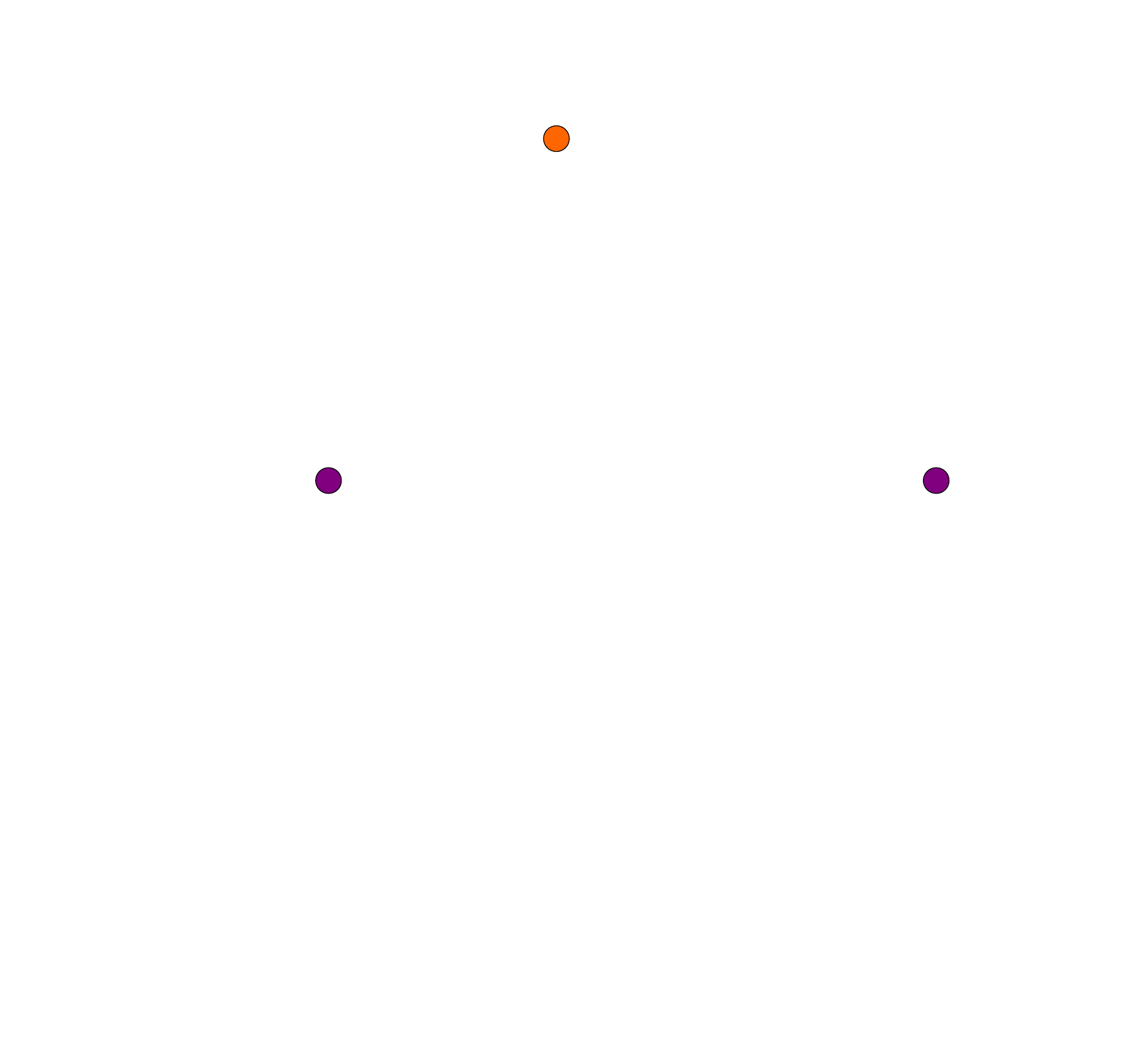
    \caption{The common adjacent vertices $\{c_1, \dots, c_k\}$ of $a$ and $b$ are only adjacent to both $a$ and $b$.}
   \label{fig:PictureCommensurator}
\end{figure} \vspace{-0.4cm}
\end{proof}

\begin{exam}
The commensurator of $W_{\{v,x\}}$ in Figure \ref{fig:ExampleNonHypGraph} is the direct product of the virtually non-abelian free group $W_{\{w,m,y\}}$ and the infinite dihedral group $W_{\{v,x\}}$.
\end{exam}

\begin{nota}
From now on we refer to the two "new" types of cylinder vertex groups and their corresponding cylinder vertex as
\begin{itemize}[noitemsep,topsep=2pt]
    \item \textit{VA}, if the cylinder vertex group is virtually $\mathbb{Z}^2$.
    \item \textit{VFD}, if the cylinder vertex group is the direct product of a virtually non-abelian free group and an infinite dihedral group.
\end{itemize}
\end{nota}

\subsection{Non-cylinder vertices} \label{SubsectionNonCylVert}
The fact that a certain collection of vertices gives a hanging or rigid vertex group in a graph of groups with respect to incident two-ended edge groups is intrinsic to this collection in the sense that it is independent of the existence of squares in the defining graph $\Gamma$. Furthermore, by Outline \ref{ExFramework}, if we see a hanging or rigid vertex in the JSJ decomposition, it transfers over to the JSJ graph of cylinders. So, the results of \cite{DaniThomasBowditch} in Theorem \ref{ToChypRACG} translate to the general setting:

\begin{prop} \label{Charhanging}
Let $A \subseteq V(\Gamma)$ be a set of vertices such that the $A$-induced subgraph $\Gamma_A$ is not a complete graph and $A$ satisfies either the following conditions:
    \begin{itemize}[leftmargin=1.5cm]
        \item[(A1)] \label{A1n} Elements of $A$ pairwise separate the geometric realization $|\Gamma|$. 
        \item[(A2)] \label{A2n} If any subgraph $\Gamma^\prime$ of $\Gamma$ that is a subdivided $K_4$ contains more than 2 vertices of $A$, all vertices of $A$ lie on the same branch of $\Gamma^\prime$.
        \item[(A3)] \label{A3n} The set $A$ is maximal among all sets satisfying \hyperref[A1n]{(A1)} and \hyperref[A2n]{(A2)}.
    \end{itemize} 
Or $A$ satisfies the following condition:
    \begin{itemize}[leftmargin=1.5cm]
        \item[(A*)] \label{A*} The set $A$ is a maximal collection of pairwise crossing cut triples.
    \end{itemize}
If $A$ is not contained in a vertex set corresponding to a cylinder vertex, then $A$ corresponds to a vertex in the JSJ graph of cylinders $\Lambda_c$ with hanging vertex group $W_A$.
\end{prop}

\begin{sproof}
Recall that the assumption that $\Gamma_A$ is not a complete graph ensures that $W_A$ is infinite. Now we want to give some motivation on how the graph theoretic conditions \hyperref[A1n]{(A1)}, \hyperref[A2n]{(A2)} and \hyperref[A3n]{(A3)} and the graph theoretic condition \hyperref[A*]{(A*)} produce a hanging subgroup $W_A$. The intuitive picture to have in mind as the hanging subgroup is a surface with boundary.

Let us first consider crossing cut pairs. By the proof of Proposition \ref{ThmUncrossedCutPairs}, they are not universally elliptic and thus belong in a hanging subgroup. They give crossing curves corresponding to different interfering splittings, which are thus not part of a JSJ decomposition. In particular, a collection of pairwise separating vertices as forced by condition \hyperref[A1]{(A1)} contains all pairwise crossing cut pairs within a branch and at least one uncrossed essential cut pair (cf. Remark \ref{hangingUncrossedEssCutPair}). Such an uncrossed cut pair then generates precisely a universally elliptic subgroup as the boundary component.
If we see however a subdivided $K_4$ in $\Gamma$, we could choose three or all four corner vertices as a collection of pairwise $|\Gamma|$-separating vertices. But then this collection cannot contain any non-essential vertex contained in a branch. So there are no crossing cut pairs in the collection producing crossing curves. Thus the resulting group is not a hanging, but rather a candidate for a rigid vertex group. Therefore we need to exclude such a collection by condition \hyperref[A2]{(A2)}. Maximality needs to be ensured by condition \hyperref[A3]{(A3)}, since a JSJ decomposition is maximal (cf.\ Definition \ref{defHanging} and Lemma \ref{EqDefJSJdecomp}).

\smallskip
Consider now crossing cut triples contained in a collection $A$ satisfying condition \hyperref[A*]{(A*)}. Again, by the proof of Proposition \ref{ThmUncrossedCutPairs}, they are not universally elliptic and thus belong in a hanging subgroup. By the definition of a JSJ decomposition, we again need maximality.

Since all the cut triples in $A$ cross pairwise, all share their "middle" common adjacent vertex, call it $c$. Thus the subgraph induced on the collection $A$ is a graph theoretical star based at $c$. Since by \hyperref[StandingAssumption2]{Standing Assumption 2}, $\Gamma$ has no triangles and no separating edge, for every pair $\{x,y\} \in A\setminus\{c\}$ of leaves, $x$ and $y$ are not adjacent and there are at least three disjoint paths connecting $x$ and $y$: One is the segment $\{x,y,c\}$ and two paths do not contain $c$, call them $p_1$ and $p_2$. 

We claim that either $\{x,y\}$ is an uncrossed cut pair or $\{x,y,c\}$ is a cut triple: If every path connecting the interior of $p_1$ (or analogously $p_2$) with $c$ passes through $x$ or $y$, removing $\{x,y\}$ separates the interior of $p_1$ from $\Gamma_A \setminus \{x,y\}$. Hence $\{x,y\}$ is a cut pair. In fact $\{x,y\}$ is an uncrossed cut pair, because no matter which other cut pair is removed, $x$ and $y$ will stay connected with each other via either $p_1$, $p_2$ or the segment $\{x,y,c\}$. Thus the cut pair $\{x,y\}$ generates a universally elliptic subgroup representing the boundary component of the surface.

If both $p_1$ and $p_2$ contain an interior vertex that is connected to $c$ via a path not passing through $x$ or $y$, then we need to show that removing $\{x,y,x\}$ separates $\Gamma$. Since $x$ is leaf in $A \setminus \{c\}$, there exists $x^\prime \in A\setminus\{c,x\}$ such that $\{x,x^\prime,c \}$ is a cut triple separating $\Gamma$ into two connected components $C$ and $C^\prime$ of $\Gamma \setminus \{x,x^\prime,c \}$. Then either the interior of $p_1$ or $p_2$ must be contained in one of the connected components of $\Gamma \setminus \{x,x^\prime,c \}$. Without loss of generality assume that the interior of $p_1$ is contained in $C$. Thus there is a vertex $l_1$ in $V(C) \cap V(p_1)$ which is not connected in $\Gamma \setminus \{x,x^\prime,c \}$ to some $l_2 \in V(C^\prime)$. However, since $l_1$ is in the interior of $p_1$, $l_1$ will also not be connected to $l_2$ in  $\Gamma \setminus \{x,y,c\}$. Thus also $\{x,y,c\}$ is a cut triple contained in $A$. If it is uncrossed it represents a boundary component.
\end{sproof}

\begin{cor} \label{HangStar}
    Let $A \subseteq V(\Gamma)$ be a set of vertices satisfying the condition \hyperref[A*]{(A*)}. Then the $A$-induced subgraph $\Gamma_A$ of $\Gamma$ is a star.
\end{cor}

\begin{prop} \label{charrig}
For any set $B \subseteq EV(\Gamma)$ of essential vertices in $\Gamma$ satisfying the properties \hyperref[B1n]{(B1)}, \hyperref[B2n]{(B2)} and \hyperref[B3n]{(B3)}, there is a vertex in the JSJ graph of cylinders $\Lambda_c$ with rigid vertex group $W_B$, where the properties \hyperref[B1n]{(B1)}, \hyperref[B2n]{(B2)} and \hyperref[B3n]{(B3)} are the following:
        \begin{itemize}[leftmargin=1.5cm]
            \item[(B1)] For any set $C$ that is a pair $\{c,d\} \subseteq EV(\Gamma)$ of essential vertices or a path $\{c,d,e\} \subseteq EV(\Gamma)$ of length 2 of essential vertices, $B \setminus C$ is contained in one single connected component of $\Gamma \setminus C$. \label{B1n}
            \item[(B2)] The set $B$ is maximal among all sets satisfying \hyperref[B1n]{(B1)}. \label{B2n}
            \item[(B3)] $|B| \geq 4$. \label{B3n}
        \end{itemize}
\end{prop}

\begin{sproof}
Again, we want to give some motivation on how the graph theoretic conditions \hyperref[B1n]{(B1)}, \hyperref[B2n]{(B2)} and \hyperref[B3n]{(B3)} produce a rigid subgroup $W_B$.
The key feature of a rigid vertex group is that it cannot be split any further. This is precisely captured by condition \hyperref[B1n]{(B1)}: We consider the collection of cut pairs and cut triples which are pairwise not separating the collection. We want a maximal such collection and thus impose condition \hyperref[B2n]{(B2)}. Suppose now we find a collection $B =\{x,y,z\}$ with only three essential vertices satisfying conditions \hyperref[B1n]{(B1)} and \hyperref[B2n]{(B2)}. Then, since we restrict to special subgroups, the virtually free RACG $W_B$ can have the adjacent edge groups $W_{\{x,y\}}$, $W_{\{y,z\}}$ and $W_{\{x,z\}}$. However, such a group is then virtually a surface with boundary and thus not considered to be rigid. This case is excluded by condition \hyperref[B3n]{(B3)}.
\end{sproof}

\subsection{Edges} \label{SubsectionEdges}
It remains to be determined which vertices in the JSJ graph of cylinders are connected by an edge:
\begin{lem} \label{EdgeIffIntersection}
For any pair of vertices in the JSJ graph of cylinders $\Lambda_c$ there is an edge connecting them if and only if the pair consists of one cylinder vertex corresponding to the cut collection $\{a-b\}$ and one non-cylinder vertex and their vertex groups intersect in a special subgroup containing $W_{\{a-b\}}$. The edge group is the special subgroup generated by the intersection of the corresponding vertex sets.
\end{lem}

\begin{proof} Since the JSJ graph of cylinders $\Lambda_c$ is bipartite, edges can only connect cylinder with non-cylinder vertices. Suppose there is an edge connecting a cylinder vertex corresponding to a cut collection $\{a-b\}$ and a non-cylinder vertex. By definition of the fundamental group of a graph of groups the edge group is the intersection of their vertex groups. If this intersection would be a finite group, the group $W$ cannot be one-ended by Stallings' Theorem \cite{Sta71}, in contradiction to the \hyperref[StandingAssumption2]{Standing Assumption 2}. Thus the intersection is infinite. Furthermore, the vertex groups are special subgroups, thus so is their intersection by \cite[Theorem 4.1.6]{DavisCoxeter}. Since the structure of $\Lambda_c$ comes from a JSJ decomposition $\Lambda$, by Proposition \ref{ThmUncrossedCutPairs} the edge group in $\Lambda$ must contain $W_{\{a-b\}}$, thus so does the edge group in $\Lambda_c$.

Assume conversely that the vertex group of a cylinder vertex corresponding to a cut collection $\{a-b\}$ and a non-cylinder vertex intersect in a special subgroup containing $W_{\{a-b\}}$. Then they are connected by an edge by the definition of the fundamental group of a graph of groups.
\end{proof}

\begin{exam} \label{ExampleOverlap}
For the graph $\Gamma$ shown in Figure \ref{fig:Overlap}, which satisfies \hyperref[StandingAssumption2]{Standing Assumption 2}, we can construct the corresponding JSJ graph of cylinders in Figure \ref{fig:OverlapGoC} by reading off the following collections of vertices according to Proposition \ref{PropUncrossedCutPairsCylVert}, Theorem \ref{Charhanging}, Theorem \ref{charrig} and connect them with edges according to Lemma \ref{EdgeIffIntersection}:

\begin{center}
\begin{tabular}[t]{c|c||c||c}
   uncrossed cut collection & commensurator & hanging & rigid \\\hline
   $x,w$ & $x,w,k,d$ & & \\
   $v,w$ & $v,w,d$ & $v,w,l_1,l_2$ \\
   $v,y$ & $v,y,m,d$ & & \\
   $x,y,b$ | $x,y,c$ | $x,y,d$ & $x,y,a,b,c,d$ & $v,w,x,y,d$ & \\
   $c,d$ & $c,d,x,y$ & $c,d,n_1,n_2$ & $c,x,d,y$ \\
   $b,c$ & $b,c,x,o,y$ & & $b,x,c,y$ \\
   $a,b$ & $a,b,x,y$ & $a,b,p_1,p_2$ & $a,x,b,y$
\end{tabular}
\end{center}

\begin{figure}[ht]
    \centering
    \scriptsize
    \def\svgwidth{440pt}
    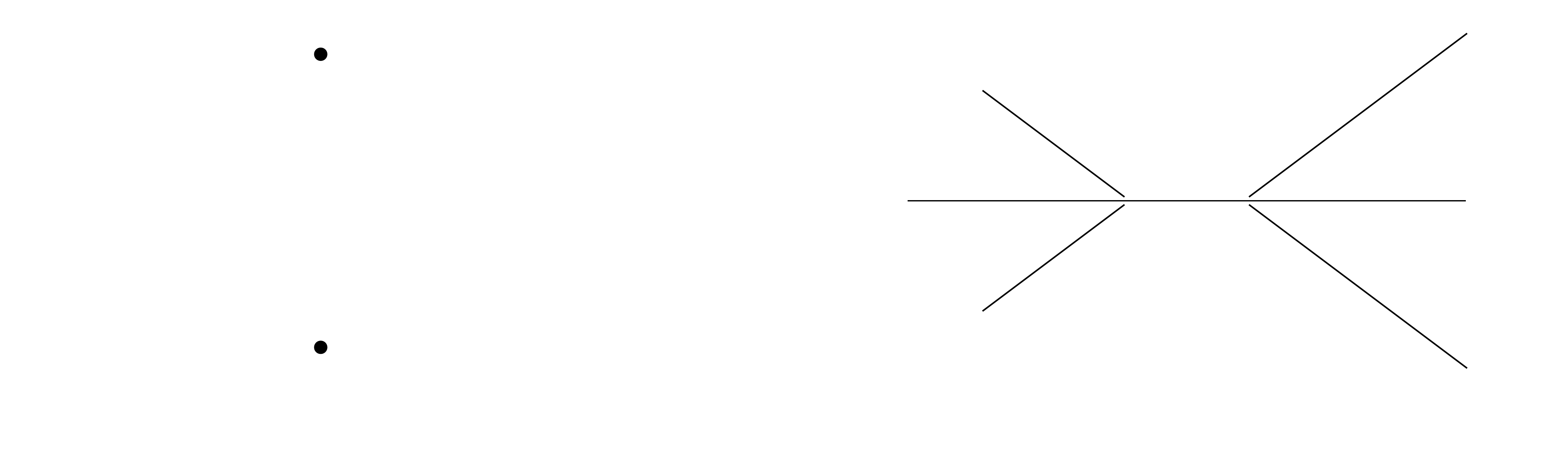
    \captionsetup{belowskip=-0.6cm}
    \caption{$\Lambda_c$ is the JSJ graph of cylinders of the RACG $W_\Gamma$.}
    \label{fig:OverlapGoC}
\end{figure}
\end{exam}

\subsubsection{Two-ended edge groups} \label{TwoEndedEdgeGroups}

As indicated in Remark \ref{RemCollapsedToC}, we aim to restrict to trees of cylinders that are $\mathcal{VC}$-trees themselves. This is not always the case, as we can see in  Example \ref{ExampleOverlap}, Figure \ref{fig:OverlapGoC}: The vertex set $\{x,y,a,b,c,d\}$ generating the commensurator of the uncrossed cut triples $\{x,y,b\}$, $\{x,y,c\}$ and $\{x,y,d\}$ contains for instance the collection $\{c,x,d,y\}$, which corresponds to an adjacent rigid vertex. Thus, the connecting edge group generated by the intersection by Lemma \ref{EdgeIffIntersection} is $W_{\{c,x,d,y\}} = D_\infty \times D_\infty$, which is not two-ended.


Therefore we need to impose assumptions on the defining graph $\Gamma$ to ensure that the intersection of vertex groups is two-ended. Recall that by the bipartiteness of the JSJ tree of cylinders, the only intersections we need to consider are between cylinder and non-cylinder vertices.

\begin{lem} \label{HangCylEdges2Ended}
If the intersection of a cylinder vertex group $G_Y$ corresponding to a cut collection $\{a-b\}$ and a hanging vertex group $W_A$ contains $W_{\{a-b\}}$ and
\begin{enumerate}[label=(\alph*)]
    \item the set $A$ satisfies the conditions \hyperref[A1n]{(A1)}, \hyperref[A2n]{(A2)} and \hyperref[A3n]{(A3)}, then the intersection is the infinite dihedral group $D_\infty$. \label{HangCylEdges2EndedCase1}
    \item the set $A$ satisfies the condition \hyperref[A*]{(A*)}, then the intersection is two-ended. \label{HangCylEdges2EndedCase2}
\end{enumerate} 
\end{lem}

The proof of Lemma \ref{HangCylEdges2Ended} relies on the following properties:

\begin{lem} \label{LemVerticesInA} If $A \subseteq V(\Gamma)$ is a set that satisfies conditions \hyperref[A1n]{(A1)}, \hyperref[A2n]{(A2)} and \hyperref[A3n]{(A3)}, the $A$-induced subgraph $\Gamma_A$ is not a complete graph and $W_A$ is not contained in a cylinder vertex group, then:
\begin{enumerate}[nosep, label=(\arabic*)]
\item $A$ does not contain a cut triple.
\item $A$ does not contain two branches which share a common endpoint.
\end{enumerate}
\end{lem}

\begin{proof} 
By definition, the non-adjacent vertices $a$ and $b$ of a cut triple $\{a-b\}$ are not a cut pair and by \hyperref[StandingAssumption2]{Standing Assumption 2} $a$ and $b$ do not share an edge. Thus $a$ and $b$ do not separate $|\Gamma|$. Therefore a set satisfying \hyperref[A1n]{(A1)} cannot contain a cut triple, implying (1). 

For (2) suppose that two vertices $x$ and $y$ of degree 2 lie in different branches contained in $A$ meeting at an essential vertex $a$. Let $b_x$ and $b_y$ be the other endpoint of these branches respectively. Since $a$, $b_x$ and $b_y$ are essential, $a$ is connected to both $b_x$ and $b_y$ via a path neither passing through $x$ nor $y$. Thus $|\Gamma| \setminus \{x,y\}$ is connected, in contradiction to condition \hyperref[A1]{(A1)}.
\end{proof}

\begin{proof} \textit{(of Lemma \ref{HangCylEdges2Ended})} For the proof of \hyperref[HangCylEdges2EndedCase1]{(a)} assume that $A$ satisfies the conditions \hyperref[A1n]{(A1)}, \hyperref[A2n]{(A2)} and \hyperref[A3n]{(A3)}. Let $W_A$ be the corresponding hanging vertex group on the defining graph $\Gamma_A$ intersecting the cylinder vertex group $G_Y$ corresponding to the cut collection $\{a-b\}$ non-trivially. Recall that by Lemma \ref{EdgeIffIntersection} the intersection must contain $W_{\{a-b\}}$. Then by Lemma \ref{LemVerticesInA}.1. $\{a-b\}$ cannot be a cut triple, so $G_Y$ must correspond to a cut pair $\{a,b\}$.

If $G_Y$ is $W_{\{a,b\}} = D_\infty$, so is the intersection. Therefore we can assume that $G_Y$ is not $D_\infty$. Thus the cylinder vertex group is the special subgroup on the defining graph $\Gamma_Y$ consisting of the pair $\{a,b\}$ with a non-empty common adjacent vertex set $\{c_1, c_2, \dots, c_k\}$ for $k \geq 1$ and the degree of every vertex in $\mathcal{C}$ in $\Gamma_Y$ is 2. Since by \cite[Theorem 4.1.6]{DavisCoxeter}, the intersection of two special subgroups is the special subgroup defined on the induced graph given by the intersection, we need to consider how $\Gamma_A \cap \Gamma_Y$ can look like. 
Recall that the intersection $I = V(\Gamma_A \cap \Gamma_Y)$ contains $a$ and $b$. We distinguish the following cases:
\begin{enumerate}
    \item $I = \{a,b\}$: The corresponding group $W_{\{a,b\}} \cong D_\infty$ is two-ended.
    \item $c_i \in I$ and $c_i$ has degree 2 in $\Gamma$ for $i \in \{1,\dots,k\}$: Then $I$ contains the whole branch $\{a,c_i,b\}$. No other branch in $A$ can be attached at $a$ or $b$ by Lemma \ref{LemVerticesInA}.2., implying $I = \{a,c_i,b\}$. But $A$ cannot be equal to $I = \{a,c_i,b\}$, since the hanging vertex corresponding to $A$ is not two-ended. However, supposing that there is another vertex $v \in A\setminus \{a,c_i,b\}$ such that $|\Gamma| \setminus \{c_i,v\}$ is separated, contradicts the fact that $a$ and $b$ are uncrossed and $c_i$ has degree 2. This implies that this case cannot occur.
    \item $c_i \in I$ for $i \in \{1,\dots,k\}$ and $c_i$ essential in $\Gamma$: We argue as in case 2.\ that there must exist a $v \in A\setminus \{a,c_i,b\}$ such that $|\Gamma| \setminus \{c_i,v\}$ is separated. Since $a$ and $b$ are an uncrossed cut pair, there is a path between $v$ and $c_i$ not passing through $a$ nor $b$ and another path connecting $a$ and $b$ not passing through $c_i$ nor $v$. This means that we have a subdivided $K_4$ with corners $a,c_i,b$ and $v$, in contradiction to \hyperref[A2]{(A2)}. So again, this case cannot occur.
\end{enumerate}
To conclude, in case \hyperref[HangCylEdges2EndedCase1]{(a)} the special subgroup $W_I$ generated by the intersection $I$ is always $D_\infty$.
\smallskip

Assume now for the proof of \hyperref[HangCylEdges2EndedCase2]{(b)} that $A$ satisfies the condition \hyperref[A*]{(A*)}, and that the corresponding vertex group $W_A$ on the defining graph $\Gamma_A$ is infinite. By Corollary \ref{HangStar} $\Gamma_A$ is a graph theoretical star based at the vertex $c$ where all the cut triples contained in $A$ meet. Suppose that $W_A$ intersects the cylinder vertex group $G_Y$ corresponding to a cut collection $\{a-b\}$ in a subgroup containing $W_{\{a-b\}}$. Thus, if $G_Y$ is two-ended, so is the intersection and we can assume that $G_Y$ is not two-ended. That means that the cylinder vertex group is the special subgroup on the defining graph $\Gamma_Y$ consisting of the two non-adjacent vertices $\{a,b\}$ of the cut collection and their common adjacent vertices $\mathcal{C}=\{c_1, \dots, c_k\}$ with $k\geq 2$, which all have degree 2 in $\Gamma_Y$. As above, we need to consider the graph $\Gamma_A \cap \Gamma_Y$. Define $I = V(\Gamma_A) \cap V(\Gamma_Y)$ which contains $\{a-b\}$ by Lemma \ref{EdgeIffIntersection} and consider the following cases for $I$:
\begin{enumerate}
    \item $I = \{a,b\}$: In this case $W_I$ is $D_\infty$ thus two-ended.
    \item $I = \{a,b,c_i\}$ for some $i \in \{1, \dots, k\}$: In this case $W_I$ is $D_\infty \times \mathbb{Z}_2$ and thus two-ended.
    \item $\{a,b\} \in I$ and $|I \cap \mathcal{C}| \geq 2$: Then $I$ contains a square, thus so does $\Gamma_A$ in contradiction to the fact that $\Gamma_A$ is a triangle-free star. Thus this case cannot occur.
\end{enumerate}
In conclusion, also in case \hyperref[HangCylEdges2EndedCase2]{(b)} the special subgroup $W_I$ generated by the intersection $I$ is always two-ended. This finishes the proof.
\end{proof}

\begin{rem} \label{hangingUncrossedEssCutPair} By the \hyperref[StandingAssumption2]{Standing Assumption 2} the defining graph $\Gamma$ is never a cycle. Thus in case $A$ satisfies the conditions \hyperref[A1n]{(A1)}, \hyperref[A2n]{(A2)} and \hyperref[A3n]{(A3)} by Lemma \ref{LemVerticesInA}.2, $A$ cannot contain a cycle. This is also true in case \hyperref[HangCylEdges2EndedCase2]{(b)}, where $A$ satisfies $\hyperref[A*]{(A*)}$ since the $A$-induced subgraph $\Gamma_A$ is a graph theoretical star by Corollary \ref{HangStar}. Therefore Theorem \ref{VFree} implies that any hanging vertex group is virtually free.
\end{rem}

\begin{thm} \label{ThmK4Assump}
Let $G_Y$ be the vertex group of the cylinder vertex $v_Y$ in $\Lambda_c$ corresponding to the cut collection $\{a-b\}$ with defining graph $\Gamma_Y \subseteq \Gamma$ on the vertex set $V(\Gamma_Y) = \{a,b\} \cup \mathcal{C}$, where $\mathcal{C}$ is the set of common adjacent vertices of $a$ and $b$. Then every rigid vertex group $W_B$ adjacent to the cylinder vertex group $G_Y$ intersects $G_Y$ in a two-ended subgroup if and only if for any pair of vertices in $\mathcal{C}$ every path connecting them in $\Gamma$ passes through $a$ or $b$.
\end{thm}

\begin{exam}
In Figure \ref{fig:OverlapGoC}, the rigid vertex group generated by $\{c,d,x,y\}$ is adjacent to the cylinder vertex group corresponding to the cut triple $\{x-y\}$, which has $\{c,d\}$ as common adjacent vertices. Because there is a path through the vertices $\{c,n_1,n_2,d\}$ connecting $c$ and $d$ without passing through $x$ nor $y$, they intersect in the non-two-ended edge group generated by $\{c,d,x,y\}$.
\end{exam}

\begin{proof} Suppose first that there is a pair $\{c_i,c_j\} \subseteq \mathcal{C}$ of distinct vertices that are connected by a path in $\Gamma$ not passing through $a$ nor $b$ nor any other common adjacent vertex of $a$ and $b$. There must be a path between $a$ and $b$ not passing through $c_i$ nor $c_j$ as otherwise $\{a-b\}$ would be crossed by $\{c_i-c_j\}$. However this implies that the vertex collection $\{a,b,c_i,c_j\}$ forming a square in $\Gamma_Y$ satisfies condition \hyperref[B1]{(B1)}. While this set might not be maximal with respect to this condition, it is for sure contained in a maximal collection $B$ satisfying conditions \hyperref[B1]{(B1)}, \hyperref[B2]{(B2)} and \hyperref[B3]{(B3)}, corresponding to a rigid vertex group $W_B$. Thus, $G_Y$ is adjacent to the rigid vertex group $W_B$ which it intersects in a subgroup containing $W_{\{a,b,c_i,c_j\}} = D_\infty \times D_\infty$. Hence the intersection is not two-ended.

Assume conversely that $G_Y$ is adjacent to the rigid vertex group $W_B$ and that no pair of vertices in $\mathcal{C}$ is connected by a path in $\Gamma$ that is not passing through $a$ nor $b$. Then each such pair is separated when $a$ and $b$ are removed. Thus at most one of the $c \in \mathcal{C}$ can be contained in $B$ as otherwise condition \hyperref[B1]{(B1)} would be violated. Since the intersection of $G_Y$ and $W_B$ must be infinite, we conclude that $\{a,b\} \subseteq B$. Thus the intersection is either $W_{\{a,b\}}$ or $W_{\{a,b,c\}}$ and therefore two-ended.
\end{proof}

\begin{rem}
Combining Lemma \ref{HangCylEdges2Ended} and Theorem \ref{ThmK4Assump} implies Theorem \ref{ThmEdgeStabTwoEnded}, stating that the JSJ tree of cylinders has two-ended edge stabilizers if and only if there is no uncrossed cut collection containing the corners of a square in the defining graph where the other two corners are connected by a subdivided diagonal. Note that this can be interpreted as a condition about a subdivided $K_4$.
\end{rem}

\begin{rem} \label{NoOverlappingCutTriples}
If $\Gamma$ contains two overlapping cut triples $\{a,b,c\}$ and $\{a,b,c^\prime\}$, then $c$ and $c^\prime$ are connected by a path not passing through $a$ or $b$. Otherwise $a$ and $b$ would be a cut pair, in contradiction to the definition of a cut triple. Thus, if we only consider graphs $\Gamma$ where the JSJ graph of cylinders has two-ended edge groups, overlapping cut triples do not occur in $\Gamma$.

This has an impact on Proposition \ref{PropUncrossedCutPairsCylVert}: Recall that the set $\mathcal{S}$ contains as elements the sets $\{a,b\} \cup \mathcal{C}$, where $\{a-b\}$ is an uncrossed cut collection and $\mathcal{C}$ is the set of common adjacent vertices of $a$ and $b$. Excluding overlapping cut triples implies that every uncrossed cut collection contributes a new element to $\mathcal{S}$. Hence every uncrossed cut collection is in one-to-one correspondence with a cylinder vertex.
\end{rem}

To obtain a JSJ graph of cylinders with two-ended edge groups, we need to refine the \hyperref[StandingAssumption2]{Standing Assumption 2} to:
\begin{Assump} \label{StandingAssumption3}
The defining graph $\Gamma$
\begin{enumerate}[label=(\arabic*)]
    \item has no triangles.
    \item is connected and has neither a separating vertex nor a separating edge.
    \item has a cut collection. 
    \item is not a cycle.
    \item has only uncrossed cut collections $\{a-b\}$ for which for any pair $\{c_1,c_2\} \in \mathcal{C}$ of common adjacent vertices of $a$ and $b$, every path in $\Gamma$ connecting $c_1$ and $c_2$ passes through $a$ or $b$.
\end{enumerate}
\end{Assump}

\begin{rem} \label{IntersectionAlwaysUncrossedCutPair}
Under \hyperref[StandingAssumption3]{Standing Assumption 3}, the proofs of Lemma \ref{HangCylEdges2Ended} and Theorem \ref{ThmK4Assump} imply that the edge stabilizers are either of the form $W_{\{a,b\}}$ or $W_{\{a,b\}} \times W_{\{c\}}$, where $\{a-b\}$ is an uncrossed cut collection and $c$ is a common adjacent vertex of $a$ and $b$. In particular, the latter case can only happen when $a,b$ and $c$ are the corners of a subdivided $K_4$.
\end{rem}

To conclude, we summarize the construction of the JSJ graph of cylinders $\Lambda_c$ in the following theorem:

\begin{thm} \label{SummaryMainThm}
Let $W_\Gamma$ be a RACG with $\Gamma$ satisfying the \hyperref[StandingAssumption3]{Standing Assumption 3}. Then its JSJ graph of cylinders $\Lambda_c$ consists of the following vertices:
\begin{itemize}
    \item For any uncrossed cut collection $\{a-b\} \subseteq EV(\Gamma)$ there is a cylinder vertex with vertex group $W_{\{a,b\} \cup \mathcal{C}}$, where $\mathcal{C}$ is the collection of common adjacent vertices of $a$ and $b$ in $\Gamma$. All the cylinder vertices are either two-ended, VA or VFD.
    \item For any set $A \subseteq V(\Gamma)$ of vertices such that $W_A$ is infinite, $A$ satisfies either conditions \hyperref[A1n]{(A1)}, \hyperref[A2n]{(A2)} and \hyperref[A3n]{(A3)} or condition \hyperref[A*]{(A*)} and $A$ is not contained in a vertex set corresponding to a cylinder vertex group, there is a hanging vertex with vertex group $W_A$. The vertex group is virtually free.
    \item For any set $B \subseteq EV(\Gamma)$ of essential vertices in $\Gamma$ satisfying the conditions \hyperref[B1n]{(B1)}, \hyperref[B2n]{(B2)} and \hyperref[B3n]{(B3)}, there is a rigid vertex with vertex group $W_B$.
\end{itemize}
Furthermore a pair of vertices is connected by an edge if and only if the pair consists of one cylinder vertex corresponding to the cut collection $\{a-b\}$ and one non-cylinder vertex and their vertex groups intersect in a special subgroup containing $W_{\{a-b\}}$. The edge group is the special subgroup generated by the intersection of the corresponding vertex sets. It is two-ended.
\end{thm}

\begin{exam} For the graph $\Gamma$ shown in Figure \ref{fig:ExGoC}, which satisfies \hyperref[StandingAssumption3]{Standing Assumption 3}, we can construct the corresponding JSJ graph of cylinders by reading off the following collections of vertices:
\begin{center}
\begin{tabular}[t]{c|c||c||c}
   uncrossed cut collection & commensurator & hanging & rigid \\\hline
   $a,b$ & $a,m_1,m_2,m_3,b$ & &$a,b,c,d$ \\
   $a,c$ & $a,c$ & $a,k_1,k_2,c$ \\
   $a,d$ & $a,d$ & $a,l_1,l_2,d$ & \\
   $b,c$ & $b,o,c$ & & \\
   $b,d$ & $b,p,d$ & & \\
   $c,d$ & $c,n_1,n_5,d$ & $c,n_2,n_3,n_4,d$ &
\end{tabular}
\end{center}
\begin{figure}[ht] 
    \centering
    \scriptsize
    \def\svgwidth{420pt}
    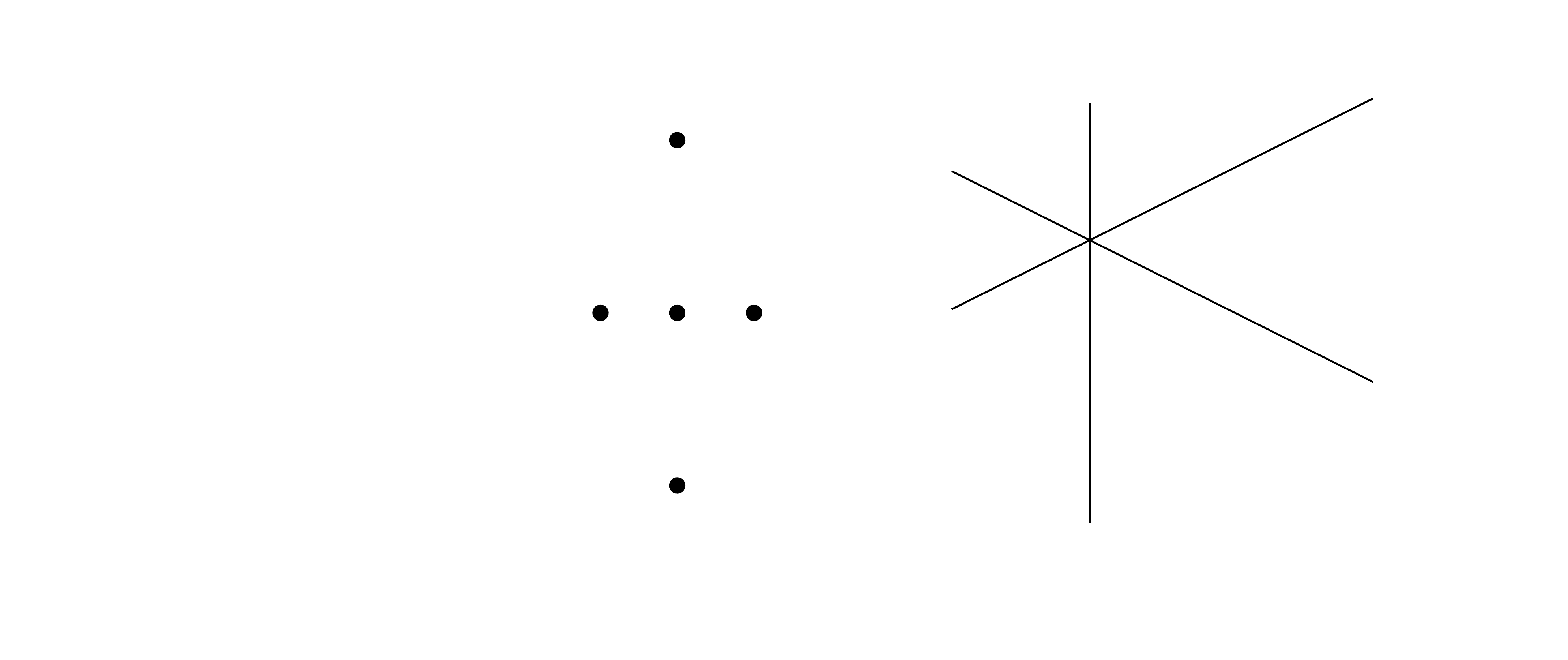
    \caption{$\Lambda_c$ is the JSJ graph of cylinders of the RACG $W_\Gamma$.}
   \label{fig:ExGoC}
\end{figure}
\end{exam}

\section{QI invariance} \label{SectionQIinvariant}
As discussed in Section \ref{SectionPreliminaries}, the feature of interest of the JSJ graph of cylinders is that it can give insight on whether two groups can be QI or not. In the case of certain hyperbolic RACGs, we know by Theorem \ref{ToChypRACG} that all the two-ended cylinder vertices have finite valence in the JSJ tree of cylinders. Thus, if two groups exhibit different valencies at their cylinder vertices, the JSJ trees of cylinders are not isomorphic and thus by Theorem \ref{QIinvar} the groups are not QI.

However, this argument is not applicable in general, since cylinder vertices with one-ended vertex groups do not have finite valence. Nonetheless, we still might be able to distinguish trees of cylinders with infinite valence cylinder vertices, and thus produce an obstruction for a QI, by taking the additional structure coming from the vertex groups and their interplay via edge groups into account. This idea was formalized by Cashen and Martin in \cite{CashenMartinStructureInvar} by the introduction of the so-called \textit{structure invariant}. We first recall its definition and illustrate when it can distinguish two RACGs up to QI and when it cannot. In a second step, we aim to produce a QI between certain groups from identical structure invariants, making the structure invariant a \textit{complete} QI-invariant. 

\subsection{The structure invariant} \label{SectionStrucInvar}
We fix $T$ to be a simplicial tree of countable valence and $G$ to be a group acting on $T$ cocompactly and without edge inversions. We introduce some terminology following \cite[Section 3]{CashenMartinStructureInvar}.

\begin{defin}
For some arbitrary set $\mathcal{O}$ of \textit{ornaments}, a $G$-invariant map $\delta\colon V(T) \rightarrow \mathcal{O}$ is called a \textit{decoration}. The tree $T$ is said to be \textit{decorated}.
\end{defin}

\begin{exam} \label{RelQIType} A classical set of ornaments for a JSJ tree of cylinders $T_c$ is the \textit{vertex type}, that is $\mathcal{O} = \{\textup{‘cylinder'}, \textup{‘hanging'}, \textup{‘rigid'}\}$. A possibly finer decoration is obtained by equipping each vertex $v$ with the ornament consisting of the vertex type and the so-called \textit{relative QI-type} of the corresponding vertex group $G_v$. This relative QI-type is determined as follows: Given the vertex group $G_v$, we consider the set $\mathcal{P}_v$ of distinct Hausdorff equivalence classes in $G_v$ of $G_v$-conjugates of images of the edge injections $\alpha_e\colon G_e \hookrightarrow G_v$ where $e \in E(T_c)$ is an edge incident to $v$. $\mathcal{P}_v$ is often referred to as the \textit{peripheral structure of $G_v$ coming from incident edge groups} or just as the \textit{peripheral structure of $G_v$}. Then the relative QI-type $\llbracket(G_v,\mathcal{P}_v)\rrbracket$ of $G_v$ is the set of all pairs $(Y,P)$, where $Y$ is a geodesic metric space and $P$ is a collection of Hausdorff equivalence classes of subsets of $Y$ such that there is a QI from $G_v$ to $Y$ inducing a bijection from $\mathcal{P}_v$ to $P$. Thus, the relative QI-type captures the structure of the vertex group with respect to its incident edge groups up to QI.
\end{exam} 

\begin{defin}
A decoration $\delta^\prime\colon V(T) \rightarrow \mathcal{O}^\prime$ is called a \textit{(strict) refinement} of the decoration $\delta\colon V(T) \rightarrow \mathcal{O}$ if the $\delta^\prime$-partition $\bigsqcup_{o^\prime \in \mathcal{O}^\prime} {(\delta^\prime)}^{-1}(o^\prime)$ of $V(T)$ (strictly) refines the $\delta$-partition  $\bigsqcup_{o \in \mathcal{O}} \delta^{-1}(o)$. A non-strict refinement is called \textit{trivial}.
\end{defin}

The refinement process used to obtain the structure invariant is the \textit{neighbor refinement}, which is an idea generalizing the degree refinement algorithm known from graph theory. It works as follows:

\begin{constr} \label{neighborRefinement}
Let $\bar{\mathbb{N}} = \mathbb{N} \cup \{\infty\}$ and call $\mathcal{O}_0 = \mathcal{O}$ the \textit{initial set of ornaments} and $\delta_0 = \delta$ the \textit{initial decoration}. Starting from $i=0$, we define for each $i \in \mathbb{N}$ and each $v \in V(T)$ the map
\begin{center}
\begin{tabular}[t]{rccl} 
    $f_{v,i}\colon$ & $\mathcal{O}_i$ & $\rightarrow$ & $\bar{\mathbb{N}}$ \\
   $ $  & $o$ & $\mapsto$ & $|\{w \in \delta_i^{-1}(o) \mid (w,v) \in E(T)\}| \, .$
\end{tabular}
\end{center}
Define $\mathcal{O}_{i+1}$ as $\mathcal{O}_0 \times \bar{\mathbb{N}}^{\mathcal{O}_i}$ and $\delta_{i+1}$ as the pair $(\delta_0(v),f_{v,i}) \in \mathcal{O}_0 \times \bar{\mathbb{N}}^{\mathcal{O}_i}$.
\end{constr}

Cashen and Martin prove the following facts about the maps defined in Construction \ref{neighborRefinement}:
\begin{lem}\textup{\cite[Lemma 3.2, Proposition 3.3]{CashenMartinStructureInvar}}
The map $\delta_{i+1}\colon V(T) \rightarrow \mathcal{O}_{i+1}$ is a decoration refining $\delta_i\colon V(T) \rightarrow \mathcal{O}_i$ for all $i \in \mathbb{N}$. Furthermore, this refinement process stabilizes. That is, there is an $s \in \mathbb{N}$ such that for any $i + 1 \leq s$, the decoration $\delta_{i+1}$ is a strict refinement of $\delta_i$, but for any $i \geq s$, the refinement $\delta_{i+1}$ is trivial.
\end{lem}

\begin{defin}
The decoration $\delta_s\colon V(T) \rightarrow \mathcal{O}_s$ at which the neighbor refinement process stabilizes, is called the \textit{neighbor refinement} of $\delta$.
\end{defin}

To capture the information contained in the neighbor refinement, we define $\pi_0\colon \mathcal{O}_s \rightarrow \mathcal{O}$ to be the projection to the first coordinate. After choosing an ordering on the image $\delta(V(T))$, we denote the $j$-th element as $\mathcal{O}[j]$. Then we can choose an ordering of $\pi_0^{-1}(\mathcal{O}[j]) \cap \delta_s(V(T))$. We order $\delta_s(V(T))$ lexicographically and denote the $i$-th element as $\mathcal{O}_s[i]$.

\begin{defin}
A \textit{structure invariant} $S = S(T,\delta, \mathcal{O})$ is the $|\delta_s(V(T))|^2$-matrix, where
\begin{gather*}
    S_{j,k} = (n_{j,k}, \pi_0(\mathcal{O}_s[j]), \pi_0(\mathcal{O}_s[k])) \, ,
\end{gather*}
with $n_{j,k}$ the number of vertices in $\delta_s^{-1}(\mathcal{O}_s[j])$ adjacent to $\delta_s^{-1}(\mathcal{O}_s[k])$. The second entry of the tuple $S_{j,k}$ is the \textit{row} and the third entry the \textit{column ornament}.
\end{defin}

We can view $S(T,\delta,\mathcal{O})$ as a block matrix, which is well defined up to block permutations and the choice of ordering on $\delta(V(T))$ and $\pi_0^{-1}(\mathcal{O}[j])$.  For economy of notation, we will denote a structure invariant in a table with entries $n_{j,k}$, whose rows and columns are labelled by the initial decoration $\delta(V(T))$, as illustrated in Example \ref{ExampleStructureInvar} or labelled by the vertex orbit representatives carrying the same ornaments, as illustrated in Example \ref{ExRefinement}.

As indicated in the definition, a structure invariant depends on the initial choice of ornaments and decoration.
When we refer to \textit{the} structure invariant, the initial decoration is the one introduced in Example \ref{RelQIType}: the ornaments consist of vertex type and relative QI-type. We call two vertices in the JSJ graph of cylinders \textit{indistinguishable} if they have the same image under $\delta_s$.

By construction, the structure invariant relates to the existence of a tree isomorphim between the JSJ trees of cylinders:
\begin{prop} \textup{\cite[cf.][Proposition 3.7]{CashenMartinStructureInvar}} \label{decopresTreeIsom} Given two groups $G $ and $G^\prime$ with JSJ trees of cylinders $T_c$ and $T_c^\prime$, and $G$- and $G^\prime$-invariant decorations $\delta\colon V(T_c) \rightarrow \mathcal{O}$ and $\delta^\prime\colon V(T_c^\prime) \rightarrow \mathcal{O}$ respectively, there is a decoration-preserving isomorphism $\phi\colon T_c \rightarrow T_c^\prime$ if and only if up to permuting rows and columns within $\mathcal{O}$-blocks, $S(T_c,\delta,\mathcal{O}) = S(T_c^\prime,\delta^\prime,\mathcal{O})$.
\end{prop}

Since the ornaments on a JSJ tree of cylinders consisting of vertex type and relative QI-type determine the structure of the group, we can refine our search to decoration-preserving tree isomorphisms. Hence, by Proposition \ref{decopresTreeIsom}, the structure invariant is indeed a QI-invariant for RACGs with defining graph satisfying \hyperref[StandingAssumption3]{Standing Assumption 3} \cite[cf.][Theorem 3.8]{CashenMartinStructureInvar}. 

\begin{exam}
The two groups with defining graphs illustrated in Figure \ref{fig:EasyEx} serve as an introductory example as they are easily distinguished as non-QI by use of the structure invariant. While the commensurator of the cut pairs $\{a,b\}$ and $\{c,d\}$ in $\Gamma_1$ both give a VFD vertex group, in $\Gamma_2$ the commensurator of $\{c,d\}$ corresponds to a VA vertex group. Since both graphs only have those two uncrossed cut collections, the initial decoration consisting of vertex and relative QI-type already shows that the groups cannot be QI.

\begin{figure}[ht]
    \centering
    \scriptsize
    \def\svgwidth{320pt}
    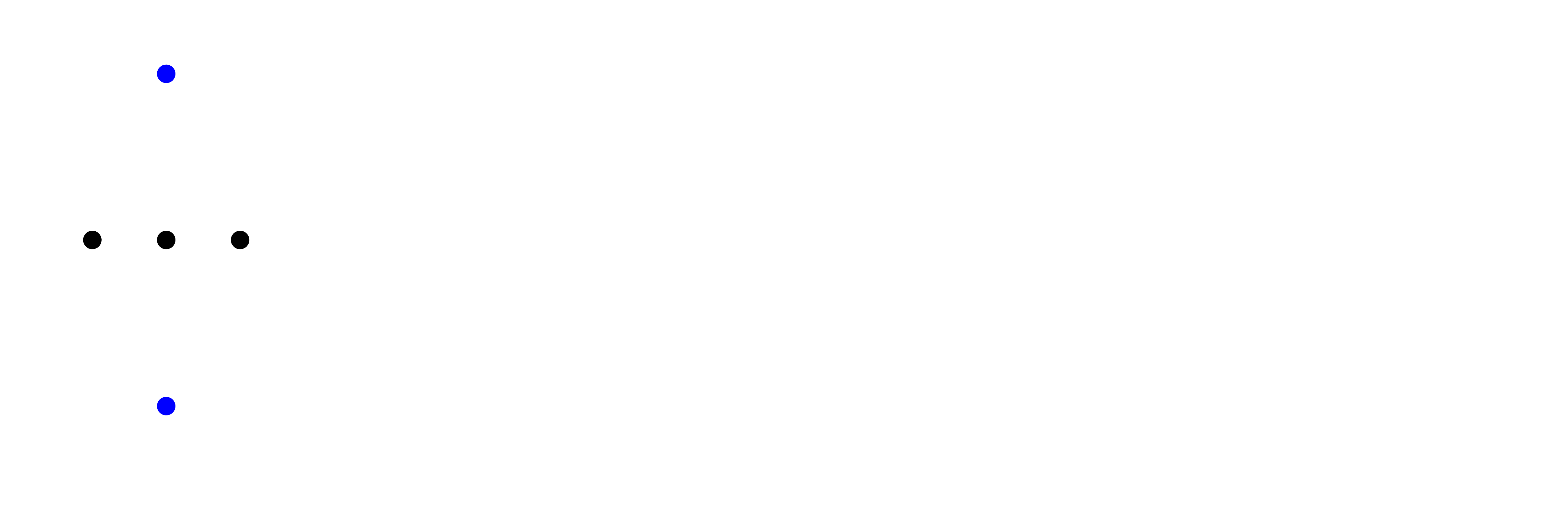
    \caption{The RACGs on the graphs $\Gamma_1$ and $\Gamma_2$ are not QI to each other.}
   \label{fig:EasyEx}
\end{figure}
\end{exam}

\begin{exam} \label{ExampleStructureInvar} To obtain the structure invariants of the JSJ graphs of cylinders $\Lambda_{c,1}$ and $ \Lambda_{c,2}$ for the two RACGs $W_1$ and $W_2$ on the defining graphs $\Gamma_1$ and $\Gamma_2$ respectively, illustrated in Figure \ref{fig:ExampleStructureInvar}, we start with the following initial decorations:
\begin{table}[H]
\small
\begin{tabular}[H]{rclrrcl} 
$\delta\colon V(\Lambda_1)$ & $\rightarrow$ & $\mathcal{O}$ & & $\delta^\prime\colon V(\Lambda_2)$ & $\rightarrow$ & $\mathcal{O}$ \\
$c$ & $\mapsto$ & $(\textup{‘cyl'},\llbracket(\textup{‘VA'},\mathcal{P}_c)\rrbracket)$ & & $c^\prime$ & $\mapsto$ & $(\textup{‘cyl'},\llbracket(\textup{‘VA'},\mathcal{P}_c)\rrbracket)$ \\
$h$ & $\mapsto$ & $(\textup{‘hang'},\llbracket(\textup{‘VF'},\mathcal{P}_h)\rrbracket)$ & & $h_1$, $h_2$ & $\mapsto$ & $(\textup{‘hang'},\llbracket(\textup{‘VF'},\mathcal{P}_h)\rrbracket)$ \\
$r$ & $\mapsto$ & $(\textup{‘rig'},\llbracket (W_{\{a,b,c,d,e,f,g,h\}}, \mathcal{P}_r) \rrbracket)$ & & $r^\prime$ & $\mapsto$ & $(\textup{‘rig'},\llbracket (W_{\{a,b,c,d,e,f,g,h\}}, \mathcal{P}_r) \rrbracket)$ \\
\end{tabular}
\end{table}

\begin{figure}[ht]
    \centering
    \scriptsize
    \def\svgwidth{340pt}
    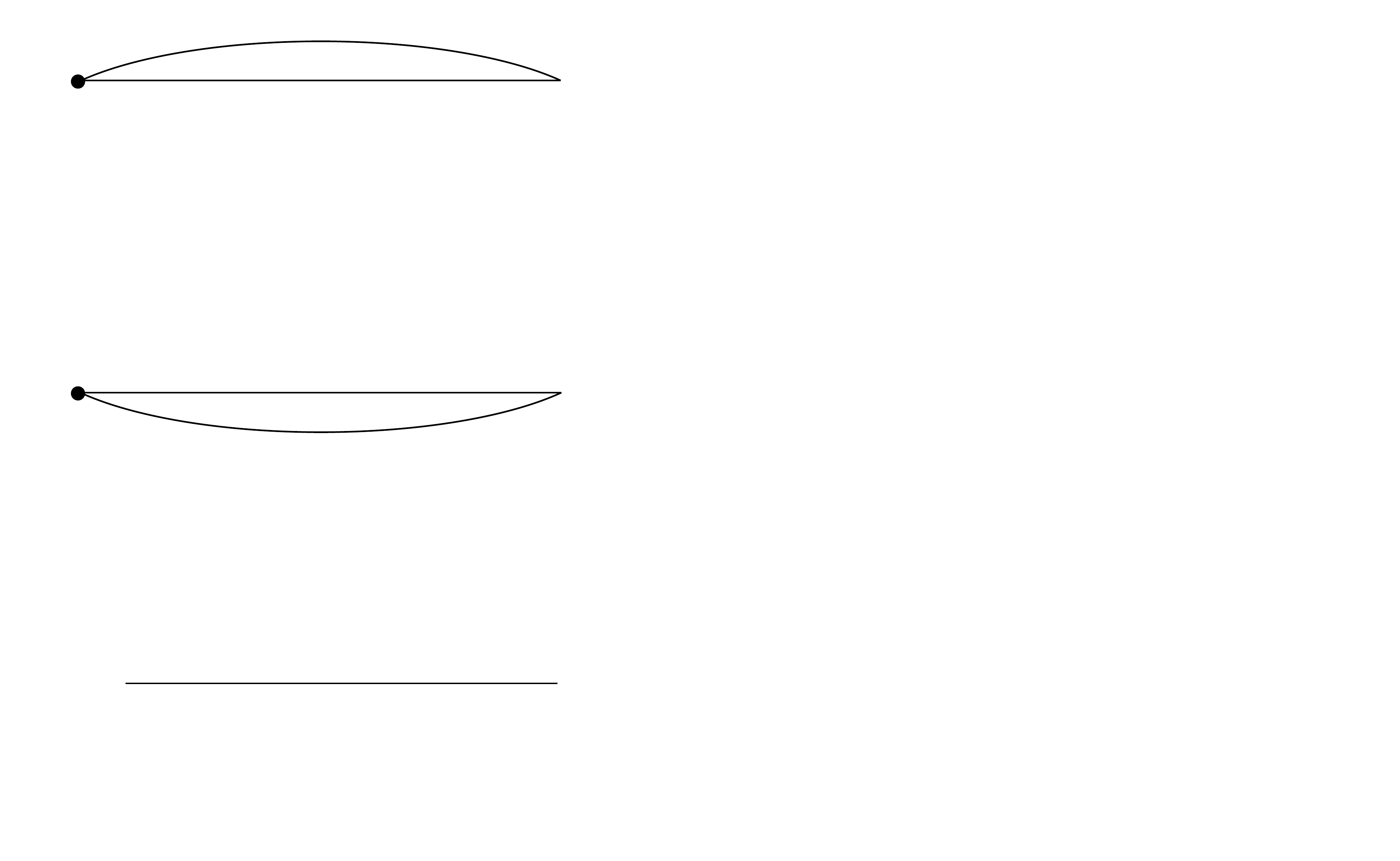
    \caption{We compare the JSJ graphs of cylinders $\Lambda_{c,1}$ and $\Lambda_{c,2}$ of the RACGs $W_{\Gamma_1}$ and $W_{\Gamma_2}$, respectively.}
   \label{fig:ExampleStructureInvar}
\end{figure}

We immediately see that no refinement is possible, the vertices $h_1$ and $h_2$ are indistinguishable and thus the following structure invariant is the same for both JSJ trees of cylinders:

\begin{table}[H]
\footnotesize
\centering
{\renewcommand{\arraystretch}{1.4}
\begin{tabular}[H]{r|c|c|c} 
 & $(\textup{‘cyl'},\llbracket(\textup{‘VA'},\mathcal{P}_c)\rrbracket)$ & $(\textup{‘hang'},\llbracket(\textup{‘VF'},\mathcal{P}_h)\rrbracket)$ & $(\textup{‘rig'},\llbracket (W_{\{a,b,c,d,e,f,g\}}, \mathcal{P}_r) \rrbracket)$ \\\hline
$(\textup{‘cyl'},\llbracket(\textup{‘VA'},\mathcal{P}_c)\rrbracket)$ & 0 & $\infty$ & $\infty$ \\\hline
$(\textup{‘hang'},\llbracket(\textup{‘VF'},\mathcal{P}_h)\rrbracket)$ & $\infty$ & 0 & 0 \\\hline
$(\textup{‘rig'},\llbracket (W_{\{a,b,c,d,e,f,g\}}, \mathcal{P}_r) \rrbracket)$ & $\infty$ & 0 & 0
\end{tabular}}
\end{table}
\end{exam}

Thus by Proposition \ref{decopresTreeIsom}, there is a decoration preserving tree isomorphism between the respective trees of cylinders $T_1$ and $T_2$. This leads to the question whether $W_1$ and $W_2$ are QI, which we will answer in the negative in Example \ref{ExamplVAnotQI}. 

\subsection{Promoting to a QI} \label{SubsectionLocalQIs}
Given two groups $G$ and $G^\prime$ with identical structure invariants and thus with a decoration-preserving isomorphism between their respective JSJ trees of cylinders, we want to determine when we can promote this isomorphism to a QI of the groups. Since any QI between $G$ and $G^\prime$ needs to restrict to a QI locally at each vertex group by Theorem \ref{QIinvar}, the general idea is the following: Start with any local QI between two cylinder vertex groups with the same entry in the structure invariant, which is bijective on the peripheral structures coming from the incident edge groups, and extend it piece by piece from there. By Lemma \ref{LemTypesofCylinderGroups}, we know that in our setting we can encounter either two-ended, VA or VFD cylinder vertex groups. Hence we first determine the possible local QIs for these different cases separately and combine the respective results to find a global QI in a next step.
\begin{rem} \label{technicalities}
All arguments will work along the lines of the ones used in \cite{CashenMartinStructureInvar}, where the case of two-ended cylinder vertices is dealt with. However, at this point we need to clarify three technicalities:
\begin{itemize}
    \item \textit{Rigid vertices need to be handled with special care:}
    \begin{itemize}
        \item While the relative QI-type of rigid vertices might be hard to determine, it can be the crucial ingredient to distinguish groups. In the case of RACGs for instance, this is illustrated by Cashen, Dani and Thomas in \cite{DaniThomasBowditch}. Their Theorem B.1 states that, while all RACGs on $3$-convex subdivided complete graphs with at least 4 essential vertices have isomorphic JSJ trees of cylinders, they are pairwise non-QI. The reason for that lies in the relative QI-type of the rigid vertices. To have more control over rigid vertices, Cashen and Martin restrict to those that have the property of being \textit{quasi-isometrically rigid relative to the peripheral structure} \cite[Definition 4.1]{CashenMartinStructureInvar}. For example, free rigid vertex groups will have this property by \cite{cashenmacura}. Under this additional assumption, another ornament, the \textit{relative stretch factor}, can be introduced to decorate edges and help distinguish rigid vertices \cite[Section 4]{CashenMartinStructureInvar}. However, whether rigid vertices in JSJ trees of cylinders of RACGs have this or a similar sufficient property (for example the related \textit{right-angled Artin groups} splitting over cyclic groups do \cite[cf.][Section 6]{margolis2018quasi}) is yet to be determined.
        \item Another issue caused by rigid vertices is that they might have adjacent edges whose edge groups are not two-ended as shown in Theorem \ref{ThmK4Assump}. Nguyen and Tran give in \cite{NguyenTranPlanar} a complete QI-classification of a class of RACGs with such edge groups: The defining graphs are connected, trianlge-free and planar and have more than 4 vertices, no seperating edge or vertex and a property called $\mathcal{CFS}$ \textit{(constructed from squares}, \cite[Definition 1.3]{BehrstockCFS}). In the proof they use the \textit{maximal suspension decomposition} and the properties of its corresponding vertex groups. However, for the groups they consider, this turns out the be in correspondence with the JSJ graph of cylinders and the decoration consisting of the relative QI-type. 
        \item Moreover, in \cite{bounds2020quasi}, Bounds and Xie show that RACGs whose defining graphs are \textit{generalized thick $m$-gons} exhibit a strong form of QI-rigidity: They are QI if and only if their defining graphs are isomorphic.
    \end{itemize}
    For simplicity, we focus on groups without any rigid vertices or on pairs of groups which have isomorphic rigid vertex groups as in Examples \ref{ExampleStructureInvar} and \ref{ExQI2}.
    \item \textit{Work on the geometric tree of spaces:} To make technical details more economic, instead of working on graphs of groups, Cashen and Martin state their results for a slightly modified space, the \textit{geometric tree of spaces} $X$ of $G$ over $T_c$. The construction of $X$ is standard and useful as $X$ is QI to $G$. Essentially, $X$ is produced from the JSJ graph of cylinders $\Lambda_c$ by substituting all groups of the same relative QI-type by a uniform \textit{model space} representing the equivalence class. Thus, instead of a subgroup $G_t$ we have a subspace $X_t$ for every $t \in T_c$. Most importantly, if two groups $G$ and $G^\prime$ exhibit subgroups $G_t$ and $G^\prime_{t^\prime}$ with equivalent relative QI-types in their JSJ graphs of cylinders, we choose the same model space $X_t$ for both $G_t$ and $G^\prime_{t^\prime}$. If convenient, we will state results in terms of the geometric tree of spaces $X$, but spare the bookkeeping, which is done thoroughly in Sections 7.2 and 2.5 of \cite{CashenMartinStructureInvar}.
    \item \textit{Partial orientations can be omitted:} For the sake of completeness it should be mentioned that, apart from the neighbor refinement, Cashen and Martin introduce the \textit{cylinder} and the \textit{vertex refinement}, depending on a \textit{partial orientation} chosen essentially on all two-ended spaces. However, since all infinite RACGs and thus all edge groups in $\Lambda_c$ contain an infinite dihedral group $D_\infty$, the orientation can always be reversed. Thus, the refinement processes become trivial and shall therefore be left out of our considerations.
\end{itemize}
\end{rem}

\subsubsection{Two-ended cylinder vertices}
In case all cylinder vertex groups are two-ended, like for instance for hyperbolic groups, Cashen and Martin give a structure invariant, which is a complete QI-invariant. Their result, stated for RACGs splitting over two-ended subgroups and thus refining Theorem \ref{QIinvar}, is reiterated in the following proposition:
\begin{prop} \textup{\cite[Theorem 7.5]{CashenMartinStructureInvar}} \label{CompleteInvarVZcylinders}
Let $W$ and $W^\prime$ be two finitely presented, one-ended RACGs with non-trivial JSJ decomposition over two-ended subgroups such that cylinder stabilizers are two-ended and all non-cylinder vertex groups are either hanging or quasi-isometrically rigid relative to the peripheral structure. Define $T$ to be the JSJ tree of cylinders of $W$ and $X$ to be the geometric tree of spaces of $W$ over $T$. The initial decoration $\delta_0$ on $T$ takes vertex type, relative QI-type and the relative stretch factor into account. Let $\delta$ be the neighbor refinement of $\delta_0$. Analogously, we define $T^\prime$, $X^\prime$, $\delta_0^\prime$ and $\delta^\prime$ for $W^\prime$. Then $W$ and $W^\prime$ are QI if and only if there is a bijection $\beta\colon \delta(T) \rightarrow \delta^\prime(T^\prime)$ such that
\begin{enumerate}
    \item $\delta_0 \circ \delta^{-1} = \delta_0^\prime \circ (\delta^\prime)^{-1} \circ \beta$;
    \item $S(T,\delta,\mathcal{O}) = S(T^\prime, \delta^\prime, \mathcal{O}^\prime)$ in the $\beta$-induced ordering;
    \item for every ornament $o \in \mathcal{O}$ with $\delta^{-1}(o)$ containing non-cylinder vertices, there is a vertex $v \in \delta^{-1}(o)$ and a vertex $v^\prime \in (\delta^\prime)^{-1}(\beta(o))$  such that there is a QI between the vertex spaces $X_v$ and $X^\prime_{v^\prime}$, which is bijective on the peripheral structures $\mathcal{P}_v$ and $\mathcal{P}^\prime_{v^\prime}$ and respecting the decorations $\delta$ and $\delta^\prime$ respectively.
\end{enumerate}
\end{prop}

The inductive construction of the QI in their proof will serve as a blueprint for the proof of the general Theorem \ref{finalresult}.

\subsubsection{VFD cylinder vertices}

It turns out that VFD cylinder vertex groups have enough flexibility to always find a QI between cylinder vertices with this same entry in the structure invariant. We construct this local QI in the following simplest setting:

\begin{prop} \label{VFD}
Let $W_1$ and $W_2$ be two RACGs on defining graphs satisfying \hyperref[StandingAssumption3]{Standing Assumption 3} with identical structure invariants and one single cylinder vertex $v_1$ and $v_2$ in the JSJ graph of cylinders $\Lambda_1$ and $\Lambda_2$ respectively. Let the vertex groups $V_1$ and $V_2$ of $v_1$ and $v_2$ respectively be VFD. Then there is a QI between $V_1$ and $V_2$ that is bijective on the respective peripheral structures.
\end{prop}

\begin{proof}
The set-up is the following: Both JSJ graphs of cylinders $\Lambda_1$ and $\Lambda_2$ look like stars, with the cylinder vertex in the middle and their neighbors grouped into $j < \infty$ classes of indistinguishable vertices. Suppose at first that $j=1$.

Thus for $i \in \{1,2\}$, each $\Lambda_i$ consists of one cylinder vertex $v_i$, which has a vertex group of the form $V_i = W_{\mathcal{C}_i} \times D_\infty$. The copy of $D_\infty$ is generated by non-adjacent vertices of a cut collection and $W_{\mathcal{C}_i}$ is generated by the set $\mathcal{C}_i$ of their common adjacent vertices. By assumption $W_{\mathcal{C}_i}$ is virtually free, thus $|\mathcal{C}_i|>2$. At the cylinder vertex $v_i$ in the middle, there is a set $\mathcal{N}_i$ of $e_i$ indistinguishable non-cylinder vertex groups of the same relative QI-type attached along a two-ended edge group. These edge groups are either a copy of $D_\infty$ or of $D_\infty \times \mathbb{Z}_2$ with $\mathbb{Z}_2 = W_{\{c\}}$ for some $c \in \mathcal{C}_i$ by Remark \ref{IntersectionAlwaysUncrossedCutPair}. Thus, in the corresponding JSJ tree of cylinders, the vertex $1 \cdot V_i$ has infinitely many adjacent vertex groups corresponding to cosets of the form $gN$: The group $N$ is an element of $\mathcal{N}_i$ and $g \in V_i$ is either any word in $W_{\mathcal{C}_i}$ or a word in $W_{\mathcal{C}_i}$ not ending on $c$, depending on whether the edge group along which $N$ attaches is $D_\infty$ or $D_\infty \times W_{\{c\}}$.

We want to interpret this set-up in terms of Cayley graphs in order to prove Claim \ref{ClaimReductionToTree}. Before, we need to fix some terminology:

As a \textit{graph $\Delta$ with tangling edges $E$} we understand some \textit{base} graph $\Delta$, where at each vertex in $V(\Delta)$ we add some additional neighbors, all of valence 1. Each such additional edge is labelled by an element of $E$ and is called a \textit{tangling edge}. In the new graph, we can think of each tangling edge as the pair $(v,i)$, where $i \in \mathbb{N}$ counts the edges attaching at the base vertex $v \in V(\Delta)$ in $\Delta$. We denote the resulting graph as $\Delta \cup E$, where the union happens via the implicit attaching map. Let $k(v)$ be the number of edges \textit{tangling} at $v \in V(\Delta)$. Then we can interpret the set $E$ as $E = \{t_{v,i} \mid v \in V(\Delta), \, i \in \{0, \dots, k(v)-1\}\} $ where $t_{v,i}$ denotes the $i$-th tangling edge at vertex $v \in V(\Delta)$.

\begin{claim} \label{ClaimReductionToTree}
The problem of finding a QI between $V_1$ and $V_2$ that is bijective on the respective peripheral structures can be reduced to finding a QI between two identical infinite, regular trees $T$ with tangling edge sets $E_1$ and $E_2$ such that the occurring numbers of tangling edges $\{k_i(v) \mid v \in V(T)\}$ in $T \cup E_1$ and $T \cup E_2$ differ. In addition, this QI must be bijective on the tangling edges.
\end{claim}

\begin{proof} (of Claim \ref{ClaimReductionToTree})
The idea of the reduction is the following: The Cayley graph of $W_{\mathcal{C}_i}$ reduces to the base tree and the tangling edges are in correspondence with the different cosets $gN$.

We start the reduction process with the object $X_i$, illustrated in Figure \ref{fig:TanglingEdgesFaces}, constructed as follows: Note first that the Cayley graph of the cylinder vertex group $V_i = W_{\mathcal{C}_i} \times D_\infty$ is the direct product of the Cayley graph $T_i$ of $W_{\mathcal{C}_i}$ and the line $D$ that is the Cayley graph of $D_\infty$. This is true because with the correct choice of generating sets, the Cayley graph of a direct product is the direct product of the Cayley graphs. Since all the vertices in $\mathcal{C}_i$ are pairwise non-connected in the defining graph, the Cayley graph $T_i$ of $W_{\mathcal{C}_i}$ is a $|\mathcal{C}_i|$-regular tree. We can think of each coset $gN$ adjacent to the vertex $1 \cdot V_i$ as attaching in this Cayley graph. If $g$ can be any word in $W_{\mathcal{C}_i}$, the coset attaches at the vertex $g$ in $T_i$ and along the line $D$. If $g$ is a word in $W_{\mathcal{C}_i}$ not ending on $c$, the coset attaches along the edge $c$ starting at the vertex $g$ in $T_i$ and along the line $D$. Either way, we can think of the $e_i$ different cosets $gN$ as $e_i$ possibly thickened half-planes at the vertex $g$ in $T_i$ attached along the line $D$. Note that at one vertex $g$ it can happen that there attach both \textit{thick} half-planes along an edge and \textit{thin} half-planes at the vertex. We call the constructed object $X_i$.

\begin{figure}[ht]
    \centering
    \footnotesize
    \def\svgwidth{430pt}
    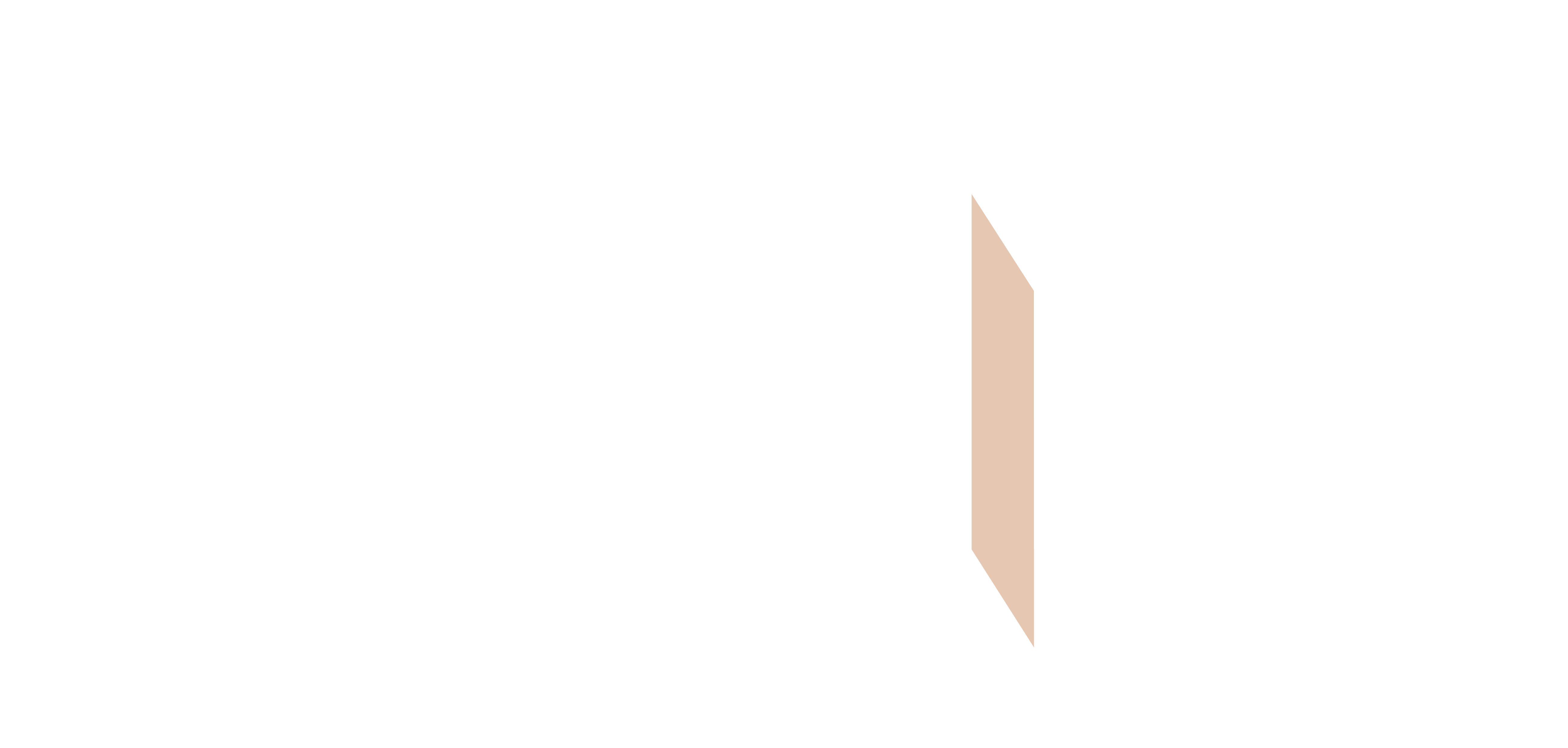
    \caption{We replace thick half-planes by thin half-planes and collapse the $D$-direction.}
   \label{fig:TanglingEdgesFaces}
\end{figure}

Since $X_i$ captures the structure of the group, the task of finding a QI between $V_1$ and $V_2$ that is bijective on the peripheral structure is done if we can show that there is a QI between $X_1$ and $X_2$ that is bijective on the half-planes corresponding to the cosets. For the reduction, we squish for any $N \in \mathcal{N}_i$ attaching along a $D_\infty \times W_{\{c\}}$ the corresponding thick half-plane: Replace each such thick half-plane attaching along an edge by a thin half-plane attached at the terminal vertex of the attaching edge that has less half-planes attached. Then, at each vertex in $T_i$ there attaches some positive number of thin tangling half-planes corresponding to the cosets $gN$. We reinterpret this object as $(T_i \cup E_i) \times D$, where $E_i$ is a set of tangling edges. It suffices to find a QI between $(T_1 \cup E_1) \times D$ and $(T_2 \cup E_2) \times D$ that is bijective on the tangling half-planes, because this immediately implies that we can find a QI between the trees with thick tangling half-planes simply by extending the map along the attaching edges via the identity.

However, now it is enough to find a QI between $T_1 \cup E_1$ and $T_2 \cup E_2$ which is bijective on the tangling edges, because again, this immediately implies that we can find a QI between $(T_1 \cup E_1) \times D$ and $(T_2 \cup E_2) \times D$ this time by extending the map to $D$ via the identity.

So, to find this QI, recall the well-known fact that two regular trees are QI to each other by contracting or inserting one edge path of a certain finite length at each vertex of the first tree to turn it into the second. Thus, we start with the tree with smaller regularity, say without loss of generality $T_1$ and perform this operation on the edges to obtain a tree $T$ that is QI to $T_1$ and isomorphic to $T_2$. While this contraction and insertion of edges redistributes the tangling edges of $T_1$, since we started with the tree with smaller regularity, the resulting $T$ also has at least one tangling edge at each vertex. Hence there is a QI between $T_1 \cup E_1$ and $T_2 \cup E_2$ that is bijective on the tangling edges if we can find a QI between $T \cup E_1$ (with some adjusted tangling edge set $E_1$) and $T \cup E_2$ that is bijective on the tangling edges.

We could keep track of the exact number $k_i(v)$ of tangling edges at each vertex $v$ in $V(T)$ of $T \cup E_i$. However, since this would require a technical case distinction, we suppress the details. In general, the number $k_i(v)$ of tangling edges at each vertex $v$ varies. However, most importantly, we see from an analysis of the reduction process that all vertices have a bounded number of tangling edges, that is all vertices have at least $y_i > 0$ and at most $x_i < \infty$ tangling edges, i.e. $0 < y_i \leq k_i(v) \leq x_i < \infty $ for all $v \in V(T)$.

With this process, we have reduced the problem of finding a QI between $V_1$ and $V_2$ that is bijective on the respective peripheral structures to finding a QI between two copies of an infinite, $r$-regular tree $T$ with differing occurring numbers $\{k_i(v) \mid v \in V(T)\}$ of tangling edges at its vertices, that is bijective on the tangling edges.
\end{proof}

So, as Claim \ref{ClaimReductionToTree} suggests, we aim to find a QI $q$ from $T \cup E_1$ to $T \cup E_2$, where the base graph $T$ is an infinite $r$-regular tree with distinguished base vertex and $q$ is bijective on the tangling edges. Without loss of generality we set the maximal number $x_1$ of edges attaching at a base vertex in $T \cup E_1$ to be greater than the maximal number $x_2$ of edges attaching at a base vertex in $T \cup E_2$.\medskip

We define the following notion on $T \cup E_i$: Given an edge $e \in E(T)$ with some tangling edge $t \in E_i$ at the vertex $o(e)$, we call it a \textit{slide} along $e$ if we detach $t$ from $o(e)$ and reattach it at $t(e)$.

\begin{claim} \label{UniformBoundonTheSlide}
There is a constant $d \in \mathbb{N}$ such that in $T \cup E_1$ every tangling edge at each vertex of $T$ needs to slide at most along $d$ edges of $T$ away from the distinguished base vertex, such that the resulting graph is isomorphic to $T \cup E_2$.
\end{claim}

Now we get the desired QI $q$ from $T \cup E_1$ to  $T \cup E_2$ which is bijective on the tangling edges: The sliding process in Claim \ref{UniformBoundonTheSlide} defines a bijective map $q^\prime: E_1 \rightarrow E_2$ mapping each edge in $E_1$ to the edge in $E_2$ on whose position it is slid to. We define $q: T \cup E_1 \rightarrow T \cup E_2$ to be the map that is the identity on $T$ and $q^\prime$ on the elements of $E_1$ as "half-open" edges without the endpoint contained in $T$. 

Let $t_{v,j}$ and $t_{v^\prime,j^\prime}$ in $E_1$ be two tangling edges based at $v$ and $v^\prime$ in $T$ respectively. Since tangling edges are always slid away from the distinguished base vertex, their images can get at most $d$ edges closer to each other than $v$ and $v^\prime$ are. Hence
$$d(t_{v,j}, t_{v^\prime,j^\prime}) - d \leq d(q(t_{v,j}),q(t_{v^\prime,j^\prime})) \, ,$$ which gives the lower QI-bound.
For the upper QI-bound note that since both tangling edges are slid at most along $d$ edges, their distance can grow at most by $2d$, that is
$$ d(q(t_{v,j}),q(t_{v^\prime,j^\prime})) \leq  d(t_{v,j}, t_{v^\prime,j^\prime}) + 2d \, .$$
Since a vertex $w \in V(T)$ is not moved by $q$, analogous bounds hold for $d(q(t_{v,j}),q(w))$. This implies that $q$ is a quasi-isometric embedding. The bijectivity of $q^\prime$ ensures the quasi-surjectivity of $q$. Therefore $q$ is a QI and the only thing left to prove is Claim \ref{UniformBoundonTheSlide}:
\begin{proof} (of Claim \ref{UniformBoundonTheSlide}) 
For simplicity we want to define the graphs $T_i^\prime$, which are identical to $T \cup E_i$ with the exception that in $T_i^\prime$, the base vertex of $T$ does not carry any tangling edges. Since the number of tangling edges we remove is bounded by $k < \infty$, the claim remains true if we can prove it on $T_i^\prime$ and the argument works analogously. However, if we work on $T_i^\prime$, we can give $d$ explicitly in terms of the maximal number $x_1$ of tangling edges at a vertex in $T \cup E_1$, the minimal number $y_2$ of tangling edges at a vertex in $T \cup E_2$ and the degree $r$ of the regularity of $T$ as follows: $$d = \left\lceil \tfrac{\log(\frac{x_1}{y_2})}{\log(r-1)}\right\rceil \, .$$ If the base vertex carries at most $k$ tangling edges as well, $d$ is bounded by $\left\lceil \tfrac{\log(\frac{x_1}{y_2})}{\log(r-1)}\right\rceil + k$, a complication we avoid without any loss of generality. Now the key feature of the proof is the following algorithm:
\begin{alg} \label{AlgVFree}
Since we always need to slide away from the base vertex, we can reduce the problem by dividing the tree $T$ into $r$ subtrees by removing the base vertex, which does not have any tangling edges. Then we have $r$ rooted trees, where the root has $r-1$ outgoing edges. We consider one rooted tree $R$, which is oriented away from the root $*$. A vertex is at level $l$ of $R$ if it has distance $l$ to the root $*$. Note that every vertex in $R$, with exception of the root, has one incoming, $(r-1)$ outgoing and at least $y_i$ and at most $x_i$ tangling edges. The root has $r-1$ outgoing and at least $y_i$ and at most $x_i$ tangling edges.

Now the idea is the following: Every vertex receives some tangling edges via a slide along its incoming edge and superfluous tangling edges leave the vertex via a slide along the outgoing edges. The sliding process follows two rules:
\begin{enumerate}
    \item The distribution of the superfluous tangling edges along the $r-1$ outgoing edges is uniform. \label{Rule1}
    \item The edges that are kept at each vertex are always the ones that have been slid the furthest. \label{Rule2}
\end{enumerate}
Without loss of generality, we can assume that all vertices in $T_1^\prime$ have the maximal number of $x_1$ tangling edges and all vertices in $T_2^\prime$ with the minimal number of $y_2$ tangling edges. If the bound $d$ works for these special cases, it works for the numbers of tangling edges in between.

We show by induction on the level $l$ of $R$ that $d$ satisfying $\left\lceil \frac{x_1}{({r-1})^d} \right\rceil \leq y_2$ works as a uniform bound. Consider the root of $R$ at level $0$ as the base case. We need to keep $y_2$ edges at the root and by rule \ref{Rule2}, we keep at most a total of $(r-1) \cdot y_2$ edges coming from the root at level 1. In general we keep at most a total of $(r-1)^{i}\cdot y_2$ edges coming from the root at level $i$. But since
$$\sum\limits_{i=0}^d (r-1)^i \cdot y_2 \geq (r-1)^d \cdot y_2 \geq x_1 \, ,$$
it is immediate that none of the $x_1$ edges coming from the root will be slid more than $d$ steps.

For the inductive step, suppose that each edge up to level $l$ will be slid at most along $d$ edges. We consider a vertex $v$ at level $l+1$. If we slide its tangling edges along $d$ edges in $R$, they are now attached at a vertex at level $l+1+d$. But by hypothesis any edge slid away from a vertex at level $l$ or any level above cannot be attached at level $l+1+d$. Thus the edges from level $l+1$ are the ones that have been slid the furthest, so by rule \ref{Rule2} they are the ones that need to stay. However, by choice of $d$, there are at most $y_2$ edges coming from $v$ per vertex at level $l+1+d$. Therefore, no edges coming from level $k+1$ are slid any further, proving that the chosen $d$ gives a uniform bound.
\end{alg}
\end{proof}

The way to interpret Algorithm \ref{AlgVFree} is that in the JSJ graphs of cylinders, we can duplicate or collapse the neighboring vertices of $v_1$ of the same QI-type to match the neighbors of $v_2$.

In order to produce a QI between $V_1$ and $V_2$ when $\Lambda_1$ and $\Lambda_2$ have $j\geq 2$ classes of indistinguishable vertices attached at the cylinder vertex, we apply the Claims \ref{ClaimReductionToTree} and \ref{UniformBoundonTheSlide} and execute Algorithm \ref{AlgVFree} for each class individually.
\end{proof}

\subsubsection{VA cylinder vertices}
The flexibility of the VA cylinder vertices lies in between the flexibility of the other two types: In the tree of cylinders, they have infinite valence like the VFD cylinder vertices. However, in order to get a QI from one VA cylinder vertex group to another, the different classes of indistinguishable neighboring vertex groups must occur with matching densities in the respective JSJ graphs of cylinders. This behaviour is similar to the two-ended cylinder vertices. The robustness comes from the fact that the QI cannot be of any type, but it must be bounded distance from scaling by precisely the density. Shepherd and Woodhouse also make use of these densities in \cite[Section 5.6]{shepherd2020quasi}.

As for the VFD cylinders, we construct the local QI in the following simplest setting:

\begin{prop} \label{1VACyl}
Let $W_1$ and $W_2$ be two RACGs on defining graphs satisfying \hyperref[StandingAssumption3]{Standing Assumption 3} with identical structure invariants and one single cylinder vertex $v_1$ and $v_2$ in the JSJ graph of cylinders $\Lambda_1$ and $\Lambda_2$ respectively. Let the cylinder vertex group $V_1 \cong D_\infty \times D_\infty$ and $V_2 \cong D_\infty \times D_\infty$ at $v_1$ and $v_2$ respectively be VA. Suppose at $v_1$ and $v_2$ attach $e_1$ and $e_2$ neighbors of the same class of indistinguishable vertices respectively. The number $e_1$ decomposes as the sum of $m_1$ vertices attaching along a $D_\infty$-edge and $n_1$ vertices attaching along a $D_\infty \times \mathbb{Z}_2$-edge. Analogously, $e_2 = m_2 + n_2$.

There is a QI from $V_1$ to $V_2$ that is bijective on the respective peripheral structures if and only if there is a QI that is the identity map on the first $D_\infty$-copy of $V_1$ and $V_2$ and that scales under the natural identification with $\mathbb{Z}$ the second $D_\infty$-copy of $V_1$ to the second $D_\infty$-copy of $V_2$ by $$\frac{2\, m_1 + n_1}{2 \, m_2 + n_2} \, .$$

Furthermore, every QI between $V_1$ and $V_2$ that is bijective on the respective peripheral structures is bounded distance from one of the form
\begin{center}
\begin{tabular}[t]{rrcl}
     $\psi \colon$ & $D \times L$ &$\rightarrow$ & $D \times L$  \\
     & $(x,y)$ & $\mapsto$ & $(\psi^\prime(x,y),\psi^{\prime\prime}(x,y))$ \, , 
\end{tabular}
\end{center}

where $D$ and $L$ are Cayley graphs of $D_\infty$ and $\psi^{\prime\prime}$ is scaling by  $\tfrac{2\, m_1 + n_1}{2 \, m_2 + n_2}$.  
\end{prop}

\begin{proof} The proof resembles the proof of Proposition \ref{VFD}, we have a similar set-up: For $i \in \{1,2\}$, the JSJ graph of cylinders $\Lambda_i$ looks like a star with one VA cylinder vertex $v_i$ in the middle. The VA vertex group $V_i$ at $v_i$ corresponds to the uncrossed cut collection $\{a_i-b_i\}$ with common adjacent vertex set $\{s_i,t_i\}$, thus $V_i = W_{\{a_i,b_i\}} \times W_{\{s_i,t_i\}} = D_\infty \times D_\infty$.

As before, in $\Lambda_i$, at the cylinder vertex $v_i$ there is a set $\mathcal{N}_i$ of $e_i$ indistinguishable non-cylinder vertex groups of the same relative QI-type. Of these, $m_i$ are attached along a $W_{\{a_i,b_i\}} = D_\infty$-edge group and $n_i$ are attached along a $W_{\{a_i,b_i,s_i\}} = D_\infty \times \mathbb{Z}_2$-edge group or a $W_{\{a_i,b_i,t_i\}} = D_\infty \times \mathbb{Z}_2$-edge group (cf. Remark \ref{IntersectionAlwaysUncrossedCutPair}). In the corresponding JSJ tree of cylinders, the vertex $1 \cdot V_i$ has infinitely many adjacent vertex groups corresponding to cosets of the form $gN$. The group $N$ is an element of $\mathcal{N}_i$ and $g \in V_i$ is either any word in $W_{\{s_i,t_i\}}$ or any word in $W_{\{s_i,t_i\}}$ not ending on $s_i$ or on $t_i$, depending on whether $N$ attaches along the edge group $W_{\{a_i,b_i\}}$ or $W_{\{a_i,b_i,t_i\}}$ or $W_{\{a_i,b_i,s_i\}}$ respectively. \\

Again, we want to interpret this set-up in terms of Cayley graphs in order to prove the following claim:

\begin{claim} \label{ReduceZwithTangling}
The problem of finding a QI between $V_1$ and $V_2$ that scales $W_{\{s_1,t_1\}}$ to $W_{\{s_2,t_2\}}$ by $\tfrac{2\, m_1 + n_1}{2 \, m_2 + n_2}$ and that is bijective on the respective peripheral structures can be reduced to finding a QI between two copies of the number line with different occurring numbers of tangling edges that scales the number line by $\tfrac{2\, m_1 + n_1}{2 \, m_2 + n_2}$.
\end{claim}

\begin{proof} (of Claim \ref{ReduceZwithTangling})
Analogous to the procedure for a VFD cylinder vertex in the proof of Proposition \ref{VFD}, we use the Cayley graph of $V_i$. It is given by the $a_i b_i \times s_i t_i$-grid. Note that the bi-labelled $s_it_i$-line $L$ corresponds to $T_i$ in the proof of Proposition \ref{VFD} and the bi-labelled $a_ib_i$-line corresponds to $D$. If $g$ can be any word in $W_{\{s_i,t_i\}}$, the coset $gN$ attaches at vertex $g$ in $L$ along $D$. If $g$ is a word in $W_{\{s_i, t_i\}}$ not ending on $s_i$, the coset $gN$ attaches along the edge $s_i$ starting at the vertex $g$ in $L$ and along $D$ and if $g$ is a word in $W_{\{s_i, t_i\}}$ not ending on $t_i$, the coset $gN$ attaches along the edge $t_i$ starting at the vertex $g$ in $L$ and along $D$. Either way, we can think of the $e_i$ different cosets $gN$ as $e_i$ possibly thickened half-planes at the vertex $g$ in $L$ attached along the line $D$. Using the notation from the proof of Proposition \ref{VFD}, we call this object $X_i$, it is illustrated in Figure \ref{fig:VATanglingFaces}.

\begin{figure}[ht]
    \centering
    \scriptsize
    \def\svgwidth{430pt}
    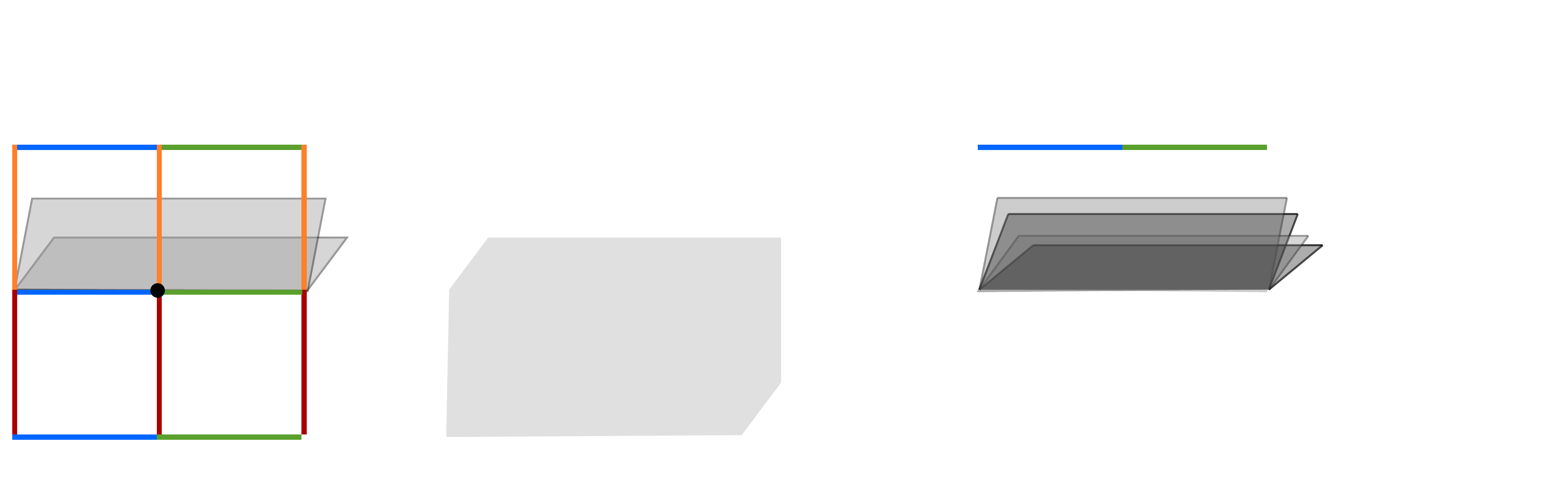
    \caption{We replace thick half-planes by thin half-planes and collapse the $D$-direction.}
   \label{fig:VATanglingFaces}
\end{figure}

Again, we obtain a QI between $V_1$ and $V_2$ that is bijective on the peripheral structure if we can find a QI between $X_1$ and $X_2$ that is bijective on the half-planes corresponding to the cosets. Unlike in the proof of Proposition \ref{VFD}, we need to make sure that the following reductions work both ways in order to prove the fact about the scaling, that is, we show that we find the desired QI between the reduced objects if and only if we find one between the original ones. 

First, we get rid of the thick half-planes, as in the proof of Claim \ref{ClaimReductionToTree}. For any $N \in \mathcal{N}_i$ we find attaching along a $W_{\{a_i,b_i,s_i\}}$- or $W_{\{a_i,b_i,t_i\}}$-edge, we need to squish the corresponding thick half-planes: Replace half of the thick half-planes attaching along an edge $e \in \{s_i,t_i\}$ by a thin half-plane attached at $o(e)$ and the other half at $t(e)$. Then, at each vertex in $L$ there attach $\frac{n_i}{2}$ thin tangling half-planes coming from thick ones and $m_i$ originally thin ones. In total, there attach $m_i + \frac{n_i}{2}$ thin half-planes at each vertex. Of course, it can happen that $n_i$ is odd and we have produced "half a half-plane" with this procedure. However, this will not affect the rest of the argument.

We reinterpret this object as $(L \cup E_i) \times D$, where $E_i$ is the set of tangling edges. It suffices to find a QI between $(L \cup E_1) \times D$ and $(L \cup E_2) \times D$ that scales $L$ by $\tfrac{2\, m_1 + n_1}{2 \, m_2 + n_2}$ that is bijective on the tangling half-planes. The existence of such a QI immediately implies that we can find a QI between the grids with thick tangling half-planes simply by extending the map along the attaching edge via the identity. Conversely, if we find a QI between two grids with thick tangling half-planes which is bijective on all tangling half-planes, this means that the horizontal $D$-lines are preserved. Thus we can restrict this QI to obtain the desired QI between the grids with thin tangling edges.

In the second reduction step we check that it is enough to find a QI between $L \cup E_1$ and $L \cup E_2$ that is bijective on the tangling edges and scaling $L$ by $\tfrac{2\, m_1 + n_1}{2 \, m_2 + n_2}$. Again, given such a QI, we can find the desired QI between $D \times (L \cup E_1)$ and  $D \times (L \cup E_2)$ by extending the map via the identity to $D$.

For the converse, suppose that $\psi\colon D \times (L \cup E_1) \rightarrow D \times (L \cup E_2)$ is a QI which is bijective on the tangling half-planes and scaling $L$ by $\tfrac{2\, m_1 + n_1}{2 \, m_2 + n_2}$. Restricted to the grid, $\psi$ is of the form:
\begin{center}
\begin{tabular}[t]{rrcl} 
    $\psi\colon$ & $D \times L$ & $\rightarrow$ & $D \times L$ \\
   $ $  & $(x,y)$ & $\mapsto$ & $(\psi^\prime(x,y),\psi^{\prime\prime}(x,y))$\, .
\end{tabular}
\end{center}

However, the bijectivity on the tangling half-planes implies that $\psi$ coarsely preserves the copies of $D$, that is the $a_ib_i$-horizontal lines. This means that, given the pair $(x_0,y_0)$ in the grid of $D \times (L \cup E_1)$, with image $\psi((x_0,y_o)) = (x_0^\prime,y_0^\prime)$, any other pair $(x,y_0)$ is mapped to $(x^\prime,y_0^\prime)$. This means that $\psi^{\prime\prime}(x,y)$ is independent of the input of $x$, hence we can interpret $\psi^{\prime\prime}$ as follows:

\begin{center}
\begin{tabular}[t]{rccl} 
    $\psi^{\prime\prime}\colon$ & $L$ & $\rightarrow$ & $L$ \\
   $ $  & $y$ & $\mapsto$ & $\psi^{\prime\prime}(y)$\, .
\end{tabular}
\end{center}

Via $\psi$ we extend $\psi^{\prime\prime}$ again to the tangling edges, that is we find map $\psi^{\prime\prime} \colon L \cup E_1 \rightarrow L \cup E_2$ that is bijective on the tangling edges.

Lastly note, that $\psi^{\prime\prime}$ is a QI. Indeed, given two pairs $(x_0,y_0)$,$(x_1,y_1)$ in the grid $D \times L$, we can decompose their distance as follows: $$d_{D \times L}((x_0,y_0),(x_1,y_1)) = d_D(x_0,x_1) + d_L(y_0,y_1) \, .$$
This implies that a QI-inequality for $\psi$ also holds for $\psi^{\prime\prime}$.

For convenience, we include a third reduction step. As in the proof of Claim \ref{ClaimReductionToTree} of Proposition \ref{VFD}, we contract every other edge of $L$. This way, we have $2 \, m_i + n_i$ edges at each vertex, removing the issue with the "half-edges".

So we have a QI between $V_1$ and $V_2$ which is bijective on the respective peripheral structures and which scales one copy of $D_\infty$ to the other by $\tfrac{2\, m_1 + n_1}{2 \, m_2 + n_2}$ if and only if we find a QI between two copies of the line $L$ with tangling edge set $E_1$ and $E_2$ with $2\, m_1 + n_1$ and $2 \, m_2 + n_2$ tangling edges at each vertex respectively that scales $L$ by $\tfrac{2\, m_1 + n_1}{2 \, m_2 + n_2}$ and is bijective on the tangling edges. 

\end{proof}

As Claim \ref{ReduceZwithTangling} suggests, we need to find a QI from $L \cup E_1$ to $L \cup E_2$ where the base graph $L$ is a line whose vertex set we can identify with $\mathbb{Z}$ and the QI is bijective on the tangling edges and scaling by $\tfrac{2\, m_1 + n_1}{2 \, m_2 + n_2}$. Without loss of generality we set $2 \, m_1 + n_1 > 2 \, m_2 + n_2$.

The first step is to define for $i \in \{1,2\}$ the following map
\begin{center}
\begin{tabular}[t]{rccl} 
    $\phi_i\colon$ & $L \cup E_i$ & $\rightarrow$ & $\mathbb{Z}$ \\
    $ $ & $l$ & $\mapsto$ & $l \cdot (2m_i + n_i)$  \\
   $ $  & $t_{z,j}$ & $\mapsto$ & $\begin{cases} 
     z \cdot (2 \, m_i + n_i) + j & \text{if } z \geq 0 \\ 
     (z+1)\cdot (2 \, m_i + n_i) - j - 1 & \text{if } z < 0 \, ,
   \end{cases}$     
\end{tabular}
\end{center}
where $t_{z,j}$ is one of the $2 \, m_i + n_i$ tangling edges at $z \in \mathbb{Z}$, i.e. $j \in \{0,\dots, 2 \, m_i + n_i-1\}$. It is easily checked that $\phi_i$ is bijective on $E_i$ for both $i \in \{1,2\}$ and by definition $\phi_i$ scales $L$ by $2 \, m_i + n_i$.

Now it suffices to show the following claim:
\begin{claim} \label{ClaimBijectiveQIId}
Any bijective QI $f\colon \mathbb{Z} \rightarrow \mathbb{Z}$ which fixes $0$, $\infty$ and $-\infty$ is bounded distance from the identity map.
\end{claim}

By using Claim \ref{ClaimBijectiveQIId}, for any isometry $i$ composed with $f$, we obtain the commuting diagram of the form:
\begin{center}
\begin{tabular}[t]{lcrclcr} 
     & $L$ & & $\overset{\varphi}{\longrightarrow}$ & & $L$ & \\
     $\restr{\phi_1}{L}$ & $ \bigg\downarrow$ & $\cdot (2\, m_1+n_1)$  & $ $ & $\restr{\phi_2}{L}$ & $ \bigg\uparrow$ & $: (2 \, m_2 + n_2)$  \\
     & $\mathbb{Z}$ & & $\overset{i \circ f}{\longrightarrow}$ &  & $\mathbb{Z}$ &
\end{tabular}
\end{center}

This implies that $\varphi$ is a QI scaling by $\frac{2\, m_1+n_1}{2 \, m_2 + n_2}$.\\

Thus, we are left to prove Claim \ref{ClaimBijectiveQIId}:

\begin{proof}(of Claim \ref{ClaimBijectiveQIId}) Let $f$ be a $(C,D)$-QI satisfying the assumptions and suppose it is not bounded distance from the identity. Then for any $n \in \mathbb{N}$ we can find a $z_n \in \mathbb{N}$ such that $d(z_n, f(z_n)) > n$.

First we claim that there is a maximal $k \in \mathbb{N}$ such that $f(-k) \geq 0$, implying by surjectivity of $f$ that $[0,\infty) \subseteq f([-k,\infty))$. Suppose this is not true. Then for every $k \in \mathbb{N}$ with $f(-k) \geq 0$ there is a $k^\prime \in \mathbb{N}$ such that $k< k^\prime$ and $f(-k^\prime) \geq 0$. However, by the QI-property and the fact that $f$ fixes $0$ we have
$$\frac{1}{C} \cdot d(-k,0) - D \leq d(f(-k),0) = f(-k) $$ for every $k \in \mathbb{N}$. Thus, with $k \in \mathbb{N}$ tending to $\infty$, so does $f(-k)$, in contradiction to the assumption that $f$ fixes $- \infty$.

Now let $B_R(z_n)$ be a ball of radius $R$ around $z_n$. We want to show that $$[0,f(z_n)+\frac{R}{C}-D] \subseteq f([-k,z_n+R])$$ for any $R \in \mathbb{N}$ large enough.

Since $f$ is bijective by assumption, some elements must map onto the interval $[0,f(z_n)+\frac{R}{C}-D]$. It is indeed $[-k,z_n+R]$ by the following observations illustrated in Figure \ref{fig:IllsutrationTanglingEdges} below:
\begin{enumerate}
    \item By choice of $k$, there is no element $k^\prime < -k$ such that $f(k^\prime) \geq 0$.
    \item Since $f$ is a QI, $B_{\frac{R}{C}-D}(f(z_n)) \subseteq f(B_R(z_n))$.
    \item Any element $z > z_n+R$ maps to an element $f(z) > f(z_n)+\frac{R}{C} -D$: Pick some $a > C \cdot f(z_n) + R$, for which $$f(a) = d(f(a),0) \geq \frac{a}{C} - D > f(z_n) + \frac{R}{C} - D \, .$$ Such an $a$ must exist, since $f$ fixes $0$, $\infty$ and $-\infty$. Thus $a$ is mapped to the right side of $Z\setminus B_{\frac{R}{C}-D}(f(z_n))$. Now choose $a^\prime$ such that $d(a,a^\prime) = 1$. This implies $$d(f(a),f(a^\prime)) \leq C +D \, .$$ If we choose $R \in \mathbb{N}$ such that $2(\frac{R}{C}-D) > C +D$, then $f(a)$ and $f(a^\prime)$ cannot be mapped to different sides of $Z\setminus B_{\frac{R}{C}-D}(f(z_n))$ and not in the ball. Thus they are both mapped to the right side. Now, for any arbitrary $z > z_n + R$, we pick a sequence $(a_i)_{i=0}^k$, where $a_0 = a$, $d(a_i,a_{i+1}) = 1$ for every $i \in \{0, \dots, k-1\}$ and $a_k = z$. Then all $f(a_i)$ with $i \in \{0,\dots, k\}$, in particular $f(z)$, must be on the right side of $Z\setminus B_{\frac{R}{C}-D}(f(z_n))$, that is $f(z) > f(z_n)+\frac{R}{C} -D$.
    
\end{enumerate}
\vspace{0.5cm}
\begin{figure}[ht]
    \centering
    \scriptsize
    \def\svgwidth{470pt}
    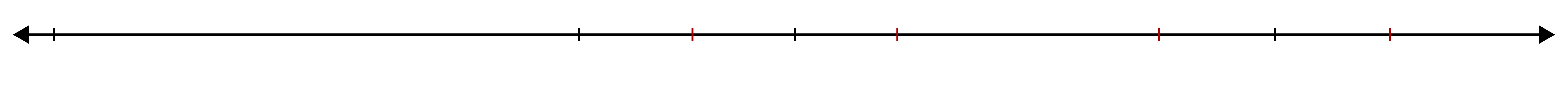
    \caption{The interval $[-k,z_n+R]$, illustrated in pink, maps onto the interval $[0,f(z_n)+\tfrac{R}{C}-D]$, illustrated in blue.}
   \label{fig:IllsutrationTanglingEdges}
\end{figure}
\vspace{0.5cm}
Hence, we ruled out all the elements outside of $[-k,z_n+R]$ to be mapped to $[0,f(z_n)+\frac{R}{C} - D]$, implying that $[0,f(z_n)+\frac{R}{C}-D] \subseteq f([-k,z_n+R])$. Thus, since $f$ is bijective, we obtain $$|[0,f(z_n)+\frac{R}{C}-D]| \leq |f([-k,z_n+R])|$$ and thus
\begin{align*}
    f(z_n) + \frac{R}{C}-D +1 &\leq z_n + R + k + 1 \\
    f(z_n) - z_n  &\leq (1-\frac{1}{C})R + k + D \, .
\end{align*}
But now choose $n > (1-\frac{1}{C})R + k + D$, then
$$  (1-\frac{1}{C})R + k + D < n \leq d(f(z_n),z_n) = f(z_n) -z_n \leq (1-\frac{1}{C})R + k + D \, , $$
which is a contradiction.
\end{proof}
Thus, to conclude, given two RACGs $W_1$ and $W_2$ with one cylinder vertex with VA vertex group, and one class of $e_1$ and $e_2$ indistinguishable non-cylinder vertices respectively, there is a QI between $W_1$ and $W_2$ and any such QI is bounded distance from scaling by $\frac{2 \, m_1 + n_1}{2 \, m_2 + n_2}$.
\end{proof}

If at the VA cylinder vertex attach $j \geq 2$ classes of indistinguishable neighbors we can apply Proposition \ref{1VACyl} to each class individually to obtain the following generalization:

\begin{cor} \label{VARatios}
Let $W_1$ and $W_2$ be two RACGs on defining graphs satisfying the \hyperref[StandingAssumption3]{Standing Assumption 3} with the same structure invariant and one single cylinder vertex $v_1$ and $v_2$ in the JSJ graph of cylinders $\Lambda_1$ and $\Lambda_2$ respectively. Let the cylinder vertex group $V_1$ and $V_2$ at $v_1$ and $v_2$ respectively be VA. Suppose at both $v_1$ and $v_2$ attach $j \geq 1$ classes of indistinguishable vertices respectively and let $e_{i,k} = m_{i,k} + n_{i,k}$ for $i \in \{1,2\}$ and $k \in \{1, \dots, j\}$ denote the number of neighbors $v_i$ in class $k$ with $m_{i,k}$ the number of neighbors attaching along a $D_\infty$-edge and $n_{i,k}$ the number of neighbors attaching along a $D_\infty \times \mathbb{Z}_2$-edge. If there is a QI between $W_1$ and $W_2$, then the ratio $\frac{2 \, m_{1,k} + n_{1,k}}{2 \, m_{2,k} + n_{2,k}}$ is the same for all $k \in \{1,\dots,j\}$.
\end{cor}

\begin{exam} \label{ExamplVAnotQI}
In Example \ref{ExampleStructureInvar} illustrated by Figure \ref{fig:ExampleStructureInvar}, we have $j = 2$ classes of indistinguishable tangling edges. Since all occurring edge groups are $D_\infty$, we have $e_{i,1} = m_{i,1}$ counting the hanging vertices and $e_{i,2} = m_{i,2}$ counting the rigid. Then we have 
\begin{center}
\begin{tabular}[t]{c|cc} 
   $m_{i,k} $ & $k=1$ & $k=2$ \\\hline
   $i=1$ & 1 & 1 \\
   $i=2$ & 2 & 1 \\\hline
   ratio & $\frac{1}{2}$ & $1$
\end{tabular}
\end{center}
implying that by Corollary \ref{VARatios} the VA cylinder vertices of $W_1$ and $W_2$ don't have the same relative QI-type and thus $W_1$ and $W_2$ are not QI by Proposition \ref{decopresTreeIsom}.
\end{exam}

\subsection{Refinement of the structure invariant} \label{SubsectionDensityRef}
In Corollary \ref{VARatios}, we have seen that the number of neighbors per class of indistinguishable vertices at a VA cylinder vertex in the JSJ graph of cylinders is an essential characteristic to determine whether two groups are QI or not. Thus, we aim to alter the structure invariant in a way such that this information is taken into account. For that purpose we introduce a process we call \textit{density refinement}.

\begin{constr} \label{ValencyRefin}
We start with an initial decoration $\delta_0$ with an initial set of ornaments consisting of the vertex and the relative QI-type. We perform the neighbor refinement, giving us a stable decoration $\delta_i$.

Now we define the map $\nu_i\colon V(T) \rightarrow \mathbb{N}^{\delta_i(V(T))}\bign/\sim \cup \, \{\#\}$, where $\sim$ is an equivalence relation defined in Step \ref{DefOfDensitiyRefStep2} below, as follows:
\begin{itemize}
    \item For any vertex $v \in V(T)$, whose vertex group is not VA, $\nu_i$ maps $v$ to $\#$.
    \item A vertex $v \in V(T)$, whose vertex group is VA, is mapped to an equivalence class of tuples with entries in $\mathbb{N}$ indexed by the image of the decoration $\delta_i$. We obtain the image $\nu_i(v)$ in two steps:
    \begin{enumerate}
        \item We associate to $v$ a tuple $\alpha$ obtained as follows: The entry indexed by $o \in \mathcal{O}_i$ is computed from the JSJ graph of cylinders $\Lambda_c$. We look at the neighbors of the vertex in $\Lambda_c$ corresponding to the orbit of $v$ with ornament $o$. Let $m$ be the number of such neighbors attached along a $D_\infty$-edge and $n$ be the number of such neighbors attached along a $D_\infty \times \mathbb{Z}_2$-edge. Then the entry is $2 \, m + n$.
        \item \label{DefOfDensitiyRefStep2} Define the image of $v$ under $\nu_i$ as the \textit{projective class} of $\alpha$, that is the equivalence class under the relation: $\alpha \sim \beta \text{ if and only if there is a } k \in \mathbb{R}^+ \text{ such that } k \cdot \alpha = \beta \, ,$ where the multiplication $\cdot$ is defined coordinate-wise.
    \end{enumerate}
\end{itemize}

With the map $\nu_i$ we provide a new decoration: The new set of ornaments is $$\mathcal{O}_i^\prime := \mathcal{O}_0 \times \mathbb{N}^{\delta_i(V(T))}\bign/\sim \cup \, \{\#\} \times \bar{\mathbb{N}}^{\mathcal{O}_i}$$ and the decoration is $\delta_i^\prime\colon T \rightarrow \mathcal{O}_i^\prime$ with $$\delta_i^\prime(v) := (\delta_0(v), \nu_i(v),f_{v,i})$$ for any $v \in V(T)$. Possibly, $\delta_i^\prime$ is a refinement of $\delta_i$ and thus we can perform the neighbor refinement on it. Again, we obtain a stable decoration $\delta_j$ for which we can define a map $\nu_j$ as above. We define a new set of ornaments $\mathcal{O}_j^\prime := \mathcal{O}_0 \times \mathbb{N}^{\delta_j(V(T))}\bign/\sim \cup \, \{\#\} \times \bar{\mathbb{N}}^{\mathcal{O}_j}$ and the decoration $\delta_j^\prime\colon T \rightarrow \mathcal{O}_j^\prime$ with $\delta_j^\prime(v) := (\delta_0(v), \nu_j(v),f_{v,j})$ for any $v \in V(T)$. We repeat this alternating refinement process. Since there are only finitely many cylinder vertices in $\Lambda_c$, this process will eventually stabilize. The resulting decoration is the \textit{density refinement} of $\delta_0$.
\end{constr}

Combining Proposition \ref{decopresTreeIsom} and Corollary \ref{VARatios} yields that two RACGs can only be QI if their structure invariants, where $\delta_s$ is stable with respect to the density refinement, are identical.

\begin{exam} \label{ExRefinement}
The original structure invariant for the group illustrated in Figure \ref{fig:ValencyRefin} with respect to only the neighbor refinement is illustrated in the following table: 
\begin{table}[H]
\footnotesize
\centering
\begin{tabular}[H]{r|c|c||ccc|cc|c} 
 & \makecell[b]{vertex\\type} & QI type & \makecell[b]{$c_1$\\$c_4$\\$c_6$} & $c_3$ & \makecell[b]{$c_2$\\$c_5$} & \makecell[b]{$h_1$\\$h_3$\\$h_4$} & $h_2$ & $r$ \\\hhline{=|=|=||===|==|=}
 $c_1,c_4,c_6$ & ‘cyl' & 2-ended & 0 & 0 & 0 & 0 & 0 & 1 \\
 $c_3$ & ‘cyl' & 2-ended & 0 & 0 & 0 & 0 & 1 & 1 \\
 $c_2,c_5$ & ‘cyl' & ‘VA' & 0 & 0 & 0 & $\infty$ & 0 & $\infty$ \\\hline
$h_1,h_3,h_4$ & ‘hang' & ‘VF' & 0 & 0 & $\infty$ & 0 & 0 & 0 \\
$h_2$ & ‘hang' & ‘VF' & 0 & $\infty$ & 0 & 0 & 0 & 0 \\\hline
$r$ & ‘rig' & & $\infty$ & $\infty$ & $\infty$ & 0 & 0 & 0 
\end{tabular}
\end{table}

We see that the vertices $c_2$ and $c_5$ are indistinguishable. However, when performing the density refinement according to Construction \ref{ValencyRefin}, the images of $c_2$ and $c_5$ under $\nu_i$ differ: $$\nu_i(c_2) = [(0,0,0,2,0,2)] \quad \text{ and } \quad  \nu_i(c_5) = [(0,0,0,4,0,2)] \, .$$ This makes it possible to further distinguish $h_1$ from $h_3$ and $h_4$. We obtain the following refined structure invariant:

\begin{table}[H]
\footnotesize
\centering
\begin{tabular}[H]{r|c|c|c||cccc|ccc|c} 
 & \makecell[b]{vertex\\type} & QI type & $\nu_{stable}$ & \makecell[b]{$c_1$\\$c_4$\\$c_6$} & $c_3$ & \makecell[b]{$c_2$} & $c_5$ & $h_1$ & \makecell[b]{$h_3$\\$h_4$} & $h_2$ & $r$ \\\hhline{=|=|=|=||====|===|=}
 $c_1,c_4,c_6$ & ‘cyl' & 2-ended & $\#$ & 0 & 0 & 0 & 0 & 0 & 0 & 0 & 1 \\
 $c_3$ & ‘cyl' & 2-ended & $\#$ & 0 & 0 & 0 & 0 & 0 & 0 & 1 & 1 \\
 $c_2$ & ‘cyl' & ‘VA' & $[(0,0,0,0,1,0,0,1)]$ & 0 & 0 & 0 & 0 & $\infty$ & 0 & 0 & $\infty$ \\
 $c_5$ & ‘cyl' & ‘VA' & $[(0,0,0,0,0,2,0,1)]$ & 0 & 0 & 0 & 0 & 0 & $\infty$ & 0 & $\infty$ \\\hline
$h_1$ & ‘hang' & ‘VF' & $\#$ & 0 & 0 & $\infty$ & 0 & 0 & 0 & 0 & 0 \\
$h_3,h_4$ & ‘hang' & ‘VF'& $\#$  & 0 & 0 & 0 & $\infty$ & 0 & 0 & 0 & 0 \\
$h_2$ & ‘hang' & ‘VF' & $\#$ & 0 & $\infty$ & 0 & 0 & 0 & 0 & 0 & 0 \\\hline
$r$ & ‘rig' & & $\#$ & $\infty$ & $\infty$ & $\infty$ & $\infty$ & 0 & 0 & 0 & 0
\end{tabular}
\end{table}

\begin{figure}[ht]
    \centering
    \footnotesize
    \def\svgwidth{480pt}
    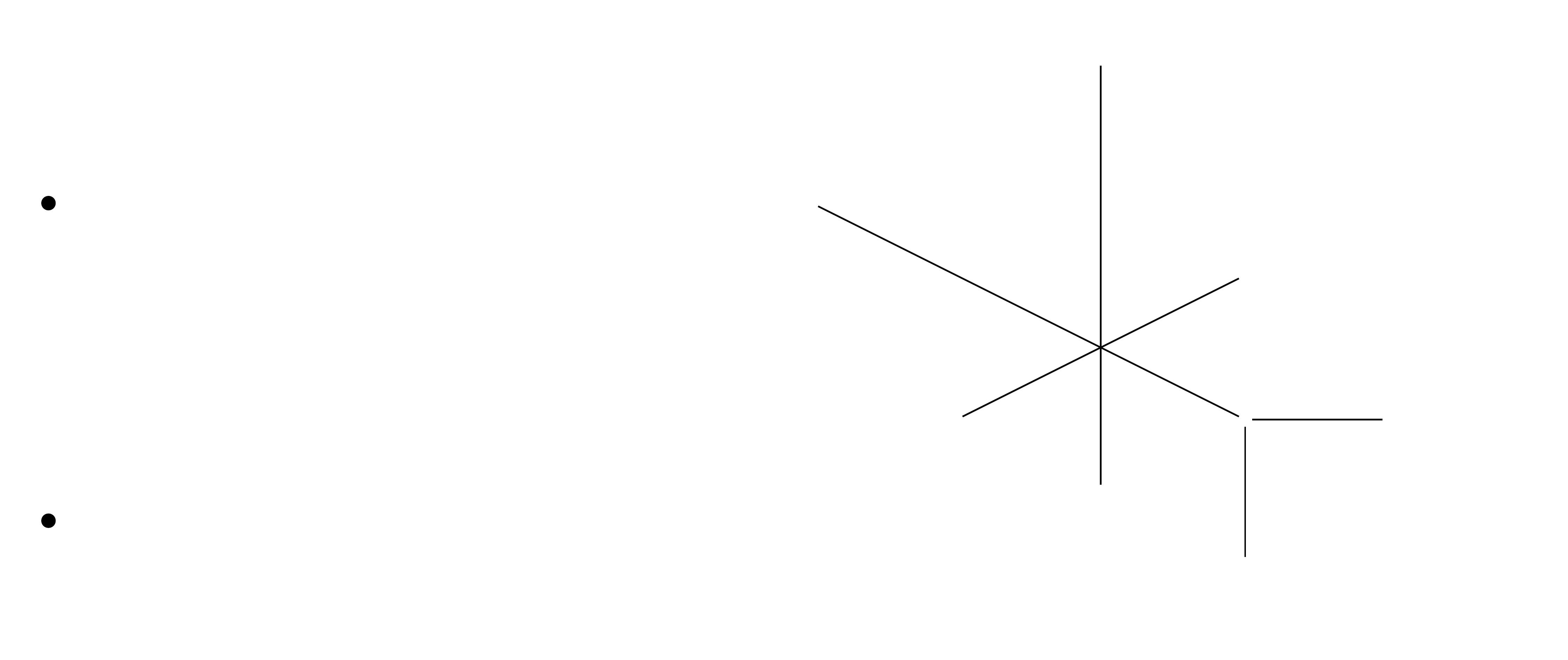
    \caption{$\Lambda_c$ is the JSJ graph of cylinders of the RACG $W_\Gamma$.}
   \label{fig:ValencyRefin}
\end{figure}
\end{exam}

\subsection{Complete QI-Invariant} \label{SectionCompleteQI}
Now we aim to put the local QIs between cylinder vertex groups together to obtain a global QI between the groups and thus have a structure invariant which is a complete QI-invariant for certain groups. As mentioned in Remark \ref{technicalities}, we exclude rigid vertices so the only missing piece are the local QIs between the hanging vertices. We see that we can choose them with a lot of flexibility: 

\begin{thm}\label{finalresult}
Let $W$ and $W^\prime$ be two finitely presented, one-ended RACGs with non-trivial JSJ decompositions over two-ended subgroups, which both have no rigid vertices. Define $T$ to be the JSJ tree of cylinders of $W$ and $X$ to be the geometric tree of spaces of $W$ over $T$. The initial decoration $\delta_0$ on $T$ takes vertex type and relative QI-type into account. Let $\delta$ be the density refinement of $\delta_0$. Analogously, we define $T^\prime$, $X^\prime$, $\delta_0^\prime$ and $\delta^\prime$ for $W^\prime$. Then $W$ and $W^\prime$ are QI if and only if there is a bijection $\beta\colon \delta(T) \rightarrow \delta^\prime(T^\prime)$ such that
\begin{enumerate}
    \item $\delta_0 \circ \delta^{-1} = \delta_0^\prime \circ (\delta^\prime)^{-1} \circ \beta$;
    \item $S(T,\delta,\mathcal{O}) = S(T^\prime, \delta^\prime, \mathcal{O}^\prime)$ in the $\beta$-induced ordering;
    \item for every ornament $o \in \mathcal{O}$, there is a vertex $v \in \delta^{-1}(o)$ and a vertex $v^\prime \in (\delta^\prime)^{-1}(\beta(o))$ such that there is a QI between the vertex spaces $X_v$ and $X^\prime_{v^\prime}$ respecting the decorations $\delta$ and $\delta^\prime$ and which is bijective on the peripheral structures $\mathcal{P}_v$ and $\mathcal{P}^\prime_{v^\prime}$  respectively.
\end{enumerate}
\end{thm}

\begin{sproof}
This is an analogue of the proof of \cite[Theorem 7.5]{CashenMartinStructureInvar}, with some generalizations and some specializations. The statement is more specialized in the two aspects laid out in Remark \ref{technicalities}: We assume that the considered groups do not have any rigid vertices, thus the relative stretch factors do not apply. Moreover, since we restrict to RACGs, partial orientations can be omitted. However, we do not assume the cylinder vertex groups to be two-ended, which makes the statement more general. 

The idea is to inductively build a tree isometry $\chi\colon T \rightarrow T^\prime$, which respects the decorations by using the local vertex QIs $\phi_v\colon X_v \rightarrow X^\prime_{\chi(v)}$ bijective on the respective peripheral structures inducing $\chi$ on the link of the vertex $v \in V(T)$. Then $\chi$ induces a global QI.

For the base case, we pick some cylinder vertex $c \in V(T)$ and some $c^\prime \in (\delta^\prime)^{-1}(\beta(\delta(c)))$ and define $\chi(c) := c^\prime$. Because the initial decoration depends on the relative QI-type, there is a QI between $X_c$ and $X^\prime_{c^\prime}$. Depending on whether $c$ has a two-ended, a VFD or a VA vertex group, we pick such a QI $\phi_c: X_c \rightarrow X^\prime_{c^\prime}$ according to Propositions \ref{CompleteInvarVZcylinders}, \ref{VFD} and \ref{1VACyl} respectively. By construction, $\phi_c$ will be bijective on the respective peripheral structures and thus defines how to pick the bijection between the edge spaces incident to $c$. Thus, we can extend $\chi$ to the link of $c$ according to this bijection.

Since the considered trees are bipartite, the inductive step consists of two parts: First we extend $\chi$ to a hanging vertex and from there we extend $\chi$ to a cylinder vertex.

Suppose there is an edge $e_1\in E(T)$ such that $o(e_1)$ is a cylinder vertex, $\tau(e_1) =: h$ is a hanging vertex and $\chi(o(e_1))$ is already defined. Then there is a QI $\phi_{o(e_1)}\colon X_{o(e_1)} \rightarrow X^\prime_{\chi(o(e_1))}$ respecting the decorations and bijective on the respective peripheral structures. Thus $\phi_{o(e_1)}|_{X_{e_1}}\colon X_{e_1} \rightarrow X^\prime_{\chi(e_1)}$ defines the QI on $X_{e_1}$. The QI on $X_h$ can now be produced as suggested in \cite[Proposition 7.1]{CashenMartinStructureInvar}, which is guided by \cite[Theorem 1.2]{behrstockNeumann}. The key feature is the following: Pick for any other edge $e$ adjacent to $h$ some real constant $\sigma_e$. The only condition is that for all edges in the same orbit the constant needs to be identical. Then we can choose a QI $\phi_{h}: X_{h} \rightarrow X^\prime_{\chi(h)}$ such that when restricted to $X_{e_1}$ it matches $\phi_{o(e_1)}|_{X_{e_1}}$ and when restricted to $X_e$ for any other $e$ adjacent to $h$, this $\phi_{h}|_{X_e}$ is a QI with multiplicative constant $\sigma_e$.

Of course, we do not pick the $\sigma_e$ randomly, but we choose them among the set $\Sigma$ of multiplicative constants occurring in the QIs produced by the Propositions \ref{CompleteInvarVZcylinders}, \ref{VFD} and \ref{1VACyl}. Since there are only finitely many orbits of cylinder vertices, this set $\Sigma$ is finite and we also only pick a finite configuration of $\sigma_e$'s from $\Sigma$. If we later see that our choice of configuration conflicts with the constants forced by the QIs of the adjacent cylinder vertices, we return to $h$ and pick a different configuration. Since the number of such different configurations is finite, we know that eventually we have found the correct QI and extend $\chi$ to the link of $h$ accordingly. Thus without loss of generality we can assume that we have picked a suitable QI at $h$ satisfying all requirements.

Suppose now that $e_2 \in E(T)$ is an edge such that $o(e_2)$ is a hanging vertex, $\tau(e_2) = c_2$ is a cylinder vertex and $\chi(o(e_2))$ is already defined in the previous step. We repeat the extension process: We know that there is a QI $\phi_{o(e_2)}\colon X_{o(e_2)} \rightarrow X^\prime_{\chi(o(e_2))}$ respecting decorations and bijective on the respective peripheral structures, which restricts to a QI on $X_{e_2}$. We can now extend $\chi$ to the link of $c_2$ and define $\phi_{c_2}\colon X_{c_2} \rightarrow X_{\chi(c_2)}$ according to the Propositions \ref{CompleteInvarVZcylinders}, \ref{VFD} and \ref{1VACyl} such that it agrees with $\phi_{o(e_2)}$ on $X_{e_2}$. 
\end{sproof}

\begin{exam}
We make the introductory example of the two groups with defining graphs illustrated in Figure \ref{fig:EasyQIEx} explicit. By Proposition \ref{VFD}, we see that the VFD cylinder vertices coming from the blue uncrossed cut pairs are QI. In both cases there is one hanging vertex group generated by the $l_i$'s, thus by Proposition \ref{1VACyl}, the VA cylinder vertices coming from the red uncrossed cut pairs are QI. Hence the structure invariants are identical and Theorem \ref{finalresult} implies that the groups are QI.

\begin{figure}[ht]
    \centering
    \scriptsize
    \def\svgwidth{350pt}
    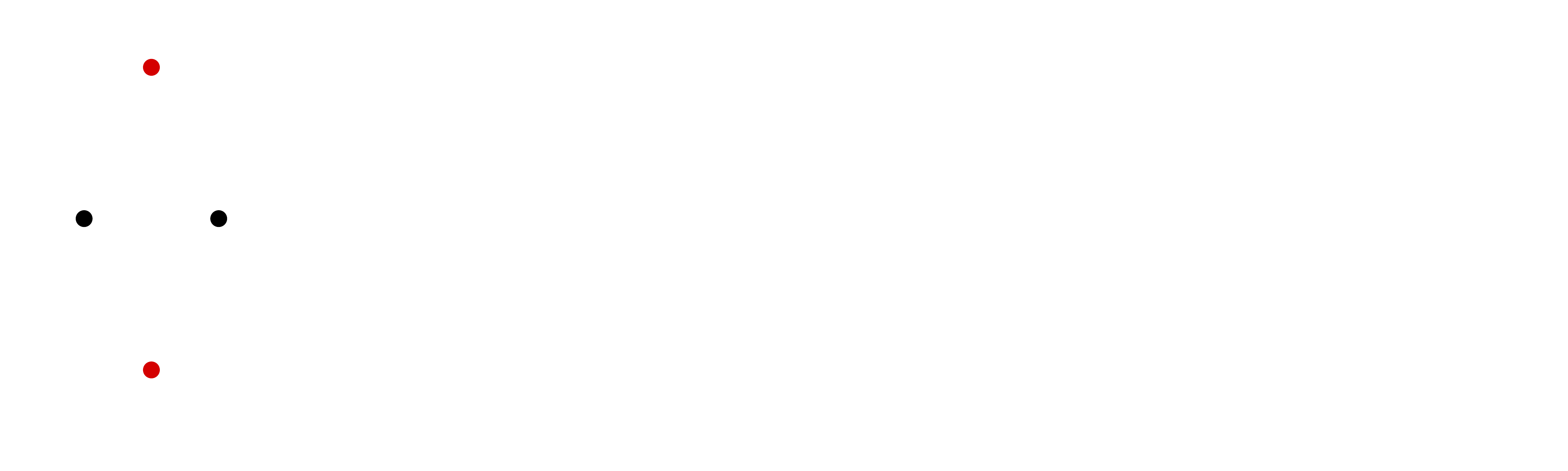
    \caption{The RACGs on the graphs $\Gamma_1$ and $\Gamma_2$ are QI to each other.}
   \label{fig:EasyQIEx}
\end{figure}
\end{exam}

\begin{rem} \label{RemWithRig}
It is discussed in Remark \ref{technicalities} that Theorem \ref{finalresult} excludes groups whose JSJ decompositions have rigid vertices. However, in certain cases we can add another induction step to the proof of Theorem \ref{finalresult} handling rigid vertices following again the proof of Theorem 7.5 of \cite{CashenMartinStructureInvar}. For instance, we can consider the subgraphs $\Lambda$ and $\Lambda^\prime$ of the graphs of cylinders $\Lambda_c$ and $\Lambda_c^\prime$ respectively, which consist of one rigid vertex and all its adjacent cylinder vertices. If there is a decoration preserving graph isomorphism $\phi$ between $\Lambda$ and $\Lambda^\prime$ and in addition, every vertex and edge group $G_t$ is isomorphic to the image vertex group $G_{\phi(t)}$, then the induction extends also to these rigid vertices. The obvious method to produce such an example is to simply use identical defining graphs for corresponding special subgroups. This is illustrated in the following Example \ref{ExQI2}.

Alternatively, if the rigid vertices are virtually free, they are quasi-isometrically rigid relative to the peripheral structure by \cite{cashenmacura} and thus relative stretch factors can be used as introduced in Section 4 of \cite{CashenMartinStructureInvar}.
\end{rem}

\begin{out} \label{methodQIEx}
Theorem \ref{finalresult} illustrates the flexibility we have to change the defining graph in a way such that the group on the resulting graph is QI to the one on the original graph. The changes happen at the cylinder vertices:
\begin{itemize}[noitemsep, topsep=2pt]
    \item At a virtually cyclic cylinder vertex coming from the uncrossed cut collection $\{a-b\}$ we can only remove or add a common adjacent vertex such that $|\mathcal{C}|\in \{0,1\}$ is maintained since the valencies in the JSJ tree of cylinders need to be preserved.
    \item At a VFD cylinder vertex coming from the uncrossed cut collection $\{a-b\}$ and its common adjacent vertices $\mathcal{C}$, we can duplicate or remove tangling pieces in the JSJ tree of cylinders that are equivalent up to QI (cf.\ Algorithm \ref{AlgVFree}). In the defining graph $\Gamma$ this corresponds to duplicating or removing any connected component of $\Gamma \setminus \{a-b\}$ disjoint from $\mathcal{C}$, and reattaching the new collection of pieces to $a$ and $b$. Note that in the reattaching the roles of $a$ and $b$ can be interchanged, thus this move can be interpreted as a reflection along the subgraph on $\{a,b\} \cup \mathcal{C}$. Additionally the number of vertices in $\mathcal{C}$ can be changed. The only restriction is that $|\mathcal{C}|>2$.
    
    Also, within a connected component containing vertices of a set $A$ contributing to a hanging vertex, the number of vertices can be altered while preserving the virtually free QI-type. That means, we can add or remove elements on a branch, as long as the resulting vertex set $A$ still produces a hanging vertex. Thus by Proposition \ref{Charhanging}, the altered set $A$ should still satisfy conditions \hyperref[A1n]{(A1)}, \hyperref[A2n]{(A2)} and \hyperref[A3n]{(A3)} and $W_A$ has to be infinite and not a cylinder vertex group.
    
    \item At a VA cylinder vertex coming from the uncrossed cut collection $\{a-b\}$ and its common adjacent vertices $\mathcal{C}$ we can perform changes similar to the ones at VFD cylinder vertices. There are only two differences: We perform the duplication or removal of pieces with a fixed ratio and the number of common adjacent vertices has to stay fixed $|\mathcal{C}|=2$.
\end{itemize}
\end{out}

These observations can be used as a method to produce examples of QI RACGs:

\begin{exam} \label{ExQI1}
The RACGs on the defining graphs $\Gamma_1$ and $\Gamma_1^\prime$ with JSJ graphs of cylinders $\Lambda_{c,1}$ and  $\Lambda_{c,1}^\prime$ respectively, illustrated in Figure \ref{fig:ExQInonRig} are QI by Theorem \ref{finalresult}.
\end{exam}

\begin{exam} \label{ExQI2}
The RACGs on the defining graphs $\Gamma_2$ and $\Gamma_2^\prime$ with JSJ graphs of cylinders $\Lambda_{c,2}$ and  $\Lambda_{c,2}^\prime$ respectively, illustrated in Figure \ref{fig:ExQIwithRig} are QI by Theorem \ref{finalresult} and Remark \ref{RemWithRig}.
\end{exam}

\begin{rem}
It would be most interesting to produce QIs that do not arise from algebraic considerations. One might guess that a simple graph operation like duplicating the complement of a subgroup corresponding to a cylinder vertex would produce either a group which is a finite index subgroup of the original one or at least produce a group which shares a common finite index subgroup with it. In this case we call the groups (abstractly) commensurable and this already implies that they are QI. However, our construction has much more flexibility than that. 

Only partial commensurability results are known, such as the commensurability classification for certain hyperbolic RACGs done by Dani, Stark and Thomas in \cite{dani2018commensurability}. Their proof is not applicable to our more general setting, as it strongly depends on the fact that the finite valence of cylinder vertices in the JSJ tree of cylinders of hyperbolic RACGs is a QI-invariant. This tool is lost for non-hyperbolic RACGs. In \cite[Section 4]{hruska2020surface}, Hruska, Stark and Tran provide examples of commensurable non-hyperbolic RACGs whose defining graphs are generalized theta graphs. However, a complete classification for some class of non-hyperbolic RACGs is yet to be stated and should be addressed separately. Nonetheless, we can show that the non-hyperbolic Examples \ref{ExQI1} and \ref{ExQI2} for which we produced QIs with our methods are not abstractly commensurable, by application of the following Lemma \ref{LemCommens}, which is guided by Lemma 7.2 of Shepherd and Woodhouse \cite{shepherd2020quasi}.
\end{rem}

\begin{lem} \label{LemCommens}
The two RACGs $W_1$ and $W_1^\prime$ in Example \ref{ExQI1} on the defining graphs $\Gamma_1$ and $\Gamma_1^\prime$ in Figure \ref{fig:ExQInonRig} are not commensurable to each other and the  two RACGs $W_2$ and $W_2^\prime$ in Example \ref{ExQI2} on the defining graphs $\Gamma_2$ and $\Gamma_2^\prime$ in Figure \ref{fig:ExQIwithRig} are not commensurable to each other.
\end{lem}

\begin{proof}\cite[cf.][Lemma 7.2]{shepherd2020quasi}
Let $W$ and $W^\prime$ be two RACGs whose JSJ graphs of cylinders have cylinder vertices $v$ and $v^\prime$ with vertex groups $W_{\mathcal{C}} \times D_\infty$ and $W_{\mathcal{C}^\prime} \times D_\infty$ respectively such that $W_{\mathcal{C}}$ and $W_{\mathcal{C}^\prime}$ are both virtually free, i.e. $|\mathcal{C}|, |\mathcal{C}^\prime| > 2$. In fact, given $\mathcal{C} = \{c_1, \dots, c_{i+1}\}$, as per the proof of Theorem B.1 of Cashen, Dani and Thomas in \cite[Appendix B]{DaniThomasBowditch}, $W_{\mathcal{C}}$ has a free subgroup $F_i$ generated by $\langle c_1c_2, \dots, c_1c_{i+1} \rangle$ of rank $i$ and index 2. Analogously $W_{\mathcal{C}^\prime}$ has a free subgroup $F_j$ of index 2 and rank $|\mathcal{C}^\prime|-1 =: j$.

Suppose that $W$ and $W^\prime$ are commensurable, that is they have isomorphic finite index subgroups. By \cite[Corollary 7.4]{JSJDecompGps} we can assume that the induced JSJ graphs of cylinders of these subgroups are identical. Call this induced JSJ graph of cylinders $\hat{\Gamma}$ with fundamental group $\hat{W}$. The idea is now to compute the degree of a vertex in $\hat{\Gamma}$, using first $W$ and then $W^\prime$ and obtain a contradiction for the groups we are interested in as the computed degrees cannot match.

Suppose there is a $\hat{v} \in V(\hat{\Gamma})$ with vertex group $\hat{G}_{\hat{v}}$ covering $v$ and $v^\prime$. Then we can embed $\hat{G}_{\hat{v}}$ into both $W_{\mathcal{C}} \times D_\infty$ and $W_{\mathcal{C}^\prime} \times D_\infty$ as a finite index subgroup.
Note that such a vertex $\hat{v}$ exists in both the examples we consider here: $\Lambda_{c,1}$ has only one VFD cylinder vertex, the vertex $c_2$. Thus any vertex $\hat{v}$ covering some VFD cylinder vertex in $\Lambda_{c,1}^\prime$ has to cover $c_2$ as well. This argument works also for $\Lambda_{c,2}$ with its only VFD vertex $c_3$ and $\Lambda_{c,2}^\prime$.

Moreover, in the considered examples all edge groups are the same $D_\infty$ generated by the cut pair. Hence the number of edges incident to $\hat{v}$ corresponds to the number of double cosets $\hat{G}_{\hat{v}} g D_\infty$ with $g$ an element in the cylinder vertex group, multiplied by the degree of the cylinder vertex. So we aim to compute $\deg(\hat{v})$ in two ways, first via $v$, then via $v^\prime$:
\begin{align*}
    \deg(\hat{v}) &= |\{\hat{G}_{\hat{v}} g D_\infty\mid g \in W_{\mathcal{C}} \times D_\infty\}| \cdot \deg(v) \\
                  &= |\{\hat{G}_{\hat{v}} g D_\infty\mid g \in W_{\mathcal{C}^\prime} \times D_\infty\}| \cdot \deg(v^\prime)
\end{align*}

In order to do this we consider $F_i \times D_\infty \leq W_{\mathcal{C}} \times D_\infty$. This is a subgroup of index 2, thus the intersection $G_i := F_i \times D_\infty \cap \hat{G}_{\hat{v}} \leq \hat{G}_{\hat{v}}$ is at most of index 2 in $\hat{G}_{\hat{v}}$. For $F_j \times D_\infty \leq W_{\mathcal{C}^\prime} \times D_\infty$ we define analogously $G_j := F_j \times D_\infty \cap \hat{G}_{\hat{v}} \leq \hat{G}_{\hat{v}}$, which is also at most of index 2 in $\hat{G}_{\hat{v}}$. Hence we have for the intersection $G := G_i \cap G_j$
$$|\hat{G}_{\hat{v}}:G| \leq |\hat{G}_{\hat{v}}:G_i| \, |\hat{G}_{\hat{v}}:G_j| \leq 2 \cdot 2 = 4 \, .$$

Now we decompose the double cosets $\hat{G}_{\hat{v}} g D_\infty$ further into double cosets of $G$ with representatives $f_i$ in $F_i \times D_\infty$. Since $D_\infty$ is central, it suffices to consider $\hat{G}_{\hat{v}} g$: Each such coset $\hat{G}_{\hat{v}} g$ consists of at most 8 cosets of the form $Gf_i$. Indeed, at most 4 cosets come from the partition of $\hat{G}_{\hat{v}}$ into $G$-cosets as the index of $G$ in $\hat{G}_{\hat{v}}$ is at most 4 and then we multiply by 2 because $F_i \times D_\infty$ is of index 2 in $W_{\mathcal{C}} \times D_\infty$. This bounds the number of double cosets $\hat{G}_{\hat{v}} g D_\infty$ by:
$$\frac{1}{8} \, |\{Gf_iD_\infty\mid f_i \in F_i \times D_\infty\}| \leq |\{\hat{G}_{\hat{v}} g D_\infty\mid g \in W_{\mathcal{C}} \times D_\infty\}| \leq 2 \, |\{Gf_iD_\infty\mid f_i \in F_i \times D_\infty\}| \, .$$

Let $\pi_i: F_i \times D_\infty \rightarrow F_i$ be the projection map. Then the image $\pi_i(G)$ is a subgroup of $F_i$ and thus free. This implies that the short exact sequence
$$1 \rightarrow G \cap \ker(\pi_i) \rightarrow G \rightarrow \pi_i(G) \rightarrow 1$$
splits, that is there is a section $\sigma_i: \pi_i(G) \rightarrow G$ with image $P_i$ isomorphic to $\pi_i(G)$. But since $D_\infty$ is central in $F_i \times D_\infty$, we know that $G = P_i \times (G \cap \ker(\pi_i))$. Thus the number of double cosets $G f_i D_\infty$ is equal to the number of cosets $\pi_i(G)\pi_i(f_i) \times D_\infty$ in $F_i \times D_\infty$. But this number is the index of $\pi_i(G)$ in $F_i$, which we compute with the Schreier-index formula
$$|F_i : \pi_i(G)| = \frac{\mathrm{rk}(\pi_i(G))-1}{i-1} \, .$$
Analogously, we perform the same argument for $F_j \times D_\infty$ and the projection map $\pi_j: F_j \times D_\infty \rightarrow F_j$ to compute the number of double cosets $Gf_jD_\infty$ with $f_j \in F_j \times D_\infty$ via
$$|F_j : \pi_j(G)| = \frac{\mathrm{rk}(\pi_j(G))-1}{j-1} \, .$$
However, we note that $$\mathrm{rk}(\pi_i(G))=\mathrm{rk}(G/\ker(\pi_i)) = \mathrm{rk}(G/\mathrm{Z}(G)) = \mathrm{rk}(G/\ker(\pi_j)) = \mathrm{rk}(\pi_j(G)) \, ,$$ that is both occurring ranks are identical, call them $r$.
Thus, when computing $\deg(\hat{v})$ via $v$, we can use the first computation to obtain the bound
$$\frac{1}{8} \, \frac{r-1}{i-1} \cdot \deg(v)  \leq \deg(\hat{v}) \leq 2 \, \frac{r-1}{i-1} \cdot \deg(v) \, .$$
When computing $\deg(\hat{v})$ via $v^\prime$, we obtain using the second computation
$$\frac{1}{8} \, \frac{r-1}{j-1} \cdot \deg(v^\prime) \leq \deg(\hat{v})  \leq 2 \, \frac{r-1}{j-1} \cdot \deg(v^\prime) \, .$$
This implies that we arrive at a contradiction, whenever 
$$2 \, \frac{r-1}{j-1}  \cdot \deg(v^\prime) < \frac{1}{8} \, \frac{r-1}{i-1} \cdot \deg(v)$$
that is, whenever
$$j > 16 \cdot (i-1) \cdot \frac{\deg(v^\prime)}{\deg(v)} + 1 \, .$$
In case of Example \ref{ExQI1} this inequality is satisfied for the two VFD cylinder vertices labelled $c_2$ and $c_2^\prime$: In $\Gamma_1$, we have $|\mathcal{C}| = 3$, thus $i = 2$ and $\deg(c_2)=1$. In $\Gamma_1^\prime$ we have $|\mathcal{C}^\prime| = 35$, thus $j = 34$ and $\deg(c_2^\prime) = 2$.
In Example \ref{ExQI2}, the condition is satisfied for the vertices labelled $c_3$. Hence $W_1$ and $W_1^\prime$ in Example \ref{ExQI1} and $W_2$ and $W_2^\prime$ in Example \ref{ExQI2} are not commensurable to each other.
\end{proof}

\begin{rem}
The proof of Lemma \ref{LemCommens} works for various other examples. In fact, it can even provide a more sensitive commensurability invariant. Recall that the argument involves computing for an edge $e$ with edge group $D_\infty$ at the vertex $v$ the number of cosets $|\{\hat{G}_{\hat{v}} g D_\infty\mid g \in W_{\mathcal{C}} \times D_\infty\}|$. Then we sum over all such edges $e$, which is the same as multiplying by the degree $\deg(v)$ of $v$. 

However, instead of summing over all edges $e$ incident to $v$, we can restrict to a certain subclass of edges. For example we can restrict to edges whose incident vertices share the same vertex types, because the vertex types of the incident vertices of a covering edge must be the same as the ones of the covered edge. Even finer than just considering the vertex type would be to restrict to edges with incident vertices sharing a particular decoration which has to be preserved by the covering.
\end{rem}

\begin{figure}[H]
    \begin{subfigure}[c]{\textwidth}
        \centering
        \scriptsize
        \def\svgwidth{380pt}
        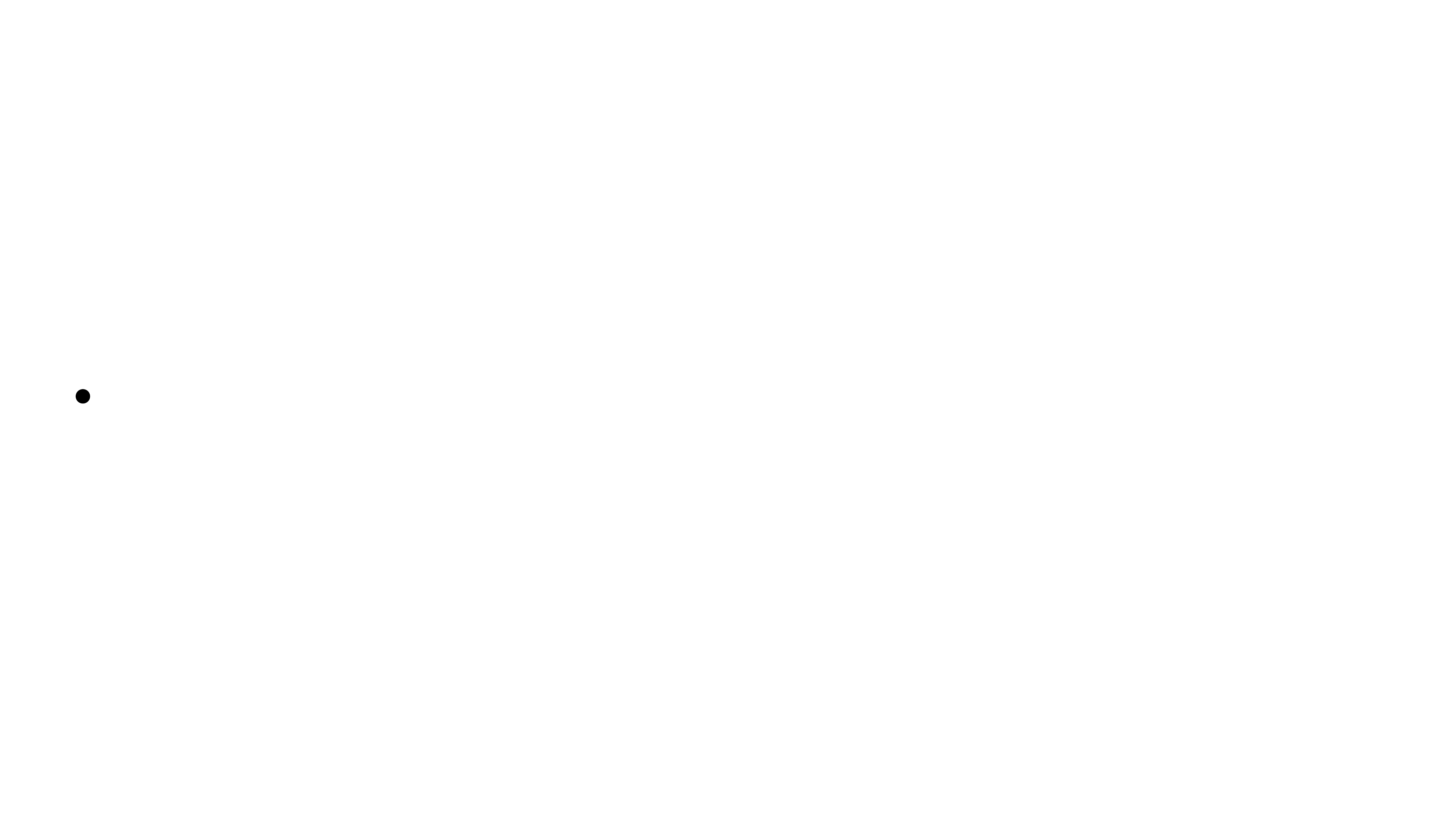
    \label{fig:sub-first}
    \end{subfigure}
\ContinuedFloat
    \begin{subfigure}[c]{\textwidth}
        \centering
        \scriptsize
        \def\svgwidth{380pt}
        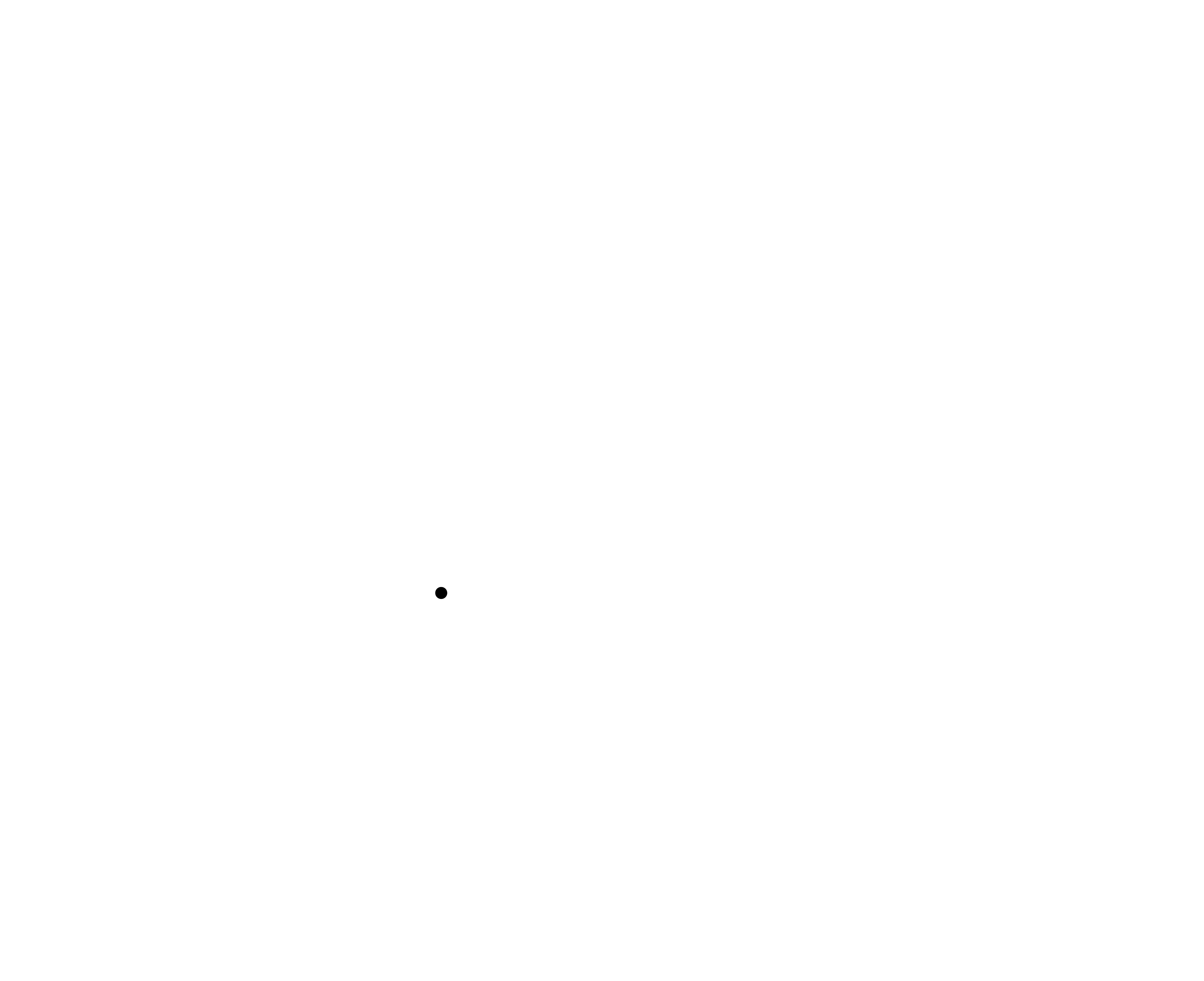
        \label{fig:sub-second}
    \end{subfigure}
\caption{The RACGs on the graphs $\Gamma_1$ and $\Gamma_1^\prime$ with their respective JSJ graphs of cylinders $\Lambda_{c,1}$ and $\Lambda^\prime_{c,1}$ are QI to each other by Theorem \ref{finalresult}. By Lemma \ref{LemCommens}, they are not commensurable.}
\label{fig:ExQInonRig}
\end{figure}

\pagebreak

\begin{figure}[H]
    \begin{subfigure}[c]{\textwidth}
        \centering
        \scriptsize
        \def\svgwidth{455pt}
        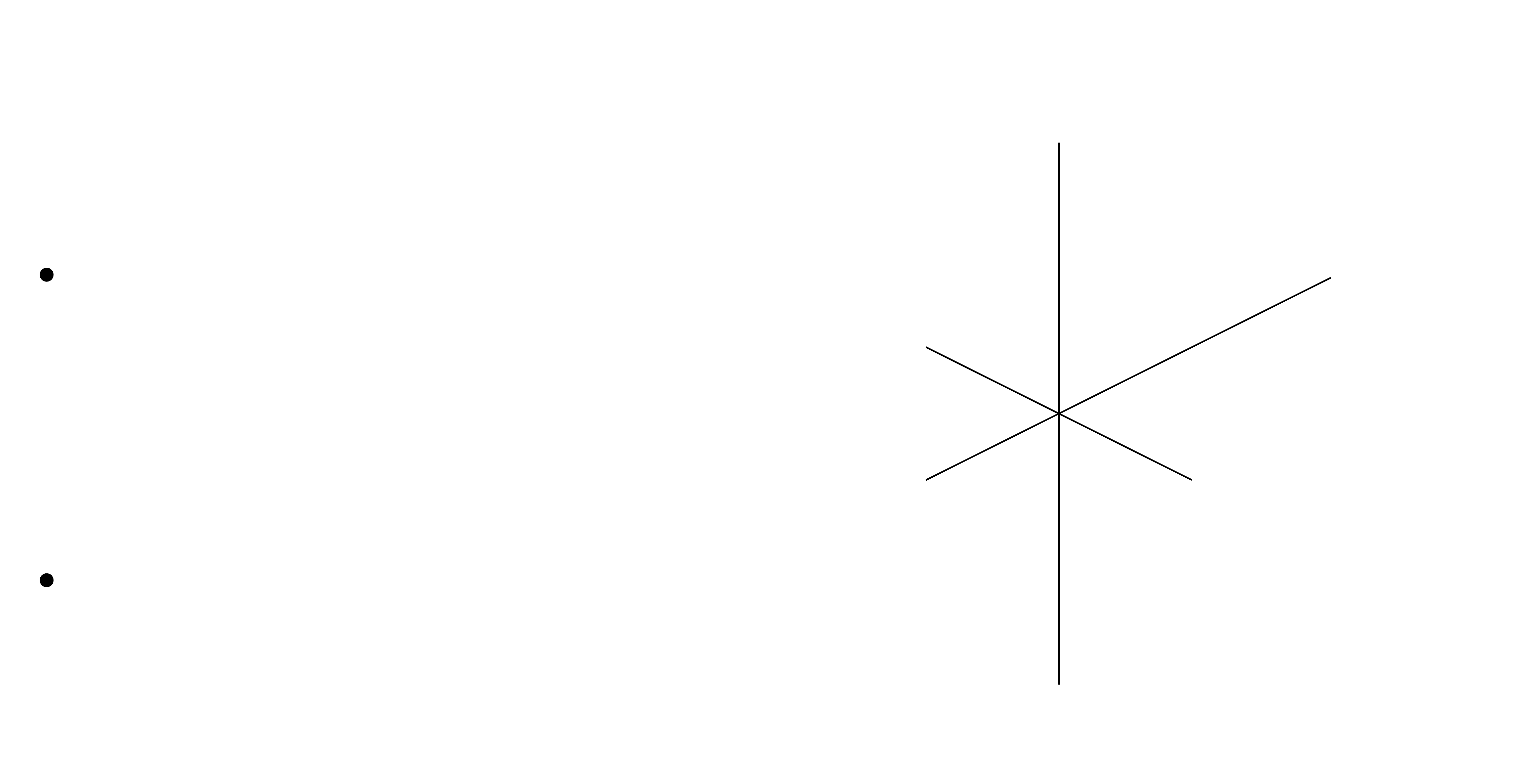
        \label{fig:ExQIwithRig1}
    \end{subfigure}\setlength{\parskip}{0pt} \setlength{\itemsep}{0pt plus 1pt}
\ContinuedFloat    
    \begin{subfigure}[c]{\textwidth}
        \centering
        \scriptsize
        \def\svgwidth{455pt}
        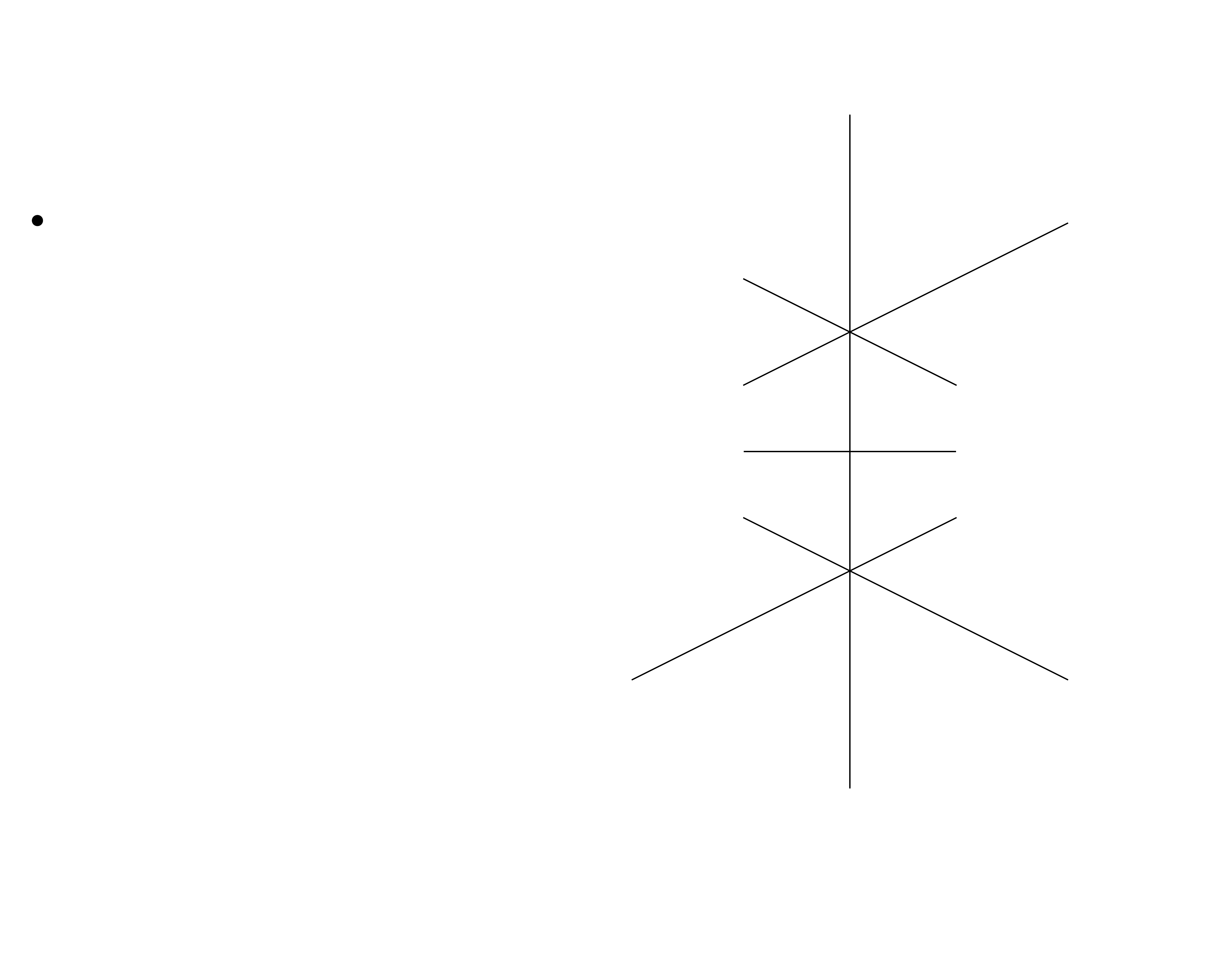
        \label{fig:ExQIwithRig2}
    \end{subfigure}
\caption{The RACGs on the graphs $\Gamma_2$ and $\Gamma_2^\prime$ with their respective JSJ graphs of cylinders $\Lambda_{c,2}$ and $\Lambda^\prime_{c,2}$ are QI to each other by Theorem \ref{finalresult}. By Lemma \ref{LemCommens}, they are not commensurable.}
\label{fig:ExQIwithRig}
\end{figure}

\newpage
\bibliography{main}

\end{document}